\documentclass[reqno]{amsart}
\usepackage{amsmath,amsfonts,amssymb,xcolor,enumitem,upgreek}
\usepackage[english]{babel}
\usepackage[colorlinks=true,linkcolor=red!60!black,citecolor=green!60!black]{hyperref} 

\newtheorem{theorem}{Theorem}[section]
\newtheorem{lemma}[theorem]{Lemma}
\newtheorem{proposition}[theorem]{Proposition}
\newtheorem{corollary}[theorem]{Corollary}

\theoremstyle{definition}
\newtheorem{definition}[theorem]{Definition}
\newtheorem{example}[theorem]{Example}

\newcommand{\w}{\omega}

\numberwithin{equation}{section}

\begin{document}
\title[Scalar D-Concave Nonautonomous Equations]{Bifurcation Theory of Attractors and Minimal Sets in D-Concave Nonautonomous Scalar Ordinary Differential Equations}
\author[J. Due\~{n}as]{Jes\'{u}s Due\~{n}as}
\author[C. N\'{u}\~{n}ez]{Carmen N\'{u}\~{n}ez}
\author[R. Obaya]{Rafael Obaya}
\address{Departamento de Matem\'{a}tica Aplicada, Universidad de Va\-lladolid, Paseo Prado de la Magdalena 3-5, 47011 Valladolid, Spain. The authors are members of IMUVA: Instituto de Investigaci\'{o}n en Matem\'{a}ticas, Universidad de Valladolid.}
\email[J.~Due\~{n}as]{jesus.duenas@uva.es}
\email[C.~N\'{u}\~{n}ez]{carmen.nunez@uva.es}
\email[R.~Obaya]{rafael.obaya@uva.es}
\thanks{All the authors were supported by Ministerio de Ciencia, Innovaci\'{o}n y Universidades (Spain) under project RTI2018-096523-B-I00 and by Universidad de Valladolid under project PIP-TCESC-2020. J. Due\~{n}as was supported by Ministerio de Universidades (Spain) under programme FPU20/01627.}
\date{}
\begin{abstract}
Two one-parametric bifurcation problems for scalar nonautonomous ordinary differential equations are analyzed assuming the coercivity of the time-dependent function determining the equation and the concavity of its derivative with respect to the state variable. The skewproduct formalism leads to the analysis of the number and properties of the minimal sets and of the shape of the global attractor, whose abrupt variations determine the occurrence of local saddle-node, local transcritical and global pitchfork bifurcation points of minimal sets and of discontinuity points of the global attractor.
\end{abstract}
\renewcommand{\subjclassname}{\textup{2020} Mathematics Subject Classification}

\maketitle
\section{Introduction}
The objective of this paper is to develop a bifurcation theory for nonautonomous scalar ordinary differential equations defined by time-dependent scalar functions with concave derivative with respect to the state variable: d-concave functions (or equations), for short.
We analyze in detail different types of d-concavity and describe rigorously conditions giving rise to saddle-node, transcritical and pitchfork nonautonomous bifurcation patterns for suitable one-parametric families of scalar d-concave ODEs.

Nonautonomous bifurcation theory is a challenging problem intensively studied during the last years. The works of Alonso and Obaya \cite{alob3}, Anagnostopoulou and J\"{a}ger \cite{anagjager1}, Anagnostopoulu et al. \cite{anjk}, Braaksma et al. \cite{brbh}, Fabbri and Johnson \cite{fabri1}, Fuhrmann \cite{fuhrmann1}, Johnson and Mantellini \cite{joma}, Johnson et al. \cite{johnson1}, Kloeden \cite{kloeden1}, Langa et al. \cite{langa1}, N\'{u}\~{n}ez and Obaya \cite{nunezobaya,nuob10}, P\"{o}tzsche \cite{potz2,potz4}, Rasmussen \cite{rasmussen2,rasmussen1}, and the references therein, offer an overview of the present state of this theory, paying special attention to scalar ordinary differential equations. Since, in general, a nonautonomous differential equation does not have constant or periodic solutions, it is not clear how to identify the objects whose bifurcation should be studied. In this paper, we make use of the skewproduct formalism, which gives some natural answers to this question: we analyze the number and shape of certain compact invariant subsets of the skewproduct flow and the changes of these structures as the parameter varies.

Let $f_0\colon\mathbb{R}\times\mathbb{R}\rightarrow\mathbb{R}$ be a continuous function, and let us denote the $t$-shift of $f_0$ by $(f_0{\cdot}t)(s,x)=f_0(t+s,x)$. Standard conditions on $f$ ensure that its hull $\Omega$, given by the closure of $\{f_0{\cdot}t \colon\, t \in \mathbb{R}\}$ in a suitable topology of $C(\mathbb{R}\times\mathbb{R},\mathbb{R})$, is a compact metric space and that the map $\sigma\colon\mathbb{R}\times\Omega\rightarrow\Omega,\;(t,\omega)\mapsto\sigma(t,\omega)=\omega{\cdot}t$ defines a continuous flow.
By defining $f(\omega,x)=\omega(0,x)$, we obtain the family of equations
\begin{equation}\label{1.primera}
 x'=f(\omega{\cdot}t,x)\,,\qquad \omega\in\Omega\,,
\end{equation}
which includes the initial one: it corresponds to $\omega_0=f_0\in\Omega$. For the first results, we will not restrict ourselves to this hull situation: $\Omega$ will be any compact metric space, $\sigma\colon\Omega\times\mathbb{R}\to\mathbb{R},\;(\omega,t)\mapsto\omega{\cdot}t$ will be a continuous flow, and $f\colon\Omega\times\mathbb{R}\rightarrow\mathbb{R}$ will be a continuous function such that the partial derivative $f_x\colon\Omega\times\mathbb{R}\rightarrow\mathbb{R}$ is jointly continuous and concave in $x$ for all $\omega\in\Omega$ (d-concave). If $u(t,\omega,x_0)$ denotes the solution of \eqref{1.primera} with $u(0,\omega,x_0)=x_0$, then $\tau\colon\mathbb{R}\times\Omega\times\mathbb{R}\rightarrow\Omega\times\mathbb{R}$, $(t,\omega,x_0)\mapsto(\omega{\cdot}t,u(t,\omega,x_0))$ defines a local skewproduct semiflow on $\Omega\times\mathbb{R}$. For our most relevant results, we will asume that $\Omega$ is minimal (and hence it is the hull of any of its elements) and that the flow $\tau$ is dissipative. Our one-parametric bifurcation analysis refers to changes in the number and characteristics of the minimal sets, and in the shape of the global attractor. This type of analysis is useful for many interesting problems.

The results are organized in five sections. Section~\ref{sec:preliminaries} contains some basic notions and properties of ergodic theory and topological dynamics.
Section \ref{sec:Dconcavity} deals with the d-concavity property. As proved by Tineo \cite{tineo1}, $f$ is d-concave if and only if, for every $\omega\in\Omega$, the second order divided differences satisfy $f(\omega,[x_1,x_0,x_2])\geq f(\omega,[x_1,x_0,x_3])$ for every $x_0,x_1,x_2,x_3\in\mathbb{R}$ with $x_1<x_2<x_3$ and $x_0\not\in\{x_1,x_2,x_3\}$. But this condition does not suffice to our purposes. Given a real interval $J=[a,b]$ and $0<\epsilon<l(J)/2=(b-a)/2$, we define the standardized $\epsilon$-module of d-concavity $b_{J,\epsilon}:\Omega \rightarrow[0,\infty)$: a function which satisfies $f(\omega,[x_1,x_0,x_2])-f(\omega,[x_1,x_0,x_3])\geq b_{J,\epsilon}(\omega)$ for every $\omega\in\Omega$, $x_0,x_1,x_2,x_3\in J$, $x_2-x_1, x_3-x_2\geq\epsilon$ and $x_0 \not\in\{x_1,x_2,x_3\}$, and which is, to some extent, optimal with this condition. We introduce different levels of strict d-concavity (abbreviated as SDC) in terms of these modules.
In particular, our function $f$ is $\mathrm{(SDC)}_m$ with respect to a $\sigma$-ergodic measure $m$ on $\Omega$ if for any $J$ there exists $\rho_J>0$ such that $m(\{\omega \in \Omega\colon\,b_{J,\epsilon}(\omega)>0\})>\rho_J$ whenever $0<\epsilon<l(J)/2$,
and $f$ is (SDC) if $m(\{\omega \in \Omega\colon\,b_{J,\epsilon}(\omega)>0\})>0$ for every $m$, $J$ and $0<\epsilon<l(J)/2$.

In Section~\ref{subsec:numberofequilibia}, we analyze some ergodic and topological properties of the compact invariant sets $\mathcal{K}\subset\Omega\times\mathbb{R}$. We prove that, if $f$ is $\mathrm{(SDC)}_m$, then there exist at most three $\tau$-ergodic measures concentrated on $\mathcal{K}$ projecting onto $m$. Similarly, if $f$ is $\mathrm{(SDC)}$, then $\mathcal{K}$ contains at most three disjoint compact $\tau$-invariant sets projecting onto $\Omega$. In addition, if $f_{xx}$ is also continuous (what we assume in the rest of this introduction), then the sum of any pair of Lyapunov exponents of $\mathcal{K}$ with respect to two different ergodic measures projecting onto the same one on $\Omega$ is lower than or equal to zero. Previous results of this same type have been obtained by Tineo \cite{tineo1} for periodic differential equations under stronger conditions of d-concavity and by J\"{a}ger \cite{jager1} for quasiperiodically forced increasing maps $T\colon\mathbb{S}^1\times[a,b] \rightarrow \mathbb{S}^1\times[a,b]$ with strictly negative Schwarzian derivative. We also explain in this section the relation between their results and ours.

In Section~\ref{sec:bifurcationproblem}, we assume that the flow $(\Omega, \sigma)$ is minimal and add to the hypotheses on $f$ its coercivity, which implies the existence of a global attractor $\mathcal{A}_\lambda$ for each one of the families $x'=f(\omega{\cdot}t,x)+\lambda$, with $\lambda\in\mathbb{R}$. We describe several possibilities for the global bifurcation diagrams of this one-parametric problem. Our main tool is to determine the number of the minimal sets and their structure for each $\lambda\in\mathbb{R}$, and the relation of these properties with the shape of $\mathcal{A}_\lambda$. We prove that, if $f$ is (SDC) and there exist three minimal sets for a point $\lambda_0\in\mathbb{R}$, then $\mathcal{A}_\lambda$ contains three minimal sets, which are hyperbolic, if and only if $\lambda$ belongs to a nondegenerate interval $(\lambda_-,\lambda_+)$, and it is simply the graph of a continuous map on $\Omega$ if $\lambda\notin[\lambda_-,\lambda_+]$. In addition, the two upper (resp.~lower) minimal sets collide on a residual invariant subset of $\Omega$ when $\lambda\downarrow\lambda_-$ (resp.~$\lambda\uparrow\lambda_+$), giving rise to two local saddle-node bifurcation points of minimal sets and two points of discontinuity of $\mathcal{A}_{\lambda}$: the global bifurcation diagram is the nonautonomous analogue of that of $x'=-x^3+x+\lambda$, with the fundamental difference of the possibility of occurrence of highly complicated dynamics on the nonhyperbolic minimal set existing at the bifurcation points. This is an extension of the saddle-node bifurcation pattern studied in \cite{alob3}, \cite{nunezobaya} and \cite{nuos4} for the concave case.
In addition, we establish conditions ensuring that, if there exist exactly two minimal sets for a value $\lambda_1$, then the situation is that described above, with $\lambda_1=\lambda_-$ or $\lambda_1=\lambda_+$. We also describe the dynamical possibilities for $\mathcal A_{\lambda}$ when there is only a minimal set for all $\lambda$ (as in $x'=-x^3+\lambda$).

In Section~\ref{sec:secondoneparametric}, also under the hypothesis of minimality of $(\Omega,\sigma)$, we analyze a problem of parametric bifurcation of a recurrent solution of a given equation. First, we show how to reformulate this problem in terms of the bifurcation analysis around the minimal set $\Omega\times\{0\}$ of a family $x'=f(\omega{\cdot}t,x)+\lambda x$, with $f(\cdot,0)\equiv 0$. A first description of the posible local transcritical, local saddle-node and global pitchfork bifurcation points which may appear for this model was given in \cite{nunezobaya} in the case of uniquely ergodic flow on $\Omega$. For an (SDC) function $f$, the results of Section~\ref{sec:secondoneparametric} extend the casuistic described in \cite{nunezobaya} to the more general case of coexistence of several Lyapunov exponents on $\Omega\times\{0\}$, which may give rise to scenarios of greater dynamical complexity. Some particular cases for which the bifurcation diagram is completely characterized are described, and an interpretation of the conclusions for the initial bifurcation problem of recurrent solutions completes the paper.

\section{Preliminaries}\label{sec:preliminaries}
Subsection~\ref{subsec:skewproductflow} explains the definition of a skewproduct flow from a family of scalar nonautonomous ordinary differential equations. Subsection~\ref{subsec:equilibria} includes the basic definitions of equilibria, semiequilibria and global upper and lower solutions, and relations among them. In Subsection~\ref{subsec:compactinvariant}, we recall some particular properties of compact invariant sets and minimal sets in the scalar case. Finally, Subsection~\ref{subsec:furstenbergremark} summarizes the required notions and properties concerning exponential dichotomy, Sacker and Sell spectrum and Lyapunov exponents of a family of scalar ODEs.

The interested reader can find in \cite{arnol}, \cite{lineardissipativescalar}, \cite{chueshov3}, \cite{coppel1}, \cite{ellis1}, \cite{furstenberg1}, \cite{jager2}, \cite{johnson1}, \cite{kloeden2}, \cite{johnson2}, \cite{mane1}, \cite{nemytskii1}, \cite{almostautomorphic}, \cite{sackersell3} and \cite{walters1} all the details of the properties summarized here.
\subsection{Scalar skewproduct flow}\label{subsec:skewproductflow}
Let $\Omega$ be a compact metric space, and let $\sigma\colon\mathbb{R}\times\Omega\rightarrow\Omega$, $(t,\omega)\mapsto \sigma_t(\omega)=\omega{\cdot}t$ define a global continuous flow on $\Omega$. All our results refer to the dynamical systems generated by the solutions of nonautonomous scalar ordinary differential equations of the family
\begin{equation}\label{eq:generalequation}
x'=f(\omega{\cdot}t,x)\,,\quad\omega\in\Omega\,,
\end{equation}
where $f\colon\Omega\times\mathbb{R}\rightarrow\mathbb{R}$ is assumed to be jointly continuous and $f_x$ is supposed to exist and to be jointly continuous. We shall denote the set of this kind of functions by $C^{0,1}(\Omega\times\mathbb{R},\mathbb{R})$. Along the paper, $C^{0,2}(\Omega\times\mathbb{R},\mathbb{R})$ and successive function spaces will appear. For example, $C^{0,2}(\Omega\times\mathbb{R},\mathbb{R})$ denotes the set of functions $f\colon\Omega\times\mathbb{R}\rightarrow\mathbb{R}$ for which $f$, $f_x$ and $f_{xx}$ exist and are jointly continuous. The family \eqref{eq:generalequation} allows us to define the local continuous flow
\begin{equation}\label{eq:flujotau}
\tau\colon\mathcal{U}\subseteq \mathbb{R}\times\Omega\times\mathbb{R}\rightarrow\Omega\times\mathbb{R}\,,\quad (t,\omega,x_0)\mapsto (\omega{\cdot}t,u(t,\omega,x_0))\,,
\end{equation}
where $\mathcal{I}_{\omega,x_0}\rightarrow\mathbb{R}$, $t\mapsto u(t,\omega,x_0)$ is the maximal solution of \eqref{eq:generalequation} with initial datum $u(0,\omega,x_0)=x_0$, and $\mathcal{U}=\bigcup_{(\omega,x_0)\in\Omega\times\mathbb{R}}(\mathcal{I}_{\omega,x_0}\times\{(\omega,x_0)\})$. That is, $u$ satisfies the \emph{cocycle property} $u(t+s,\omega,x_0)=u(t,\omega{\cdot}s,u(s,\omega,x_0))$ when the right-hand term is defined. The fiber-monotonicity of the flow $\tau$ is guaranteed by the uniqueness of solutions of the ordinary differential equation: if $x_1<x_2$, then $u(t,\omega,x_1)<u(t,\omega,x_2)$ for any $t$ in the common interval of definition of both solutions. The flow $(\Omega\times\mathbb{R},\tau)$ is a type of \emph{local skewproduct flow} on $\Omega\times\mathbb{R}$ projecting onto $(\Omega,\sigma)$, which is called the \emph{base flow} of the skewproduct. The $\omega$-\emph{section} of a set $\mathcal{K}\subseteq\Omega\times\mathbb{R}$ is defined as $(\mathcal{K})_\omega=\{x\in\mathbb{R}\colon\; (\omega,x)\in\mathcal{K}\}$, and $B_\Omega(\omega,\delta)$ represents the open ball of radius $\delta$ and centered in $\omega\in\Omega$.
\subsection{Equilibria, superequilibria and subequilibria}\label{subsec:equilibria}
We shall say that a map $\beta\colon\Omega\rightarrow\mathbb{R}$ is a $\tau$-\emph{equilibrium} if $\beta(\omega{\cdot}t)=u(t,\omega,\beta(\omega))$ for all $\omega\in\Omega$ and $t\in\mathbb{R}$, a $\tau$-\emph{subequilibrium} if $\beta(\omega{\cdot}t)\leq u(t,\omega,\beta(\omega))$ for all $\omega\in\Omega$ and $t\geq 0$, a \emph{time-reversed $\tau$-subequilibrium} if this property is satisfied when $t\leq 0$ instead of $t\geq 0$, a $\tau$-\emph{superequilibrium} if $\beta(\omega{\cdot}t)\geq u(t,\omega,\beta(\omega))$ for all $\omega\in\Omega$ and all $t\geq 0$, and a \emph{time-reversed $\tau$-superequilibrium} if this property is satisfied when $t\leq0$. Note that these definitions require $u(t,\omega,\beta(\omega))$ to be defined for either the whole $\mathbb{R}$ or a half-line. The reference to the flow $\tau$ will be frequently omitted if there is no risk of confusion. Both superequilibria and subequilibria are generally called \emph{semiequilibria}. Time-reversed semiequilibria satisfy the first semiequilibria definitions for the \emph{time-reversed flow} defined by $\widetilde\sigma\colon\mathbb{R}\times\Omega\rightarrow\Omega$, $(t,\omega)\mapsto\widetilde\sigma_t(\omega)=\omega\odot t=\omega\cdot(-t)$ and
\begin{equation*}
\widetilde\tau\colon\widetilde{\mathcal{U}}\subseteq \mathbb{R}\times\Omega\times\mathbb{R}\rightarrow\Omega\times\mathbb{R}\,,\quad (t,\omega,x_0)\mapsto (\omega\odot t,\widetilde u(t,\omega,x_0))\,,
\end{equation*}
where $\widetilde u(t,\omega,x_0)=u(-t,\omega,x_0)$ can be obtained as solutions of $x'=-f(\omega\odot t,x)$, and $\widetilde{\mathcal{U}}=\{(t,\omega,x)\colon\; (-t,\omega,x)\in\mathcal{U}\}$.

Given a Borel measure $m$ defined on $\Omega$, we shall say that an equilibrium (resp. semiequilibrium) $\beta\colon\Omega\rightarrow\mathbb{R}$ is $m$-\emph{measurable} if $\beta$ is measurable with respect to the $m$-completion of the Borel $\sigma$-algebra. We shall say that an equilibrium (resp. semiequilibrium) $\beta\colon\Omega\rightarrow\mathbb{R}$ is a \emph{semicontinuous equilibrium} (resp.~\emph{semiequilibrium}) if $\beta$ is a bounded semicontinuous map. A \emph{copy of the base for the flow $\tau$} is the graph of a continuous $\tau$-equilibrium.

Let $\beta\colon\Omega\rightarrow\mathbb{R}$ be a superequilibrium (resp.~subequilibrium). Following the construction in \cite{chueshov3} and \cite{almostautomorphic}, if we define
\begin{equation*}
\beta_s\colon\Omega\rightarrow\mathbb{R}\,,\quad\omega\mapsto u(s,\omega{\cdot}(-s),\beta(\omega{\cdot}(-s)))
\end{equation*}
for $s\geq0$, then we obtain a family $\{\beta_s\}_{s\geq0}$ of superequilibria (resp.~subequilibria) which is nonincreasing (resp.~nondecreasing) as $s$ increases. If, in addition, $\beta$ is upper (resp.~lower) semicontinuous and $\{u(t,\omega,\beta(\omega))\colon\;t\geq0,\,\omega\in\Omega\}$ is bounded, then $\beta_\infty(\omega)=\lim_{s\rightarrow\infty}\beta_s(\omega)$ is an upper (resp.~lower) semicontinuous equilibrium.

Now, let $\beta\colon\Omega\rightarrow\mathbb{R}$ be a time-reversed superequilibrium (resp.~subequilibrium). Analogously, we define
\begin{equation*}
\beta_s\colon\Omega\rightarrow\mathbb{R},\quad\omega\mapsto u(s,\omega{\cdot}(-s),\beta(\omega{\cdot}(-s)))=\widetilde u(-s,\omega\odot s,\beta(\omega\odot s)))
\end{equation*}
for $s \leq 0$, and conclude from the previous properties applied to $\widetilde\tau$ that $\{\beta_s\}_{s\leq0}$ is a family of time-reversed superequilibria (resp.~subequilibria) which is nonincreasing (resp.~nondecreasing) as $s$ decreases. If, in addition, $\beta$ is upper (resp.~lower) semicontinuous and $\{u(t,\omega,\beta(\omega))\colon\;t\leq0,\,\omega\in\Omega\}$ is bounded, then $\beta_{-\infty}(\omega)=\lim_{s\rightarrow-\infty}\beta_s(\omega)$ is an upper (resp.~lower) semicontinuous equilibrium.

A superequilibrium (resp.~subequilibrium) $\beta\colon\Omega\rightarrow\mathbb{R}$ shall be said to be \emph{strong} if there exists a time $s_*>0$ such that $\beta(\omega{\cdot}s_*)>u(s_*,\omega,\beta(\omega))$ (resp.~$\beta(\omega{\cdot}s_*)<u(s_*,\omega,\beta(\omega))$) for all $\omega\in\Omega$. Recall that the flow $(\Omega,\sigma)$ is \emph{minimal} if every $\sigma$-orbit is dense in $\Omega$. In this case, if $\beta$ is a semicontinuous strong superequilibrium (resp.~subequilibrium), then Proposition~4.3 of \cite{almostautomorphic} ensures that there exist $e_0>0$ and a time $s^*>0$ such that $\beta(\omega)\geq \beta_{s^*}(\omega)+e_0$ (resp.~$\beta(\omega)\leq \beta_{s^*}(\omega)-e_0$) for every $\omega\in\Omega$. Note that the nonincreasing (resp.~nondecreasing) monotonicity of the family $\{\beta_s\}_{s\geq0}$ ensures that $\beta(\omega)\geq \beta_{s}(\omega)+e_0$ (resp.~$\beta(\omega)\leq \beta_{s}(\omega)-e_0$) for every $\omega\in\Omega$ and all $s\geq s^*$. Strong time-reversed semiequilibria have the same definition and properties with strictly negative times $s_*$ and $s^*$.

The concepts of superequilibria and subequilibria are strongly related to those of global upper and lower solutions. An $m$-measurable map $\kappa\colon\Omega\rightarrow\mathbb{R}$ will be said to be $C^1$ \emph{along the base orbits} if, for any $\omega\in\Omega$, the map $t\mapsto \kappa_\omega(t)=\kappa(\omega{\cdot}t)$ is $C^1$ on $\mathbb{R}$. In this case, we represent $\kappa'(\omega)=\kappa'_\omega(0)$. It is clear by definition that every $m$-measurable equilibrium is $C^1$ along the base orbits. Such a map $\kappa$ shall be said to be a \emph{global upper} (resp.~\emph{lower}) \emph{solution} for a family $x'=f(\omega{\cdot}t,x)$ of differential equations if $\kappa'(\omega)\geq f(\omega,\kappa(\omega))$ (resp.~$\kappa'(\omega)\leq f(\omega,\kappa(\omega))$) for every $\omega\in\Omega$, and to be \emph{strict} if the inequalities are strict for every $\omega\in\Omega$.

If every forward $\tau$-semiorbit is globally defined, that is, if $[0,\infty)\subseteq \mathcal{I}_{\omega,x_0}$ for every $(\omega,x_0)\in\Omega\times\mathbb{R}$, a comparison argument shows that $\beta\colon\Omega\rightarrow\mathbb{R}$ is $C^1$ along the base orbits and a superequilibrium (resp.~subequilibrium) if and only if it is a global upper (resp.~lower) solution. Analogously, if $(-\infty,0]\subseteq \mathcal{I}_{\omega,x_0}$ for every $(\omega,x_0)\in\Omega\times\mathbb{R}$, then $\beta\colon\Omega\rightarrow\mathbb{R}$ is $C^1$ along the base orbits and a time-reversed subequilibrium (resp.~superequilibrium) if and only if it is a global upper (resp.~lower) solution. Moreover, in the case of globally defined forward $\tau$-semiorbits, any strict global upper (resp.~lower) solution is a strong superequilibrium (resp.~subequilibrium) and, in the case of globally defined backward $\tau$-semiorbits, it is a strong time-reversed subequilibrium (resp.~superequilibrium) (see Sections 3 and 4 of \cite{almostautomorphic}).

The following two useful propositions explore the relations between different semiequilibria which are somehow connected.
\begin{proposition} \label{prop:atravesar} Let $(\Omega,\sigma)$ be minimal. Let $\beta\colon\Omega\rightarrow\mathbb{R}$ be a semicontinuous strong superequilibrium (resp.~subequilibrium), let $\kappa\colon\Omega\rightarrow\mathbb{R}$ be a semicontinuous superequilibrium (resp.~subequilibrium) and let us assume that there exists a residual set $\mathcal{R}$ of continuity points of both maps such that $\beta(\omega)=\kappa(\omega)$ for all $\omega\in\mathcal{R}$. Then, there exist $e>0$ and $s_*>0$ such that $\kappa(\omega{\cdot}s_*)-e>u(s_*,\omega,\beta(\omega))$ (resp.~$\kappa(\omega{\cdot}s_*)+e<u(s_*,\omega,\beta(\omega))$) for all $\omega\in \Omega$.
\end{proposition}
\begin{proof} If $\beta$ is a semicontinuous strong superequilibrium, then there exist $s_0>0$ and $e_0>0$ such that $u(s_0,\omega,\beta(\omega))<\beta(\omega{\cdot}s_0)-e_0$ for all $\omega\in\Omega$. Let $\omega_0\in\mathcal{R}\cap\sigma_{-s_0}(\mathcal{R})$. Then, $\omega_0$ is a continuity point of $\beta$, $\beta\circ\sigma_{s_0}$ and $\kappa\circ\sigma_{s_0}$ and, in addition, $\beta(\omega_0{\cdot}s_0)=\kappa(\omega_0{\cdot}s_0)$. Consequently, $u(s_0,\omega_0,\beta(\omega_0))<\beta(\omega_0{\cdot}s_0)-e_0=\kappa(\omega_0{\cdot}s_0)-e_0$.

By the continuity of the involved semiequilibria at $\omega_0$ and the continuous dependence of solutions on initial data, there exists $\rho>0$ such that $u(s_0,\omega,\beta(\omega))<\kappa(\omega{\cdot}s_0)-e_0$ for all $\omega\in B_\Omega(\omega_0,\rho)$. By minimality of $(\Omega,\sigma)$, there exists $s_1>0$ such that for all $\omega\in\Omega$ there exists $0<s_\omega\leq s_1$ such that $\omega{\cdot}s_\omega\in B_\Omega(\omega_0,\rho)$. Therefore, using the cocycle property, the fiber-monotonicity, the superequilibrium property and the previous inequality, we have, for all $\omega\in\Omega$,
\begin{equation*}
\begin{split}
u(s_0+s_\omega,\omega,\beta(\omega))&=u(s_0,\omega{\cdot}s_\omega,u(s_\omega,\omega,\beta(\omega)))\\ &\leq u(s_0,\omega{\cdot}s_\omega,\beta(\omega{\cdot}s_\omega))<\kappa(\omega{\cdot}(s_0+s_\omega))-e_0\,.
\end{split}
\end{equation*}
By fiber-monotonicity, evolving both sides $s_1-s_\omega>0$ we obtain
\begin{equation}\label{eq:4.21}
u(s_0+s_1,\omega,\beta(\omega))<u(s_1-s_\omega,\omega{\cdot}(s_0+s_\omega),\kappa(\omega{\cdot}(s_0+s_\omega))-e_0)\,.
\end{equation}
Since $u(t,\omega,x)-u(t,\omega,x-e_0)>0$ for all $(t,\omega,x)\in [0,s_1]\times \Omega\times\mathrm{cl}_\mathbb{R}(\kappa(\Omega))$, which is a compact set, there exists $e>0$ such that $u(t,\omega,x-e_0)<u(t,\omega,x)-e$ for all $(t,\omega,x)\in [0,s_1]\times \Omega\times\mathrm{cl}_\mathbb{R}(\kappa(\Omega))$. Then, \eqref{eq:4.21} yields
\begin{equation*}
u(s_0+s_1,\omega,\beta(\omega))<u(s_1-s_\omega,\omega{\cdot}(s_0+s_\omega),\kappa(\omega{\cdot}(s_0+s_\omega)))-e\leq\kappa(\omega{\cdot}(s_0+s_1))-e
\end{equation*}
for all $\omega\in\Omega$, since $\kappa$ is a also a superequilibrium. Rewriting $s_*=s_0+s_1$, we obtain $u(s_*,\omega,\beta(\omega))<\kappa(\omega{\cdot}s_*)-e$ for all $\omega\in\Omega$, as we wanted to show. The subequilibrium case is proved analogously.
\end{proof}
\begin{proposition}\label{prop:atravesar2} Let $\beta\colon[0,1]\times\Omega\rightarrow\mathbb{R}$, $(\lambda,\omega)\mapsto \beta_\lambda(\omega)$ be a continuous map such that $\beta_\lambda$ is a strong superequilibrium (resp.~subequilibrium) for every $\lambda\in[0,1]$ and $\beta_\lambda(\omega)\leq\beta_\xi(\omega)$ for all $\omega\in\Omega$ if $\lambda\leq \xi$. Then, there exist $e_0>0$ and $s_0\geq 0$ such that $u(s_0,\omega,\beta_1(\omega))\leq \beta_0(\omega{\cdot}s_0)-e_0$ (resp.~$\beta_1(\omega{\cdot}s_0)+e_0\leq u(s_0,\omega,\beta_0(\omega))$) for all $\omega\in\Omega$.
\end{proposition}
\begin{proof}
We work in the superequilibrium case. Let us define $A=\{\lambda\in[0,1]\colon$ there exist $e_\lambda>0\text{ and }s_\lambda\geq 0\text{ such that }(\beta_1)_{s_\lambda}(\omega)\leq \beta_\lambda(\omega)-e_\lambda\text{ for all } \omega\in\Omega\}$. It is nonempty since $\beta_1$ is a continuous strong superequilibrium, so let us define $\lambda_0=\inf A$. As $(\lambda,\omega)\mapsto\beta_\lambda(\omega)$ is continuous, for each $\lambda\in[0,1]$, there exists a neighborhood $V_\lambda\subseteq[0,1]$ of $\lambda$ such that $|\beta_\lambda(\omega)-\beta_\xi(\omega)|<e_\lambda/2$ for every $\xi\in V_\lambda$ and $\omega\in\Omega$. This shows that $A$ is open in $[0,1]$. Now, let us prove that $\lambda_0\in A$.

Since $\beta_{\lambda_0}$ is a continuous strong superequilibrium, there exist $e>0$ and $s_*>0$ such that $(\beta_{\lambda_0})_{s^*}(\omega{\cdot}s_*)=u(s_*,\omega,\beta_{\lambda_0}(\omega))\leq\beta_{\lambda_0}(\omega{\cdot}s_*)-e$ for all $\omega\in\Omega$. Fixed $0<e_0<e$, we deduce from the uniform continuity on a compact neighborhood of $\beta_{\lambda_0}(\Omega)$ the existence of $\delta_0>0$ such that
\begin{equation}\label{eq:2.6}
u(s_*,\omega,x)\leq\beta_{\lambda_0}(\omega{\cdot}s_*)-e_0
\end{equation}
for every $\omega\in\Omega$ and $x\in B_\mathbb{R}(\beta_{\lambda_0}(\omega),\delta_0)$. Now, let us take $\lambda_1\in A$ with $|\beta_{\lambda_1}(\omega)-\beta_{\lambda_0}(\omega)|<\delta_0$ for all $\omega\in\Omega$. Then, by definition of $A$, there exists $s_{\lambda_1}\geq 0$ such that
\begin{equation*}
u(s_{\lambda_1},\omega,\beta_1(\omega))< \beta_{\lambda_1}(\omega{\cdot} s_{\lambda_1})
\end{equation*}
for every $\omega\in\Omega$. Evolving the last inequality by monotonicity a time step $s_*>0$ and applying \eqref{eq:2.6}, we obtain
\begin{equation*}
\begin{split}
(\beta_1)_{s_*+s_{\lambda_1}}(\omega{\cdot}(s_*+s_{\lambda_1}))&=u(s_*+s_{\lambda_1},\omega,\beta_1(\omega))\\
&< u(s_*,\omega{\cdot}s_{\lambda_1},\beta_{\lambda_1}(\omega{\cdot} s_{\lambda_1}))\leq\beta_{\lambda_0}(\omega{\cdot}(s_*+s_{\lambda_1}))-e_0
\end{split}
\end{equation*}
for all $\omega\in\Omega$. This shows that $\lambda_0\in A$. Therefore, $\lambda_0=0\in A$, and the statement follows easily from here. The subequilibrium case is proved analogously.
\end{proof}
\subsection{Compact invariant sets and minimality}\label{subsec:compactinvariant}
Let us recall some properties of compact $\tau$-invariant sets for a local skewproduct flow $(\Omega\times\mathbb{R},\tau)$ over a compact base $\Omega$, as well as some properties of minimal sets when $(\Omega,\sigma)$ is a minimal flow.

A set $\mathcal{K}\subseteq\Omega\times\mathbb{R}$ is $\tau$-\emph{invariant} if it is composed by globally defined $\tau$-orbits, and it is \emph{minimal} if it is compact, $\tau$-invariant and it does not contain properly any other compact $\tau$-invariant set. A compact $\tau$-invariant set $\mathcal{K}\subset\Omega\times\mathbb{R}$ is \emph{pinched} if there exists $\omega\in\Omega$ such that the section $(\mathcal{K})_\omega$ is a singleton.

We shall say that a compact $\tau$-invariant set {\em projects onto $\Omega$} if the continuous map $\pi\colon\mathcal{K}\rightarrow\Omega$, $(\omega,x)\mapsto\omega$ is surjective. In this case, $\pi$ maps orbits onto orbits preserving its direction (i.e., it is a \emph{flow epimorphism}), and
\begin{equation*}
\mathcal{K}\subseteq \bigcup_{\omega\in\Omega} \big(\{\omega\}\times [\alpha_\mathcal{K}(\omega),\beta_\mathcal{K}(\omega)]\big)\,,
\end{equation*}
where $\alpha_\mathcal{K}(\omega)=\inf\{x\in\mathbb{R}\colon\; (\omega,x)\in\mathcal{K}\}$ and $\beta_\mathcal{K}(\omega)=\sup\{x\in\mathbb{R}\colon\; (\omega,x)\in\mathcal{K}\}$. As $\mathcal{K}$ is closed, the graphs of $\alpha_\mathcal{K}$ and $\beta_\mathcal{K}$ are contained in $\mathcal{K}$, and $\alpha_\mathcal{K}$ and $\beta_\mathcal{K}$ are respectively lower and upper semicontinuous. The flow monotonicity and the $\tau$-invariance of $\mathcal{K}$ allow us to check that, in fact, $\alpha_\mathcal{K}$ and $\beta_\mathcal{K}$ are equilibria, which will be called the \emph{lower and upper delimiter equilibria} of $\mathcal{K}$. If $(\Omega,\sigma)$ is minimal, any $\tau$-invariant compact set projects onto $\Omega$.
\begin{proposition}\label{prop:precollision} Let $(\Omega,\sigma)$ be minimal, let $\beta\colon\Omega\rightarrow\mathbb{R}$ be a semicontinuous equilibrium and let $\omega_0$ be any continuity point of $\beta$. Then,
\begin{equation}\label{eq:precollision}
\mathcal{M}=\mathrm{cl}_{\Omega\times\mathbb{R}}\{(\omega_0{\cdot}t,\beta(\omega_0{\cdot}t))\colon\; t\in\mathbb{R}\}
\end{equation}
is a minimal set, and it is independent of the choice of $\omega_0$. In addition, the sections $(\mathcal{N})_\omega$ of any $\tau$-minimal set $\mathcal{N}\subset\Omega\times\mathbb{R}$ are singletons for all the points $\omega$ in a residual $\sigma$-invariant subset of $\Omega$.
\end{proposition}
\begin{proof}
It is obvious that $\mathcal{M}$ is $\tau$-invariant, and it is compact since $\beta$ is bounded. Let us deduce from the minimality of $(\Omega,\sigma)$ that $(\mathcal{M})_\omega=\{\beta(\omega)\}$ for any continuity point $\omega$ of $\beta$. Given $x\in(\mathcal{M})_\omega$, we write $(\omega,x)=\lim_{n\rightarrow\infty}(\omega_0{\cdot}t_n,\beta(\omega_0{\cdot}t_n))$ for a suitable sequence $\{t_n\}_{n\in\mathbb{N}}$. Since $\beta$ is continuous at $\omega$, then $x=\beta(\omega)$, as asserted. In particular, if $\mathcal{N}\subseteq\mathcal{M}$ is minimal, then $(\mathcal{N})_{\omega_0}=\{\beta(\omega_0)\}$ and hence $\mathcal{M}\subseteq\mathcal{N}$, which shows the minimality. The independence of the choice of $\omega_0$ follows from $\mathrm{cl}_{\Omega\times\mathbb{R}}\{(\omega{\cdot}t,\beta(\omega{\cdot}t))\colon\; t\in\mathbb{R}\}\subseteq\mathcal{M}$ for every continuity point $\omega$ of $\beta$. Finally, the last assertion is deduced by applying the previous properties to one of the delimiter equilibria of any minimal set $\mathcal{N}$.
\end{proof}
\begin{proposition} \label{prop:collision} Let $(\Omega,\sigma)$ be minimal, let $\beta_1,\beta_2\colon\Omega\rightarrow\mathbb{R}$ be semicontinuous equilibria, let $\mathcal{R}_1$, $\mathcal{R}_2$ be $\sigma$-invariant residual sets of continuity points of $\beta_1$ and $\beta_2$ respectively, and let $\mathcal{R}=\mathcal{R}_1\cap\mathcal{R}_2$. The following statements are equivalent:
\begin{enumerate}[label=\rm{(\alph*)}]
\item $\beta_1(\omega_0)=\beta_2(\omega_0)$ for a point $\omega_0\in\mathcal{R}$,
\item $\beta_1(\omega_0)=\beta_2(\omega_0)$ for every $\omega_0\in\mathcal{R}$,
\item there exist $\omega_1 \in\mathcal{R}_1$ and $\omega_2\in\mathcal{R}_2$ such that
\begin{equation*}
\mathrm{cl}_{\Omega\times\mathbb{R}}\{(\omega_1{\cdot}t,\beta_1(\omega_1{\cdot}t))\colon\; t\in\mathbb{R}\}= \mathrm{cl}_{\Omega\times\mathbb{R}}\{(\omega_2{\cdot}t,\beta_2(\omega_2{\cdot}t))\colon\; t\in\mathbb{R}\}\,.
\end{equation*}
\end{enumerate}
\end{proposition}
\begin{proof} (b)$\Rightarrow$(a) and (a)$\Rightarrow$(c) are immediate. Let us prove (c)$\Rightarrow$(b). Proposition~\ref{prop:precollision} proves that $\mathrm{cl}_{\Omega\times\mathbb{R}}\{(\omega_1{\cdot}t,\beta_1(\omega_1{\cdot}t))\colon\; t\in\mathbb{R}\}$ and $\mathrm{cl}_{\Omega\times\mathbb{R}}\{(\omega_2{\cdot}t,\beta_2(\omega_2{\cdot}t))\colon\; t\in\mathbb{R}\}$ are minimal sets, and they independent of the choice of $\omega_1\in\mathcal{R}_1$ and $\omega_2\in\mathcal{R}_2$. Hence, if $\omega_0\in\mathcal{R}$, then $\mathrm{cl}_{\Omega\times\mathbb{R}}\{(\omega_0{\cdot}t,\beta_1(\omega_0{\cdot}t))\colon\; t\in\mathbb{R}\}=\mathrm{cl}_{\Omega\times\mathbb{R}}\{(\omega_0{\cdot}t,\beta_2(\omega_0{\cdot}t))\colon\; t\in\mathbb{R}\}$. In particular, there exists a sequence $\{t_n\}_{n\in\mathbb{N}}$ such that $(\omega_0{\cdot}t_n,\beta_1(\omega_0{\cdot}t_n))\rightarrow(\omega_0,\beta_2(\omega_0))$. As $\omega_0\in\mathcal{R}$, we have $\beta_1(\omega_0)=\beta_2(\omega_0)$, as asserted.
\end{proof}
\begin{proposition} \label{prop:postcollision} Let $(\Omega,\sigma)$ be minimal and let $\beta_1,\beta_2\colon\Omega\rightarrow\mathbb{R}$ be, respectively, lower and upper semicontinuous equilibria such that $\beta_1(\omega)\leq\beta_2(\omega)$ for every $\omega\in\Omega$. If there exists $\omega_0\in\Omega$ such that $\beta_1(\omega_0)=\beta_2(\omega_0)$, then $\omega_0$ is a continuity point of both maps and hence the equivalences of Proposition~$\mathrm{\ref{prop:collision}}$ hold.
\end{proposition}
\begin{proof} If $\omega_n\rightarrow\omega_0$ as $n\rightarrow\infty$, then \begin{equation*}
\beta_1(\omega_0)\leq\liminf_{n\rightarrow\infty}\beta_1(\omega_n)\leq\limsup_{n\rightarrow\infty}\beta_1(\omega_n)\leq\limsup_{n\rightarrow\infty}\beta_2(\omega_n)\leq\beta_2(\omega_0)=\beta_1(\omega_0)\,.
\end{equation*}
The third term can be replaced by $\liminf_{n\rightarrow\infty} \beta_2(\omega_n)$. This shows the assertion.
\end{proof}

We shall say that two disjoint compact $\tau$-invariant sets $\mathcal{K}_1,\mathcal{K}_2\subset\Omega\times\mathbb{R}$ which project onto $\Omega$ are \emph{fiber-ordered}, $\mathcal{K}_1<\mathcal{K}_2$, if $x<y$ for every $(\omega,x)\in\mathcal{K}_1$ and $(\omega,y)\in\mathcal{K}_2$. When $(\Omega,\sigma)$ is minimal, two different $\tau$-minimal sets are always fiber-ordered: see Section 2.4 of \cite{lineardissipativescalar}. Finally, we shall say that a $\tau$-minimal set $\mathcal M\subset\Omega\times\mathbb R$ is \emph{hyperbolic attractive} (resp.~\emph{repulsive}) if it is uniformly exponentially asymptotically stable at $\infty$ (resp.~at $-\infty$). If none of these situations holds, then $\mathcal M$ is \emph{nonhyperbolic}.
\subsection{Ergodic measures and Lyapunov exponents}\label{subsec:furstenbergremark} Given a normalized Borel measure $m$ on $\Omega$, we shall say that it is $\sigma$-\emph{invariant} if $m(\sigma_t(\mathcal{B}))=m(\mathcal{B})$ for every $t\in\mathbb{R}$ and every Borel subset $\mathcal{B}\subseteq\Omega$, and that it is $\sigma$-\emph{ergodic} if it is $\sigma$-invariant and $m(\mathcal{B})\in\{0,1\}$ for every $\sigma$-invariant subset $\mathcal{B}\subseteq\Omega$. The nonempty sets of normalized $\sigma$-invariant and $\sigma$-ergodic Borel measures on $\Omega$ are represented by $\mathfrak{M}_\mathrm{inv}(\Omega,\sigma)$ and $\mathfrak{M}_\mathrm{erg}(\Omega,\sigma)$ respectively. The flow $(\Omega,\sigma)$ is said to be \emph{uniquely ergodic} if $\mathfrak{M}_\text{inv}(\Omega,\sigma)$ reduces to just one element $m$, in which case $m$ is ergodic; and it is said to be \emph{finitely ergodic} if $\mathfrak{M}_\mathrm{erg}(\Omega,\sigma)$ is a finite set. The support of $m\in\mathfrak{M}_\mathrm{inv}(\Omega,\sigma)$, $\mathrm{Supp}(m)$, is the complement of the largest open set with zero measure, and it is a compact invariant set. If $(\Omega,\sigma)$ is minimal, then $\text{Supp}(m)=\Omega$ for every $m\in\mathfrak{M}_\mathrm{inv}(\Omega,\sigma)$.

Let $a:\Omega\rightarrow\mathbb{R}$ be a continuous map. The \emph{Lyapunov exponent} of the family of linear differential equations $x'=a(\omega{\cdot}t)\,x$ with respect to $m\in\mathfrak{M}_\mathrm{erg}(\Omega,\sigma)$ is
\begin{equation*}
\gamma_a(\Omega,m)=\int_\Omega a(\omega)\, dm\,.
\end{equation*}
The family $x'=a(\omega{\cdot}t)\,x$ has \emph{exponential dichotomy} over $\Omega$ if there exist $k\geq1$ and $\delta>0$ such that either
\begin{equation*}
\exp\int_0^t a(\omega{\cdot}s)\; ds\leq ke^{-\delta t}\text{ whenever }\omega\in\Omega\text{ and }t\geq0
\end{equation*}
or
\begin{equation*}
\exp\int_0^t a(\omega{\cdot}s)\; ds\leq ke^{\delta t}\text{ whenever }\omega\in\Omega\text{ and }t\leq0\,.
\end{equation*}
The \emph{Sacker and Sell spectrum of} $a\colon\Omega\rightarrow\mathbb{R}$ is the set $\Sigma_a$ of $\lambda\in\mathbb{R}$ such that the family $x'=(a(\omega{\cdot}t)-\lambda)\,x$ does not have exponential dichotomy over $\Omega$. It is known (\cite{johnson2}, \cite{sackersell3}) that there exist $m^l,m^u\in\mathfrak{M}_\mathrm{erg}(\Omega,\sigma)$ such that $\Sigma_a=\left[\gamma_a(\Omega,m^l),\gamma_a(\Omega,m^u)\right]$, and that $\int_\Omega a(\omega)\,dm\in\Sigma_a$ for any $m\in \mathfrak{M}_\mathrm{inv}(\Omega,\sigma)$. 

Now, assume that the map $f$ of \eqref{eq:generalequation} belongs to $C^{0,1}(\Omega\times\mathbb{R},\mathbb{R})$. The \emph{Lyapunov exponent} of a compact $\tau$-invariant set $\mathcal{K}\subset\Omega\times\mathbb{R}$ with respect to $\nu\in\mathfrak{M}_\mathrm{erg}(\mathcal{K},\tau)$ is
\begin{equation*}
\gamma_{f_x}(\mathcal{K},\nu)=\int_\mathcal{K} f_x(\omega,x)\, d\nu\,.
\end{equation*}
We will omit the subscript $f_x$. The \emph{Sacker and Sell spectrum of $f_x$ on a compact $\tau$-invariant set $\mathcal{K}\subset\Omega\times\mathbb{R}$} is that of $f_x:\mathcal{K}\rightarrow\mathbb{R}$. 
Recall that, when the base flow is minimal, a $\tau$-minimal set $\mathcal{M}\subset\Omega\times\mathbb{R}$ is nonhyperbolic if and only if $0$ belongs to the Sacker and Sell spectrum of $f_x$ on $\mathcal{M}$. In addition, if all its Lyapunov exponents are strictly negative (resp.~positive), then $\mathcal{M}$ is an attractive (resp.~repulsive) hyperbolic copy of the base. See e.g. Proposition~2.8 of \cite{lineardissipativescalar}.

The results of Furstenberg \cite{furstenberg1} (see also Theorem~1.8.4. of \cite{arnol}) show that, given a compact $\tau$-invariant set $\mathcal{K}\subset\Omega\times\mathbb{R}$ projecting onto $\Omega$ and $\nu\in\mathfrak{M}_\mathrm{erg}(\mathcal{K},\tau)$ which projects onto $m\in\mathfrak{M}_\mathrm{erg}(\Omega,\sigma)$ (i.e., $m(A)=\nu((A\times\mathbb{R})\cap \mathcal{K})$), there exists an $m$-measurable equilibrium $\beta\colon\Omega\rightarrow\mathbb{R}$ with graph contained in $\mathcal{K}$ such that, for every continuous function $g\colon\Omega\times\mathbb{R}\rightarrow\mathbb{R}$,
\begin{equation}\label{eq:new}
\int_\mathcal{K} g(\omega,x)\; d\nu=\int_\Omega g(\omega,\beta(\omega))\; dm\,.
\end{equation}
In particular, the Lyapunov exponent of $\mathcal{K}$ for \eqref{eq:generalequation} with respect to any $\tau$-ergodic measure projecting onto $m$ is given by $\int_\Omega f_x(\omega,\beta(\omega))\, dm$ for a suitable $m$-measurable equilibrium. It is easy to check that the converse of the previous property is also true: any $m$-measurable equilibrium $\beta:\Omega\rightarrow\mathbb{R}$ with graph contained in $\mathcal{K}$ defines $\nu\in\mathfrak{M}_\text{erg}(\mathcal{K},\tau)$ projecting onto $m$ by \eqref{eq:new}. Clearly, $\beta_1$ and $\beta_2$ define the same measure if and only if they coincide $m$-a.e.
\section{D-Concavity and related properties}\label{sec:Dconcavity}
In this section, we shall define and study some properties which will assumed on the function $f$ of our family of scalar nonautonomous differential equations \eqref{eq:generalequation} in due time. All of them are related to the concavity with respect to the state variable of the function $f$.
\subsection{Divided differences and modules of d-concavity}\label{sec:Dconcavity1} A map $f\in C^{0,1}(\Omega\times\mathbb{R},\mathbb{R})$ is \emph{d-concave} ($\mathrm{(DC)}$ for short) if its derivative $f_x(\omega,\cdot)$ is concave on $\mathbb{R}$ for all $\omega\in\Omega$. As this property is not enough to achieve some of the desired results, we need to define some kind of stronger properties in terms of modules of d-concavity. To this end, we recall the definition of second order \emph{divided differences},
\begin{equation*}
f(\omega,[x_1,x_2,x_3])=\frac{f(\omega,[x_2,x_3])-f(\omega,[x_1,x_2])}{x_3-x_1}\,,
\end{equation*}
with $x_i\neq x_j$ for $i,j\in\{1,2,3\}$, where, for $i\in\{1,2\}$,
\begin{equation*}
f(\omega,[x_i,x_{i+1}])=\frac{f(\omega,x_{i+1})-f(\omega,x_i)}{x_{i+1}-x_i}\,.
\end{equation*}
\begin{lemma} \label{lem:3concave} A function $f\in C^{0,1}(\Omega\times\mathbb{R},\mathbb{R})$ is $\mathrm{(DC)}$ if and only if, for all $\omega\in\Omega$,
\begin{equation*}
f(\omega,[x_1,x_0,x_2])\geq f(\omega,[x_1,x_0,x_3])
\end{equation*}
for every $x_0,x_1,x_2,x_3\in\mathbb{R}$ satisfying $x_1<x_2<x_3$ and $x_0\neq x_i$ for $i\in\{1,2,3\}$.
\end{lemma}
\begin{proof} See Lemma~2.1 of part II of \cite{tineo1}.
\end{proof}
\begin{definition}\label{def:moduleofdconcavity} Let $f$ be a $\mathrm{(DC)}$ function and let $J\subset\mathbb{R}$ be a compact interval. We shall say that a continuous function $b^\epsilon\colon\Omega\rightarrow[0,\infty)$ is an \emph{$\epsilon$-module of d-concavity }of $f$ on $J$ if, for every $\omega\in\Omega$,
\begin{equation}\label{eq:32ref}
f(\omega,[x_1,x_0,x_2])\geq f(\omega,[x_1,x_0,x_3])+b^\epsilon(\omega)
\end{equation}
if $x_i\in J$ for $i\in\{0,1,2,3\}$, $x_0\neq x_i$ for $i\in\{1,2,3\}$ and $x_2-x_1,x_3-x_2\geq\epsilon$; and that $b^\epsilon$ is \emph{$m$-strict} for $m\in\mathfrak{M}_\mathrm{inv}(\Omega,\sigma)$ if, in addition, $m(\{\omega\in\Omega\colon\; b^\epsilon(\omega)>0\})>0$.
\end{definition}
Modules of d-concavity are only coherently defined when $0<2\epsilon\leq l(J)=\sup J-\inf J$. This restriction over the possible range of $\epsilon$ will be assumed from now on.
The linearity of divided differences ensures that the sum of $\epsilon$-modules of d-concavity on $J$ of two $\mathrm{(DC)}$ functions is an $\epsilon$-module of d-concavity on $J$ of the sum of both functions. For this reason, adding a finite number of d-concave functions ends up in a d-concave function with $m$-strict $\epsilon$-module of d-concavity on $J$ if any of the terms has an $m$-strict $\epsilon$-module of d-concavity on $J$.
\begin{proposition} \label{prop:compactminimal} Let $(\Omega,\sigma)$ be minimal. Let $\epsilon>0$ and $m\in\mathfrak{M}_\mathrm{inv}(\Omega,\sigma)$ be fixed and let $b^\epsilon\colon\Omega\rightarrow [0,\infty)$ be an $\epsilon$-module of d-concavity on a compact interval $J$ of a $\mathrm{(DC)}$ function $f$. Then, $b^\epsilon$ is $m$-strict if and only if there exists $\omega_0\in\Omega$ such that $b^\epsilon(\omega_0)>0$.
\end{proposition}
\begin{proof} If there exists $\omega_0\in\Omega$ such that $b^\epsilon(\omega_0)>0$, then $b^\epsilon(\omega)>0$ for all $\omega\in B_\Omega(\omega_0,\rho)$ for some $\rho>0$. By minimality, $\mathrm{Supp}(m)=\Omega$. So, $m(B_\Omega(\omega_0,\rho))>0$ and hence $b^\epsilon$ is $m$-strict. The converse assertion is trivial.
\end{proof}
\begin{definition} \label{def:standepsmodule} Let $J$ be a compact interval. Given a $\mathrm{(DC)}$ function $f$ and $\epsilon\in[0,l(J)]$, we shall define the \emph{standardized $\epsilon$-module of d-concavity of $f$ on $J$} as
\begin{equation*}
b_{J,\epsilon}(\omega)=\frac{\epsilon}{4\, l(J)^2}\min_{x\in J_\epsilon} \big\{ 2f_x(\omega,x)-f_x(\omega,x-\epsilon/2)-f_x(\omega,x+\epsilon/2)\}\,,
\end{equation*}
where $l(J)=\sup J-\inf J$ and $J_\epsilon=\{x\in \mathbb{R}\colon\; [x-\epsilon/2,x+\epsilon/2]\subseteq J\}$.
\end{definition}
Note that it is well defined since the minimum has a continuous argument and is taken on a compact set, that $b_{J,0}(\omega)=0$ for every $\omega\in\Omega$ and that $[0,l(J)]\times\Omega\rightarrow[0,\infty)$, $(\epsilon,\omega)\mapsto b_{J,\epsilon}(\omega)$ is continuous.
\begin{theorem} \label{th:standmodule}Let $f$ be $\mathrm{(DC)}$, let $m\in \mathfrak{M}_\mathrm{inv}(\Omega,\sigma)$, and let $J$ be a compact interval. Then,
\begin{enumerate}[label=\rm{(\roman*)}]
\item the map $[0,l(J)]\rightarrow[0,\infty)$, $\epsilon\mapsto b_{J,\epsilon}(\omega)$ is nondecreasing for all $\omega\in\Omega$.
\item The map $b_{J,\epsilon}$ is an $\epsilon$-module of d-concavity of $f$ on $J$ whenever $0<\epsilon\leq l(J)/2$, and it is $m$-strict for $\epsilon\in[\epsilon_0,l(J)/2]$ if $b_{J,\epsilon_0}$ is $m$-strict.
\item The map $\omega\mapsto(1/\epsilon)\int_0^{2\epsilon} b_{J,s}(\omega)\; ds$ is also an $\epsilon$-module of d-concavity of $f$ on $J$ and it takes the value $0$ if and only if $b_{J,2\epsilon}(\omega)=0$.
\item If $b^\epsilon\colon\Omega\rightarrow\mathbb{R}$ is any other $\epsilon$-module of d-concavity of $f$ on $J$, then
\begin{equation}\label{eq:39ineq}
b^\epsilon(\omega)\leq \frac{1}{3}\left(\frac{l(J)}{\epsilon}\right)^2 b_{J,2\epsilon}(\omega)
\end{equation}
for all $\omega\in\Omega$. In particular, if $b^\epsilon$ is $m$-strict, then $b_{J,2\epsilon}$ is $m$-strict.
\end{enumerate}
\end{theorem}
The proof requires the next result, which asserts that concavity is equivalent to a monotonicity property in the extreme points of a second order finite differences scheme.
\begin{lemma} \label{lemma:13} A function $g\colon\mathbb{R}\rightarrow\mathbb{R}$ is concave if and only if
\begin{equation}\label{eq:11}
\begin{split}
&(x_3-x_1)\,g(x_2)-(x_3-x_2)\,g(x_1)-(x_2-x_1)\,g(x_3)\\ &\qquad\qquad\leq (x_4-x_0)\,g(x_2)-(x_4-x_2)\,g(x_0)-(x_2-x_0)\,g(x_4)
\end{split}
\end{equation}
for all real $x_0\leq x_1\leq x_2\leq x_3\leq x_4$.
\end{lemma}
\begin{proof}
Firstly, concavity in $x_1\leq x_3\leq x_4$ and in $x_0\leq x_1\leq x_3$ gives
\begin{equation*}
\left\{\begin{array}{l}
(x_3-x_1)\,g(x_4)\leq (x_4-x_1)\,g(x_3)-(x_4-x_3)\,g(x_1)\,,\\[1ex]
(x_3-x_1)\,g(x_0)\leq (x_3-x_0)\,g(x_1)-(x_1-x_0)\,g(x_3)\,.
\end{array}\right.
\end{equation*}
Adding up $(x_2-x_0)$ times the first inequality plus $(x_4-x_2)$ times the second one, and dividing by $(x_3-x_1)(x_4-x_0)$, results in
\begin{equation*}
\frac{x_4-x_2}{x_4-x_0}\,g(x_0)+\frac{x_2-x_0}{x_4-x_0}\,g(x_4)\leq \frac{x_3-x_2}{x_3-x_1}\,g(x_1)+\frac{x_2-x_1}{x_3-x_1}\,g(x_3)\,.
\end{equation*}
We can write the previous inequality as
\begin{equation*}
g(x_2)-\frac{x_3-x_2}{x_3-x_1}\,g(x_1)-\frac{x_2-x_1}{x_3-x_1}\,g(x_3)\leq g(x_2)-\frac{x_4-x_2}{x_4-x_0}\,g(x_0)-\frac{x_2-x_0}{x_4-x_0}\,g(x_4)\,.
\end{equation*}
The result follows from multiplying by $x_3-x_1$ and bounding $x_3-x_1\leq x_4-x_0$, as both sides are positive by concavity. The converse assertion is obtained by writing $x_1=x_2=x_3$ in \eqref{eq:11}.
\end{proof}
\begin{proof}[Proof of Theorem~$\mathrm{\ref{th:standmodule}}$] (i)-(iii) Lemma~\ref{lemma:13} applied to $f_x(\omega,\cdot)$ ensures that, if $0<\epsilon\leq l(J)/2$, $x_2\in J_{2\epsilon}$ and $x_1\leq x_2-\epsilon\leq x_2\leq x_2+\epsilon\leq x_3$, then
\begin{equation}\label{eq:3.10th}
(x_3-x_1)f_x(\omega,x_2)-(x_3-x_2)f_x(\omega,x_1)-(x_2-x_1)f_x(\omega,x_3)\geq 2\,l(J)^2 b_{J,2\epsilon}(\omega)\,.
\end{equation}
Let $0<\epsilon\leq l(J)/2$ and $x_i\in J$ for $i\in\{0,1,2,3\}$ be fixed, with $x_0\neq x_i$ for $i\in\{1,2,3\}$, $x_2-x_1\geq\epsilon$ and $x_3-x_2\geq\epsilon$. Given $s\in(0,1]$, we replace $x_i$ in \eqref{eq:3.10th} by $x_i^s=sx_i+(1-s)x_0$ for $i\in\{1,2,3\}$ and $\epsilon$ by $s\epsilon$: since $x_2^s-x_1^s\geq s\epsilon$ and $x_3^s-x_2^s\geq s\epsilon$, we get
\begin{equation*}
\begin{split}
&(x_3-x_1)f_x(\omega,x_2^s)-(x_3-x_2)f_x(\omega,x_1^s)-(x_2-x_1)f_x(\omega,x_3^s)\\&\qquad\qquad\qquad\qquad\qquad\qquad\qquad\qquad\geq \frac{2\,l(J)^2}{s} b_{J,2s\epsilon}(\omega)\geq 2\,l(J)^2 b_{J,2s\epsilon}(\omega)\,.
\end{split}
\end{equation*}
It is easy to check that $f(\omega,[x_0,x_i])=\int_0^1 f_x(\omega,x_i^s)\,ds$ for $i\in\{1,2,3\}$ and obvious that $l(J)^2\geq (x_3-x_1)(x_2-x_1)$. So, integrating the preceding inequality yields
\begin{equation*}
\begin{split}
&(x_3-x_1)f(\omega,[x_0,x_2])-(x_3-x_2)f(\omega,[x_0,x_1])-(x_2-x_1) f(\omega,[x_0,x_3])\\ &\quad\geq(x_3-x_1)(x_2-x_1)\,2\int_0^1 \,b_{J,2s\epsilon}(\omega)\, ds=(x_3-x_1)(x_2-x_1)\,\frac{1}{\epsilon}\int_0^{2\epsilon} \,b_{J,s}(\omega)\, ds\,.
\end{split}
\end{equation*}
It follows easily that $\widetilde b^\epsilon(\omega)=(1/\epsilon)\int_0^{2\epsilon} \,b_{J,s}(\omega)\, ds$ satisfies Definition~\ref{def:moduleofdconcavity}, which proves (iii). Lemma~\ref{lemma:13} also implies that $[0,l(J)]\rightarrow[0,\infty)$, $\delta\mapsto b_{J,\delta}(\omega)$ is nondecreasing for all $\omega\in\Omega$, which proves (i) and makes easy to check that $\widetilde b^\epsilon(\omega)=0$ if and only if $b_{J,2\epsilon}(\omega)=0$.

To complete the proof of (ii), we must check that $b_{J,\epsilon}$ is also an $\epsilon$-module of d-concavity, which follows from
\begin{equation*}
2\int_0^1 \,b_{J,2s\epsilon}(\omega)\, ds\geq 2\,\int_{1/2}^1 \,b_{J,2s\epsilon}(\omega)\, ds\geq 2b_{J,\epsilon}(\omega)\int_{1/2}^1 ds=b_{J,\epsilon}(\omega)\,.
\end{equation*}
We have used $b_{J,\delta}(\omega)\geq0$ for all $\omega\in\Omega$ and the monotonicity of $\delta\mapsto b_{J,\delta}(\omega)$.

(iv) Let $x\in J_{2\epsilon}$. We call $x_1=x-\epsilon$, $x_2=x$, $x_3=x+\epsilon$, write \eqref{eq:32ref} for $b^\epsilon$ and any $x_0$, and add up the three expressions obtained by taking limits $x_0\rightarrow x_i$ ($i\in\{1,2,3\}$). We get
\begin{equation*}
2\epsilon f_x(\omega,x)-\epsilon f_x(\omega,x-\epsilon)-\epsilon f_x(\omega,x+\epsilon)\geq 6\epsilon^2\,b^\epsilon(\omega)\,,
\end{equation*}
from where \eqref{eq:39ineq} follows. The last assertion in (iv) is a consequence of \eqref{eq:39ineq}, since it implies that $\{\omega\in\Omega\colon\; b_{J,2\epsilon}(\omega)>0\}\supset\{\omega\in\Omega\colon\; b^\epsilon(\omega)>0\}$.
\end{proof}
Theorem~\ref{th:standmodule} explains the scope of Definition \ref{def:standepsmodule}: although the standardized $\epsilon$-module of d-concavity is probably not the largest $\epsilon$-module, it has the smallest set of zeros that any other $(\epsilon/2)$-module can have; and the smaller these sets are, the less restrictive the strict d-concavity conditions required by our main results.

\begin{proposition}\label{prop:nonzerob} Let $f\in C^{0,2}(\Omega\times\mathbb{R},\mathbb{R})$ be $\mathrm{(DC)}$ and let $\omega_0\in\Omega$ be fixed. Then, for any $\epsilon\in(0,l(J)]$, $b_{J,\epsilon}(\omega_0)=0$ if and only if $J$ contains a subinterval of length $\epsilon$ on which $f_{xx}(\omega_0,\cdot)$ is constant. Moreover, $f_{xx}(\omega_0,\cdot)$ is strictly decreasing on a compact interval $J$ if and only if $b_{J,\epsilon}(\omega_0)>0$ for every $\epsilon\in(0,l(J)]$.
\end{proposition}
\begin{proof} It is easy to check that
\begin{equation}\label{eq:integralbJeps}
\begin{split}
&2f_x(\omega_0,x)-f_x(\omega_0,x-\epsilon/2)-f_x(\omega_0,x+\epsilon/2)\\ &\qquad\quad=\frac{\epsilon}{2}\int_0^1\Big( f_{xx}(\omega_0,x-\epsilon(1-s)/2)-f_{xx}(\omega_0,x+\epsilon(1-s)/2)\Big)\, ds\,.
\end{split}
\end{equation}
The properties of $f$ ensure that the integrand is continuous on $s\in[0,1]$, it is never strictly negative, and it is identically zero if and only if $f_{xx}(\omega_0,x-\epsilon/2)=f_{xx}(\omega_0,x+\epsilon/2)$, as this condition is equivalent to $f_{xx}(\omega_0,\cdot)$ being constant on the interval $[x-\epsilon/2,x+\epsilon/2]$. Now, given any $\epsilon>0$, $b_{J,\epsilon}(\omega_0)=0$ if and only if there exists some $\widetilde{x}\in J_\epsilon$ such that \eqref{eq:integralbJeps} vanishes, and this occurs if and only if $f_{xx}(\omega_0,\cdot)$ is constant on $[\widetilde x-\epsilon/2,\widetilde x+\epsilon/2]$. The second statement follows immediately from the first one.
\end{proof}
\subsection{Strictly concave derivatives in measure}
Now, we shall define the properties which will be part of our main assumptions in the rest of the paper.
\begin{definition}\label{def:SDCm} We shall say that $f\in C^{0,1}(\Omega\times\mathbb{R},\mathbb{R})$ is
\begin{enumerate}[label=\rm{(\roman*)}]
\item \emph{strictly d-concave with respect to $m\in\mathfrak{M}_\mathrm{erg}(\Omega,\sigma)$} on a compact interval $J$ ($(\mathrm{SDC})_m$ on $J$) if it is $\mathrm{(DC)}$ and there exists $\rho>0$ such that $m(\{\omega\in\Omega\colon\; b_{J,\epsilon}(\omega)>0\})>\rho$ for all $\epsilon\in(0,l(J)/2]$; and \emph{strictly d-concave with respect to $m$} ($(\mathrm{SDC})_m$) if it is $(\mathrm{SDC})_m$ on $J$ for every compact interval $J\subset\mathbb{R}$;
\item \emph{strictly d-concave with respect to every measure} on a compact interval $J$ ($(\mathrm{SDC})_*$ on $J$) if it is $(\mathrm{SDC})_m$ on $J$ for every $m\in\mathfrak{M}_\mathrm{erg}(\Omega,\sigma)$; and \emph{strictly d-concave with respect to every measure} ($(\mathrm{SDC})_*$) if it is $(\mathrm{SDC})_*$ on $J$ for every compact interval $J\subset\mathbb{R}$.
\end{enumerate}
\end{definition}
In particular, the $(\mathrm{SDC})_m$ property on $J$ implies that $b_{J,\epsilon}$ is $m$-strict in some uniform way. We recall that Theorem~\ref{th:standmodule} shows that these restrictions on $f$ are weaker for the maps $b_{J,\epsilon}$ than for any other family of modules of d-concavity.
\begin{proposition}\label{prop:SDCmfxx} Let $m\in\mathfrak{M}_\mathrm{erg}(\Omega,\sigma)$, let $f\in C^{0,2}(\Omega\times\mathbb{R},\mathbb{R})$ be $\mathrm{(DC)}$ and let $J$ be a compact interval. Then, $f$ is $(\mathrm{SDC})_m$ on $J$ if and only if $m(\{\omega\in\Omega\colon\; f_{xx}(\omega,\cdot)\text{ is strictly decreasing on }J\})>0$.
\end{proposition}
\begin{proof} Necessity follows from Proposition~\ref{prop:nonzerob}. We prove sufficiency by contradiction. We define $\Omega_0=\{\omega\in\Omega\colon\; f_{xx}(\omega,\cdot)\text{ is not strictly decreasing on }J\}$ and assume $m(\Omega_0)=1$. For any $\omega\in\Omega_0$ there exists a nondegenerate interval $J_\omega\subseteq J$ on which $f_{xx}(\omega,\cdot)$ is constant. If
\begin{equation*}
\Omega_0^\epsilon=\{\omega\in\Omega_0\colon\, f_{xx}(\omega,\cdot)\text{ is constant on an interval }J_\omega\subseteq J\text{ of length }\epsilon\}\,,
\end{equation*}
then $\Omega_0=\bigcup_{\epsilon>0}\Omega_0^\epsilon$ and, if $\epsilon_1<\epsilon_2$, then $\Omega_0^{\epsilon_2}\subseteq\Omega_0^{\epsilon_1}$. Therefore, $1=m(\Omega_0)=\lim_{\epsilon\rightarrow 0} m(\Omega_0^\epsilon)$, so given any $0<\rho<1$ there exists $\epsilon_\rho>0$ such that $m(\Omega_0^{\epsilon_\rho})>1-\rho>0$. As $b_{J,\epsilon_\rho}(\omega)=0$ for every $\omega\in\Omega_0^{\epsilon_\rho}$ (see Proposition~\ref{prop:nonzerob}), we deduce that $m(\{\omega\in\Omega\colon\; b_{J,\epsilon_\rho}(\omega)>0\})\leq \rho$, which contradicts the $(\mathrm{SDC})_m$ property.
\end{proof}
\begin{definition}\label{def:SDC} We shall say that $f\in C^{0,1}(\Omega\times\mathbb{R},\mathbb{R})$ is
\begin{enumerate}[label=\rm{(\roman*)}]
\item \emph{weakly strictly d-concave} on a compact interval $J$ ($\mathrm{(sDC)}$ on $J$) if it is $\mathrm{(DC)}$ and for every $\epsilon\in(0,l(J)/2]$ there exists $m_{J,\epsilon}\in\mathfrak{M}_\mathrm{erg}(\Omega,\sigma)$ such that $b_{J,\epsilon}$ is $m_{J,\epsilon}$-strict; and \emph{weakly strictly d-concave} ($(\mathrm{sDC})$) if it is $(\mathrm{sDC})$ on $J$ for every compact interval $J$.
\item \emph{strictly d-concave} on a compact interval $J$ ($\mathrm{(SDC)}$ on $J$) if it is $\mathrm{(DC)}$ and $b_{J,\epsilon}$ is $m$-strict for every $\epsilon\in(0,l(J)/2]$ and $m\in\mathfrak{M}_\mathrm{erg}(\Omega,\sigma)$; and \emph{strictly d-concave} ($(\mathrm{SDC})$) if it is $(\mathrm{SDC})$ on $J$ for every compact interval $J$.
\end{enumerate}
\end{definition}
\begin{proposition} \label{prop:sdcimplications} The following statements hold:
\begin{enumerate}[label=\rm{(\roman*)}]
\item if there exists $m\in\mathfrak{M}_\mathrm{erg}(\Omega,\sigma)$ such that $f$ is $(\mathrm{SDC})_m$ then it is $(\mathrm{sDC})$,
\item if $f$ is $(\mathrm{SDC})$ then it is $(\mathrm{sDC})$,
\item if $f$ is $(\mathrm{SDC})_*$ then it is $(\mathrm{SDC})$, and $\mathrm{(SDC)}_m$ for all $m\in\mathfrak{M}_\mathrm{erg}(\Omega,\sigma)$.
\end{enumerate}
\end{proposition}
\begin{proof}
It follows immediately from Definitions~\ref{def:SDCm} and \ref{def:SDC}.
\end{proof}
\begin{proposition} \label{prop:SDCmadesimple} Let $J$ be a compact interval. Then,
\begin{enumerate}[label=\rm{(\roman*)}]
\item \label{apart1:SDCmadesimple} if $f\in C^{0,3}(\Omega\times\mathbb{R},\mathbb{R})$ is $\mathrm{(DC)}$ and there exists $\omega_0\in\mathrm{Supp}(m)$ for some $m\in\mathfrak{M}_\mathrm{erg}(\Omega,\sigma)$ such that $f_{xxx}(\omega_0,\cdot)$ is strictly negative on $J$, then $f$ is $(\mathrm{SDC})_m$ on $J$.
\item \label{apart2:SDCmadesimple} If $f\in C^{0,2}(\Omega\times\mathbb{R},\mathbb{R})$ is $\mathrm{(DC)}$ and there exists $\omega_0\in\bigcup_{m\in\mathfrak{M}_\mathrm{erg}(\Omega,\sigma)}\mathrm{Supp}(m)$ such that $f_{xx}(\omega_0,\cdot)$ is strictly decreasing on $J$, then $f$ is $(\mathrm{sDC})$ on $J$.
\item \label{apart3:SDCmadesimple} If $(\Omega,\sigma)$ is a minimal, $f\in C^{0,2}(\Omega\times\mathbb{R},\mathbb{R})$ is $\mathrm{(DC)}$ and there exists $\omega_0\in\Omega$ such that $f_{xx}(\omega_0,\cdot)$ is strictly decreasing on $J$, then $f$ is $(\mathrm{SDC})$ on $J$.
\end{enumerate}
\end{proposition}
\begin{proof} Recall that $m(B_\Omega(\omega_0,\rho))>0$ for all $\rho>0$ if $\omega_0\in\mathrm{Supp}(m)$. If $f\in C^{0,3}(\Omega\times\mathbb{R},\mathbb{R})$ and $f_ {xxx}(\omega_0,\cdot)$ is strictly negative on $J$, then so is $f_{xxx}(\omega,\cdot)$ for all $\omega\in B_\Omega(\omega_0,\rho)$ if $\rho>0$ is small enough. Proposition~\ref{prop:SDCmfxx} proves (i). If $f\in C^{0,2}(\Omega\times\mathbb{R},\mathbb{R})$ and $f_{xx}(\omega_0,\cdot)$ is strictly decreasing on $J$, then Proposition~\ref{prop:nonzerob} ensures that $b_{J,\epsilon}(\omega_0)>0$ for all $\epsilon>0$. So, for each $\epsilon>0$, $b_{J,\epsilon}(\omega)>0$ for all $\omega$ in an open ball $B_\Omega(\omega_0,\rho_\epsilon)$ and hence $b_{J,\epsilon}$ is $m$-strict if $\omega_0\in\mathrm{Supp}(m)$, which proves (ii). If $(\Omega,\sigma)$ is minimal, then $\Omega=\mathrm{Supp}(m)$ for all $m\in\mathfrak{M}_\mathrm{erg}(\Omega,\sigma)$ and (ii) ensures (iii).
\end{proof}
The following examples show that the subsets of $C^{0,1}(\Omega\times\mathbb{R},\mathbb{R})$ which satisfy the different described d-concavity properties are different. They also show this remarkable fact: it is possible to have some strict d-concavity properties without strictly concave derivative at any point of $\Omega$.
\begin{example} \label{example:2} Let $(\Omega,\sigma)$ be a minimal and uniquely ergodic base flow with unique measure $m$. Let $\mathcal{M}$ be a $\tau$-minimal set, and let $\alpha(\omega)=\inf (\mathcal{M})_\omega$ and $\beta(\omega)=\sup (\mathcal{M})_\omega$. Let us suppose that
\begin{enumerate}[label=\rm{(\arabic*)}]
\item $m(\{\omega\in\Omega\colon\,\alpha(\omega)\neq\beta(\omega)\})=1$,
\item there are continuous maps $\alpha_n,\beta_n\colon\Omega\rightarrow\mathbb{R}$ such that $(\alpha_n)\uparrow\alpha$ and $(\beta_n)\downarrow\beta$.
\end{enumerate}
The construction of such an example is detailed in Section 8.7 of \cite{johnson1} for a concave differential equation: the Riccati equation of the Hamiltonian system of Example 8.44. The present paper will suggest its existence in skewproduct flows coming from d-concave ordinary differential equations. Now, let us define
\begin{equation*}
f_n(\omega,x)=\left\{\begin{array}{cl}
-(x-\alpha_n(\omega))^3\,,&x\in(-\infty,\alpha_n(\omega))\,,\\[0.1cm]
0\,,&x\in [\alpha_n(\omega), \beta_n(\omega)]\,,\\[0.1cm]
-(x-\beta_n(\omega))^3\,,& x\in(\beta_n(\omega),\infty)\,.
\end{array}\right.
\end{equation*}
Let $J=[a,b]\supset [\alpha_1(\omega),\beta_1(\omega)]$ for all $\omega\in\Omega$. As $f(\omega,x)=\sum_{n=1}^\infty (1/2^n) f_n(\omega,x)$ and its term-by-term first and second derivatives with respect to $x$ converge uniformly on $\Omega\times J$, it defines a $C^{0,2}(\Omega\times J,\mathbb{R})$ function whose derivative $f_x(\omega,\cdot)$ is concave on $J$ for every $\omega\in\Omega$. Hence, the function
\begin{equation}\label{eq:foutsideextension}
\widetilde f(\omega,x)=\left\{\begin{array}{cl}
f(\omega,a)+f_x(\omega,a)(x-a)+\frac{f_{xx}(\omega,a)}{2}(x-a)^2-(x-a)^3,&x<a\,,\\[0.1cm]
f(\omega,x)\,,&a\leq x\leq b,\\[0.1cm]
f(\omega,b)+f_x(\omega,b)(x-b)+\frac{f_{xx}(\omega,b)}{2}(x-b)^2-(x-b)^3,&b<x\,,
\end{array}\right.
\end{equation}
is $C^{0,2}(\Omega\times\mathbb{R},\mathbb{R})$ and $\mathrm{(DC)}$. In addition, $\widetilde f_{xx}(\omega,\cdot)$ is strictly decreasing on the intervals $(-\infty,\alpha(\omega))$ and $(\beta(\omega),\infty)$ and zero on $[\alpha(\omega),\beta(\omega)]$. Let us define $C_\epsilon=\{\omega\colon\; \beta(\omega)-\alpha(\omega)\leq\epsilon\}$. It is clear that $C_{\epsilon_1}\subseteq C_{\epsilon_2}$ if $\epsilon_1<\epsilon_2$, and that $\bigcap_{\epsilon>0} C_\epsilon=\{\omega\colon\; \alpha(\omega)=\beta(\omega)\}$. The semicontinuity properties of $\alpha$ and $\beta$ ensure that $C_\epsilon$ contains an open set for all $\epsilon>0$, and hence $m(C_\epsilon)>0$ (since $(\Omega,\sigma)$ is minimal). Proposition~\ref{prop:nonzerob} ensures that $C_{\epsilon_1}\subseteq\{\omega\colon\; b_{J,\epsilon_2}(\omega)>0\}\subseteq C_{\epsilon_2}$ if $0<\epsilon_1<\epsilon_2\leq l(J)/2$. Consequently, $0<m(\{\omega\colon\; b_{J,\epsilon}(\omega)>0\})\rightarrow0$ as $\epsilon\downarrow0$, so $\widetilde f$ is $(\mathrm{sDC})=(\mathrm{SDC})$ but not $(\mathrm{SDC})_m=(\mathrm{SDC})_*$.
\end{example}
\begin{example} Let $\Omega=\raisebox{.2em}{$\mathbb{R}$}\left/\raisebox{-.2em}{$\pi\mathbb{Z}$}\right.$ and let $\sigma$ be the trivial identity flow over $\Omega$ defined by the differential equation $\omega'=0$. Let us define
\begin{equation}\label{eq:formulaF}
f(\omega,x)=\left\{\begin{array}{cl}
-(x+(\sin^2\omega)/2)^3\sin^2\omega\,,&x\in(-\infty,-(\sin^2\omega)/2)\,,\\[0.1cm]
0\,,&x\in[-(\sin^2\omega)/2,(\sin^2\omega)/2]\,,\\[0.1cm]
-(x-(\sin^2\omega)/2)^3\sin^2\omega\,,&x\in ((\sin^2\omega)/2,\infty)\,.
\end{array}\right.
\end{equation}
Then, $f$ belongs to $C^{0,2}(\Omega\times\mathbb{R},\mathbb{R})$ and it is $\mathrm{(DC)}$. Let us take $\omega\in(0,\pi]$. Since $f_{xx}(\omega,\cdot)$ is strictly decreasing when $|x|>(\sin^2\omega)/2$, Proposition~\ref{prop:nonzerob} ensures that $b_{J,\epsilon}(\omega)>0$ for any $\epsilon>0$ and any compact interval $J$ with $l(J)\geq\epsilon$ and $[-1/2,1/2]\subseteq J$ if and only if $\sin^2\omega<\epsilon$, that is, $A_\epsilon=\{\omega\in[0,\pi)\colon\, b_{J,\epsilon}(\omega)>0\}=(0,\delta(\epsilon))\cup(\pi-\delta(\epsilon),\pi)$, where $\delta(\epsilon)=\arcsin\sqrt{\epsilon}\in(-\pi/2,\pi/2)$.

The Dirac measure $m_{\omega}$ on $\omega$ belongs to $\mathfrak{M}_\mathrm{erg}(\Omega,\sigma)$ for every $\omega\in[0,\pi)$. Recall that $m_\omega(A_\epsilon)=1$ if and only if $\omega\in A_\epsilon$, and $m_\omega(A_\epsilon)=0$ if and only if $\omega\not\in A_\epsilon$. Hence, $f$ is not $(\mathrm{SDC})$; and since $\bigcap_{\epsilon>0} A_\epsilon=\emptyset$, $f$ is not $(\mathrm{SDC})_{m_\omega}$ for any $\omega\in[0,\pi)$. Nonetheless, $f$ is $\mathrm{(sDC)}$, since for each $\epsilon>0$ we can choose $\omega\in A_\epsilon$ in order to get $m_\omega(A_\epsilon)=1$. (Note that the measure $m_{J,\epsilon}$ in the definition of $(\mathrm{sDC})$ cannot be selected independently of $\epsilon$.) For the simple variation of \eqref{eq:formulaF} given by
\begin{equation*}
\widetilde f(\omega,x)=\left\{\begin{array}{cl}
-(x+(\sin^2\omega)/2)^3\,,&x\in(-\infty,-(\sin^2\omega)/2)\,,\\[0.1cm]
0\,,&x\in[-(\sin^2\omega)/2,(\sin^2\omega)/2]\,,\\[0.1cm]
-(x-(\sin^2\omega)/2)^3\,,& x\in((\sin^2\omega)/2,\infty)\,,
\end{array}\right.
\end{equation*}
we have $\{\omega\in\Omega\colon\,b_{J,\epsilon}(\omega)>0\}=[0,\delta(\epsilon))\cup(\pi-\delta(\epsilon),\pi)$, where $\delta(\epsilon)=\arcsin\sqrt{\epsilon}$. It follows that $\widetilde f$ is $(\mathrm{SDC})_m$ for $m=m_0$ but not for any $m=m_\omega$ with $\omega\in(0,\pi)$, so it is not $(\mathrm{SDC})_*$.
\end{example}
\section{Ergodic measures and compact invariant sets on \texorpdfstring{$\mathcal{K}$}{K}}\label{subsec:numberofequilibia}
Throughout this section, the function $f$ of the family of equations \eqref{eq:generalequation} will always be assumed to belong, at least, to $C^{0,1}(\Omega\times\mathbb{R},\mathbb{R})$ without further mention to this. We will study how the additional hypotheses on $f$ described in Section~\ref{sec:Dconcavity} limit the number of distinct measurable equilibria and compact invariant sets for the local skewproduct flow $\tau$ defined by \eqref{eq:flujotau}.
\subsection{Number of ergodic measures and compact invariant sets} Given a compact $\tau$-invariant set $\mathcal{K}\subset\Omega\times\mathbb{R}$, we define its \emph{variation interval} by $J_\mathcal{K}=[\inf_{(\omega,x)\in\mathcal{K}} x,\;\sup_{(\omega,x)\in\mathcal{K}} x]$. In the following results, we will exclude the cases in which $J_\mathcal{K}$ reduces to a point, that is, when $\mathcal{K}$ is the graph of a constant equilibrium. In that case, the dynamics is trivial and $b_{J_\mathcal{K},\epsilon}$ (see Section~\ref{sec:Dconcavity1}) is not defined for any $\epsilon>0$.
\begin{theorem} \label{th:numergmeas} Let $m\in\mathfrak{M}_\mathrm{erg}(\Omega,\sigma)$, let $\mathcal{K}\subset\Omega\times\mathbb{R}$ be a compact $\tau$-invariant set projecting onto $\Omega$ and assume that $f$ is $(\mathrm{SDC})_m$ on the variation interval $J_\mathcal{K}$ of $\mathcal{K}$. Then, there exist at most three distinct $\tau$-ergodic measures concentrated on $\mathcal{K}$ which project onto $m$.
\end{theorem}
\begin{proof} Let us suppose the existence of three different measures $\nu_i\in\mathfrak{M}_\mathrm{erg}(\mathcal{K},\tau)$, $i\in\{1,2,3\}$, projecting onto $m$ and let $\beta_i\colon\Omega\rightarrow\mathbb{R}$ be $m$-measurable equilibria with graph in $\mathcal{K}$ satisfying \eqref{eq:new} for $i\in\{1,2,3\}$. Clearly, $m(\{\omega\colon\;\beta_i(\omega)<\beta_j(\omega)\})\in\{0,1\}$ for $i\neq j$, as the set is $\tau$-invariant and $m$ is ergodic. Since $m(\{\omega\colon\;\beta_i(\omega)\neq\beta_j(\omega)\})>0$, there exists a $\sigma$-invariant set $\Omega_0\subseteq\Omega$ with $m(\Omega_0)=1$ such that $\beta_1(\omega)<\beta_2(\omega)<\beta_3(\omega)$ for all $\omega\in\Omega_0$ (by changing the indices if required). We will prove that the Lyapunov exponents of $\mathcal{K}$ for \eqref{eq:generalequation} with respect to $\nu_1$ and $\nu_3$ (resp.~$\nu_2$) are strictly negative (resp.~strictly positive). This fact precludes the existence of more than three different $\tau$-ergodic measures on $\mathcal{K}$.

Let $\rho>0$ be given by the $(\mathrm{SDC})_m$ property of $f$. We apply Lusin's Theorem to find a compact set $\Delta\subseteq\Omega_0$ satisfying $m(\Delta)> 1-\rho$ such that $\beta_i|_\Delta\colon\Delta\rightarrow\mathbb{R}$ are continuous for $i\in\{1,2,3\}$. Then, we can define
\begin{equation*}
\epsilon=\frac{1}{2}\inf\big\{\beta_3(\omega)-\beta_2(\omega),\, \beta_2(\omega)-\beta_1(\omega)\colon\;\omega\in\Delta\big\}>0
\end{equation*}
and observe that $2\epsilon\leq l(J_\mathcal{K})$. We also define
\begin{equation*}
b(\omega,x_0)=f(\omega,[\beta_1(\omega),x_0,\beta_2(\omega)])-f(\omega,[\beta_1(\omega),x_0,\beta_3(\omega)])
\end{equation*}
for $\omega\in\Omega_0$ and $x_0\neq\beta_i(\omega)$ if $i\in\{1,2,3\}$; and $b_i(\omega)=\lim_{x_0\rightarrow\beta_i(\omega)}b(\omega,x_0)$, which exists since $\lim_{x_0\rightarrow\widetilde x}f(\omega,[x_0,\widetilde x])=f_x(\omega,\widetilde x)$. Theorem~\ref{th:standmodule} and Definition~\ref{def:moduleofdconcavity} ensure that $b(\omega,x_0)\geq b_{J_\mathcal{K},\epsilon}(\omega)$ if $\omega\in\Delta$ and $x_0\in J_\mathcal{K}$. Thus, $b_i(\omega)\geq b_{J_\mathcal{K},\epsilon}(\omega)$ if $\omega\in\Delta$. In addition, the $\mathrm{(DC)}$ property of $f$ and Lemma~\ref{lem:3concave} ensure that $b(\omega,x_0)\geq0$ and hence $b_i(\omega)\geq 0$ for all $\omega\in\Omega_0$.

Let us prove that the Lyapunov exponents are negative for $\nu_1$ and $\nu_3$, following the ideas of the proof of Theorem~3.2 of part II of \cite{tineo1}. It is not hard to check that
\begin{equation}\label{eq:A}
\begin{split}
&b_1(\omega)-\left(\frac{f(\omega,\beta_2(\omega))-f(\omega,\beta_1(\omega))}{(\beta_2(\omega)-\beta_1(\omega))^2}-\frac{f(\omega,\beta_3(\omega))-f(\omega,\beta_1(\omega))}{(\beta_3(\omega)-\beta_1(\omega))^2}\right)\\
&\qquad =-\left(\frac{1}{\beta_2(\omega)-\beta_1(\omega)}-\frac{1}{\beta_3(\omega)-\beta_1(\omega)}\right)f_x(\omega,\beta_1(\omega))\,,
\end{split}
\end{equation}
\begin{equation}\label{eq:B}
\begin{split}
&\frac{\beta_2(\omega)-\beta_1(\omega)}{\beta_3(\omega)-\beta_2(\omega)}\,b_3(\omega)-\left(\frac{f(\omega,\beta_3(\omega))-f(\omega,\beta_2(\omega))}{(\beta_3(\omega)-\beta_2(\omega))^2}-\frac{f(\omega,\beta_3(\omega))-f(\omega,\beta_1(\omega))}{(\beta_3(\omega)-\beta_1(\omega))^2}\right)\\
&\qquad\qquad\; =-\left(\frac{1}{\beta_3(\omega)-\beta_2(\omega)}-\frac{1}{\beta_3(\omega)-\beta_1(\omega)}\right)f_x(\omega,\beta_3(\omega))\,,
\end{split}
\end{equation}
for $\omega\in\Omega_0$. For $\nu_1$, we define $v_1$ on $\Omega_0$ by $v_1=1/(\beta_2-\beta_1)-1/(\beta_3-\beta_1)$, so that $v_1>0$. Writing \eqref{eq:A} for $\omega{\cdot}t$ and using $f(\omega{\cdot}t,\beta_i(\omega{\cdot}t))=\beta_i'(\omega{\cdot}t)$, we get
\begin{equation*}
f_x\big(\omega{\cdot}t,\beta_1(\omega{\cdot}t)\big)=-\frac{v_1'(\omega{\cdot}t)}{v_1(\omega{\cdot}t)}-\frac{b_1(\omega{\cdot}t)}{v_1(\omega{\cdot}t)}
\end{equation*}
for all $\omega\in\Omega_0$ and $t\in\mathbb{R}$. In consequence, we have
\begin{equation}\label{eq:37}
\frac{1}{t}\int_0^{t} f_x\big(\omega{\cdot}s,\beta_1(\omega{\cdot}s)\big)\, ds=-\frac{1}{t}\log\left(\frac{v_1(\omega{\cdot}t)}{v_1(\omega)}\right)-\frac{1}{t}\int_0^{t}\frac{b_1(\omega{\cdot}s)}{v_1(\omega{\cdot}s)}\, ds
\end{equation}
for $t>0$ and $\omega\in\Omega_0$. Birkhoff's Ergodic Theorem (see Theorem~1 in Section 2 of Chapter 1 of \cite{sinai1} and Proposition~1.4 of \cite{johnson1}) implies the existence of a $\sigma$-invariant subset $\Omega_0^*\subseteq\Omega_0$ with $m(\Omega_0^*)=1$ such that, for every $\omega\in\Omega_0^*$,
\begin{equation*}
\begin{split}
\lim_{t\rightarrow\infty}\frac{1}{t}\int_0^{t} f_x\big(\omega{\cdot}s,\beta_1(\omega{\cdot}s)\big)\, ds&=\int_\Omega f_x\big(\omega,\beta_1(\omega)\big)\, dm\in\mathbb{R}\,,\\
\lim_{t\rightarrow\infty} \frac{1}{t}\int_0^{t}\frac{b_1(\omega{\cdot}s)}{v_1(\omega{\cdot}s)}\, ds&=\int_\Omega\frac{b_1(\omega)}{v_1(\omega)}\, dm\in [0,\infty]\,,\\
\lim_{t\rightarrow\infty}\frac{1}{t}\int_0^t\chi_\Delta(\omega{\cdot}s)\; ds&=m(\Delta)\,.
\end{split}
\end{equation*}
We fix $\omega\in\Omega_0^*$ and deduce from the last equality the existence of a sequence $\{t_n\}_{n\in\mathbb{N}}\uparrow\infty$ such that $\omega{\cdot}t_n\in\Delta$ for every $n\in\mathbb{N}$. Therefore, the sequence $\{\log(v_1(\omega{\cdot}t_n)/v_1(\omega))\}_{n\in\mathbb{N}}$ is bounded. Writing \eqref{eq:37} for $t=t_n$ and taking limit as $n\rightarrow\infty$, we get $\int_\mathcal{K} f_x(\omega,x)\; d\nu_1=-\int_\Omega (b_1(\omega)/v_1(\omega))\, dm\in(-\infty,0]$. Hence,
\begin{equation*}
\gamma(\mathcal{K},\nu_1)=\int_\mathcal{K} f_x(\omega,x)\; d\nu_1\leq -\int_\Delta \frac{b_1(\omega)}{v_1(\omega)}\, dm\leq -\int_\Delta\frac{b_{J_\mathcal{K},\epsilon}(\omega)}{v_1(\omega)}\, dm<0\,,
\end{equation*}
since $b_1/v_1\geq 0$ on $\Omega_0$, $b_1\geq b_{J_\mathcal{K},\epsilon}$ on $\Delta$, and the choice of $\rho$ guarantees that $m(\Delta\cap\{\omega\in\Omega\colon\; b_{J_\mathcal{K},\epsilon}(\omega)>0\})>0$. Analogously, $\gamma(\mathcal{K},\nu_3)<0$, taking $v_3=1/(\beta_3-\beta_2)-1/(\beta_3-\beta_1)$ and using \eqref{eq:B}, since $(\beta_2-\beta_1)/(\beta_3-\beta_2)>\epsilon/l(J_\mathcal{K})$ on $\Delta$. To prove that $\gamma(\mathcal{K},\nu_2)>0$ we use the same argument after getting the equality involving $f_x(\omega,\beta_2(\omega))$ analogous to \eqref{eq:B}, and with $v_2=1/(\beta_2-\beta_1)+1/(\beta_3-\beta_2)$.
\end{proof}
\begin{theorem}\label{th:new}
Let $\mathcal{K}\subset\Omega\times\mathbb{R}$ be a compact $\tau$-invariant set with variation interval $J_\mathcal{K}$ and projecting onto $\Omega$. Then, $\mathcal{K}$ contains at most three disjoint compact $\tau$-invariant sets projecting onto $\Omega$ at least in the following cases:
\begin{enumerate}[label=\rm{(\roman*)}]
\item if $f\in C^{0,1}(\Omega\times\mathbb{R},\mathbb{R})$ is $\mathrm{(sDC)}$ on $J_\mathcal{K}$, or
\item if $f\in C^{0,1}(\Omega\times\mathbb{R},\mathbb{R})$ is $\mathrm{(SDC)}$ on $J_\mathcal{K}$. In this case, if $\mathcal{K}$ contains three such sets and they are ordered, then they are hyperbolic copies of the base, attractive the upper and lower ones, and repulsive the middle one.
\end{enumerate}
\end{theorem}
\begin{proof}
(i) We assume the existence of such three subsets of $\mathcal{K}$, take $\epsilon>0$ smaller than the distances between them, call $m$ the measure $m_{J_\mathcal{K},\epsilon}$ provided by the $\mathrm{(sDC)}$ property of $f$, and define $\rho=m(\{\omega:\; b_{J_\mathcal{K},\epsilon}(\omega)>0\})$. Each one of the compact $\tau$-invariant subsets of $\mathcal{K}$ concentrates a $\tau$-ergodic measure projecting onto $m$, given by the expression \eqref{eq:new} corresponding to, for instance, its upper delimiter equilibrium. We repeat the proof of Theorem~\ref{th:numergmeas} applying Lusin's Theorem only to $\rho$, and conclude that the three $\tau$-ergodic measures are the unique ones projecting onto $m$. The assertion follows from here.

(ii) The first assertion follows from (i) and Proposition~\ref{prop:sdcimplications}. Now we assume the existence of three disjoint $\tau$-invariant ordered compact subsets $\mathcal{K}_1<\mathcal{K}_2<\mathcal{K}_3$ of $\mathcal{K}$, with $\beta_1<\beta_2<\beta_3$ for the corresponding upper equilibria. Let us choose $i\in\{1,2,3\}$ and $\nu_i\in\mathfrak{M}_\mathrm{erg}(\mathcal{K}_i,\tau)$, let $m\in\mathfrak{M}_\mathrm{erg}(\Omega,\sigma)$ be the projection of $\nu_i$, and define $\nu_j\in\mathfrak{M}_\mathrm{erg}(\mathcal{K}_j,\tau)$ for $j\neq i$ from $\beta_j$ by \eqref{eq:new}, so that $\nu_j$ projects onto $m$. The proof of Theorem~\ref{th:numergmeas} shows that $\gamma(\mathcal{K},\nu_j)<0<\gamma(\mathcal{K},\nu_2)$ for $j=1,3$. Since this happens independently of the choices of $i$ and $\nu_i$, we conclude that all the Lyapunov exponents of $\mathcal{K}_1$ and $\mathcal{K}_3$ (resp.~$\mathcal{K}_2$) are strictly negative (resp.~strictly positive).

Note that, since there are at most three ergodic measures on $\mathcal{K}$ projecting onto any fixed $m\in\mathfrak{M}_\mathrm{erg}(\Omega,\sigma)$, the upper and lower equilibria $\beta_i$ and $\alpha_i$ of $\mathcal{K}_i$ coincide on a set $\Omega_0$ with $m(\Omega_0)=1$ for all $m\in\mathfrak{M}_\mathrm{erg}(\Omega,\sigma)$, for $i\in\{1,2,3\}$. It follows easily that the projection of any compact $\tau$-invariant set contains points of $\Omega_0$. Let us prove that $\mathcal{K}_i$ is an attractive hyperbolic copy of the base if $i\in\{1,3\}$. We take $(\omega,\alpha_i(\omega))\in\mathcal{K}_i$, a point $(\widetilde\omega,\alpha_i(\widetilde\omega))=(\widetilde\omega,\beta_i(\widetilde\omega))$ in its $\boldsymbol\upalpha$-limit set, with $\widetilde\omega\in\Omega_0$, and a sequence $\{t_n\}_{n\in\mathbb{N}}\downarrow-\infty$ with $(\widetilde\omega,\alpha_i(\widetilde\omega))=\lim_{n\rightarrow\infty}\tau(t_n,\omega,\alpha_i(\omega))$ and such that there exists $\lim_{n\rightarrow\infty}\tau(t_n,\omega,\beta_i(\omega))$. Then, this last limit is also $(\widetilde\omega,\alpha_i(\widetilde\omega))$. This property allows us to check the assertion by repeating the proof of Proposition~2.8 of \cite{lineardissipativescalar}, which makes use of First Approximation Theorem and the strictly negative character of all the Lyapunov exponents of $\mathcal{K}_i$: just replace the points $(\omega_1,x_1)$ and $(\omega_1,x_2)$ of that proof by $(\omega,\alpha_i(\omega))$ and $(\omega,\beta_i(\omega))$. The proof is analogous for $\mathcal{K}_2$, working now with uniform exponential stability at $-\infty$ and with analogous sets.
\end{proof}
\begin{corollary} Let $(\Omega,\sigma)$ be a minimal flow, let $\mathcal{K}\subset\Omega\times\mathbb{R}$ be a compact $\tau$-invariant set with variation interval $J_\mathcal{K}$, and let $f$ be a $\mathrm{(DC)}$ function such that for every $\epsilon>0$ there exists $\omega_{J_\mathcal{K},\epsilon}\in\Omega$ with $b_{J_\mathcal{K},\epsilon}(\omega_{J_\mathcal{K},\epsilon})>0$. Then, there exist at most three disjoint compact $\tau$-invariant sets contained in $\mathcal{K}$ and, if there exist three, then they are hyperbolic copies of the base ordered as in Theorem~$\mathrm{\ref{th:new}(ii)}$.
\end{corollary}
\begin{proof}
Combine Theorem~\ref{th:new}(ii) with Proposition~\ref{prop:compactminimal}.
\end{proof}
The last result of this section shows that two different ergodic measures concentrated on a compact $\tau$-invariant set $\mathcal{K}$ provide two Lyapunov exponents of $\mathcal{K}$ with nonpositive sum. This property will be fundamental in some of the main proofs.
\begin{proposition} \label{prop:lyapexp1} Let $f\in C^{0,2}(\Omega\times\mathbb{R},\mathbb{R})$ be $\mathrm{(DC)}$, $m\in\mathfrak{M}_\mathrm{erg}(\Omega,\sigma)$, $\mathcal{K}\subset\Omega\times\mathbb{R}$ a compact $\tau$-invariant set projecting onto $\Omega$ with variation interval $J_\mathcal{K}$ and $\beta_\nu,\beta_\mu\colon\Omega\rightarrow\mathbb{R}$ two different $m$-measurable $\tau$-equilibria with graphs contained in $\mathcal{K}$. Then,
\begin{equation*}
\int_\Omega f_x(\omega,\beta_\nu(\omega))\; dm+\int_\Omega f_x(\omega,\beta_\mu(\omega))\; dm\leq 0\,.
\end{equation*}
In addition, the previous inequality is strict in any of the following cases:
\begin{enumerate}[label=\rm{(\roman*)}]
\item if $f$ is $\mathrm{(SDC)}_m$ on the variation interval $J_\mathcal{K}$ of $\mathcal{K}$,
\item if the graphs of $\beta_\nu$ and $\beta_\mu$ are contained in disjoint compact $\tau$-invariant subsets of $\mathcal{K}$ and if, in addition, $f$ is $\mathrm{(SDC)}$ or, more generally, if $m(\{\omega\in\Omega\colon\; b_{J_\mathcal{K},\epsilon}(\omega)>0\})>0$ for every $\epsilon>0$.
\end{enumerate}
\end{proposition}
\begin{proof} There is no loss of generality in assuming $\beta_\nu(\omega)<\beta_\mu(\omega)$ for $m$-a.e. $\omega\in\Omega$, and hence on a $\sigma$-invariant set $\Omega_0\subseteq\Omega$ with $m(\Omega_0)=1$. Let us define $k(\omega)=\beta_\mu(\omega)-\beta_\nu(\omega)$. Then,
\begin{equation}\label{eq:44eq1}
k'(\omega{\cdot}t)=f(\omega{\cdot}t,k(\omega{\cdot}t)+\beta_\nu(\omega{\cdot}t))-f(\omega{\cdot}t,\beta_\nu(\omega{\cdot}t))=k(\omega{\cdot}t)F(\omega{\cdot}t,k(\omega{\cdot} t))\,,
\end{equation}
where $F(\omega,y)=\int_0^1 f_x(\omega,sy+\beta_\nu(\omega))\, ds$. As $f$ is $\mathrm{(DC)}$, $F_y(\omega,\cdot)$ is a nonincreasing function for any $\omega\in\Omega$, so
\begin{equation}\label{eq:laanterior}
F(\omega,k(\omega))=F(\omega,0)+\int_0^1 k(\omega) F_y(\omega,sk(\omega))\, ds\geq F(\omega,0)+k(\omega)F_y(\omega,k(\omega))\,.
\end{equation}
Deriving the equality $yF(\omega,y)=f(\omega,y+\beta_\nu(\omega))-f(\omega,\beta_\nu(\omega))$ with respect to $y$ and evaluating at $y=k(\omega)$ yields $F(\omega,k(\omega))+k(\omega)F_y(\omega,k(\omega))=f_x(\omega,\beta_\mu(\omega))$. This equality combined with $F(\omega,0)=f_x(\omega,\beta_\nu(\omega))$ and \eqref{eq:laanterior} provides
\begin{equation}\label{eq:326}
\begin{split}
&\int_\Omega f_x(\omega,\beta_\nu(\omega))\; dm+\int_\Omega f_x(\omega,\beta_\mu(\omega))\; dm \\&\quad= \int_{\Omega} \big( F(\omega,0)+F(\omega,k(\omega))+k(\omega) F_y(\omega, k(\omega))\big)\, dm \leq 2\int_{\Omega} F(\omega,k(\omega))\, dm\,.
\end{split}
\end{equation}
According to \eqref{eq:44eq1}, $k'(\omega)/k(\omega)=F(\omega,k(\omega))$ for all $\omega\in\Omega_0$, and hence Birkhoff's Ergodic Theorem ensures that the right-hand side of \eqref{eq:326} is zero. This proves the first assertion.

Now, let us check that the hypotheses of (i) or (ii) ensure that $F_y(\omega,0)>F_y(\omega,k(\omega))$ on a subset of $\Omega$ with positive $m$-measure. Hence, the inequality \eqref{eq:laanterior} is strict on that subset, and therefore \eqref{eq:326} is also strict, which completes the proof. To this end, it is enough to show that $f_{xx}(\omega,\beta_\nu(\omega))>f_{xx}(\omega,\beta_\mu(\omega))$ on a set with positive $m$-measure, since $F_y(\omega,0)=\int_0^1 sf_{xx}(\omega,\beta_\nu(\omega))\, ds$ and $F_y(\omega,k(\omega))=\int_0^1 sf_{xx}(\omega,s\beta_\mu(\omega)+(1-s)\beta_\nu(\omega))\, ds$, both integrands are continuous in $s$, and $sf_{xx}(\omega,\beta_\nu(\omega))\geq sf_{xx}(\omega,s\beta_\mu(\omega)+(1-s)\beta_\nu(\omega))$ for all $s\in[0,1]$. Under the hypothesis of (i), the required inequality follows from Proposition~\ref{prop:SDCmfxx}. If the situation is that of (ii), then $\epsilon=\inf\{\beta_\mu(\omega)-\beta_\nu(\omega)\colon\;\omega\in\Omega\}$ is strictly positive; and, since $l([\beta_\nu(\omega),\beta_\mu(\omega)])\geq\epsilon$ for all $\omega\in\Omega$, Proposition~\ref{prop:nonzerob} ensures that $f_{xx}(\omega,\beta_\nu(\omega))>f_{xx}(\omega,\beta_\mu(\omega))$ for all the points $\omega$ of the set on which $b_{J_\mathcal{K},\epsilon}(\omega)>0$, which has positive $m$-measure.
\end{proof}
Previous results of this same type were obtained in \cite{tineo1} for periodic differential equations under stronger conditions of d-concavity and in \cite{jager1} for quasiperiodically forced increasing maps $T\colon\mathbb{S}^1\times[a,b]\rightarrow\mathbb{S}^1\times[a,b]$ with strictly negative Schwarzian derivative. In what follows, for the case of our real flow, we will show the relation between some properties of d-concavity on the function $f$ and the condition of negative Schwarzian derivative of its flow at discrete time.
\subsection{Negative Schwarzian derivative}
Let us suppose that $f\in C^{0,3}(\Omega\times\mathbb{R},\mathbb{R})$. Recall that $\mathcal{U}\supset\{0\}\times\Omega\times\mathbb{R}$ is the (open) domain of the flow $\tau$: see \eqref{eq:flujotau}. Since, for all $(\omega,x)\in\Omega\times\mathbb{R}$, $t\mapsto u_x(t,\omega,x)$ solves the variational equation of $x'=f(\omega{\cdot}t,x)$ and satisfies $u_x(0,\omega,x)=1$, we have
$u_{x}(t,\omega,x)=\exp\big(\int_0^t f_x(\omega{\cdot}s,u(s,\omega,x))\,ds\big)$ for all $(t,\omega,x)\in\mathcal{U}$. It follows easily that $u_{xx}(0,\omega,x)=0$, $u_{xxx}(0,\omega,x)=0$ and $u_{xxxt}(0,\omega,x)=f_{xxx}(\omega,x)$ for every $(\omega,x)\in\Omega\times\mathbb{R}$. Since $u_{x}(t,\omega,x)$ never vanishes, the \emph{Schwarzian derivative} with respect to the space variable, given by
\begin{equation*}
S_x u(t,\omega,x)=\frac{u_{xxx}(t,\omega,x)}{u_x(t,\omega,x)}-\frac{3}{2}\left(\frac{u_{xx}(t,\omega,x)}{u_x(t,\omega,x)}\right)^2\,,
\end{equation*}
is well defined on $\mathcal{U}$.
\begin{proposition} \label{prop:sch1}
Let $f\in C^{0,3}(\Omega\times\mathbb{R},\mathbb{R})$. Then,
\begin{enumerate}[label=\rm{(\roman*)}]
\item $S_x u(0,\omega,x)=0$ for all $(\omega,x)\in \Omega\times\mathbb{R}$.
\item The partial derivative of $S_x u$ with respect to $t$ exists and is continuous on $\mathcal{U}$, and $(S_x u)_t(0,\omega,x)=f_{xxx}(\omega,x)$ for every $(\omega,x)\in\Omega\times\mathbb{R}$.
\end{enumerate}
\end{proposition}
\begin{proof} Property (i) is trivial, and the existence and continuity of $(S_xu)_t$ on $\mathcal{U}$ follows from the regularity properties of all the involved functions. The last assertion is proved by straight computation and evaluation at $t=0$.
\end{proof}
\begin{proposition} Let $f\in C^{0,3}(\Omega\times\mathbb{R},\mathbb{R})$. Let us suppose that $f_{xxx}(\omega,x)<0$ for every $(\omega,x)\in\Omega\times\mathbb{R}$. Then, $S_x u(t,\omega,x)<0$ for every $(t,\omega,x)\in\mathcal{U}$ with $t>0$.
\end{proposition}
\begin{proof} Let us fix $(t_0,\omega_0,x_0)\in\mathcal{U}$ with $t_0>0$. We set $k=\sup_{t\in [0,t_0]}|u(t,\omega_0,x_0)|$ and note that $x_0\in[-k,k]$. Proposition \ref{prop:sch1}(ii) ensures that $(S_xu)_t(0,\omega,x)<0$ for every $(\omega,x)\in\Omega\times[-k,k]$. Combining the continuity of $s\mapsto (S_xu)_t(s,\omega,x)$ with open character of $\mathcal{U}$ and the compactness of $\Omega\times[-k,k]$, we find $0<\delta\leq t_0$ and $l<0$ such that $(S_x u)_t(s,\omega,x)\leq l$ for every $(s,\omega,x)\in [0,\delta]\times\Omega\times[-k,k]$. Consequently, it follows from Proposition \ref{prop:sch1}(i) that
\begin{equation}\label{eq:free1}
S_x u(s,\omega, x)=\int_0^s (S_x u)_t(r,\omega,x)\; dr\leq l s<0
\end{equation}
for every $(s,\omega,x)\in (0,\delta]\times\Omega\times[-k,k]$. Let $s_0\in(0,\delta]$ and $n_0\in\mathbb{N}$ be fixed with $s_0n_0=t_0$. Section 6 of Chapter 2 of \cite{demelo} gives the formula of Schwarzian derivative of a composition, which yields
\begin{equation}\label{eq:free2}
\begin{split}
& S_x u(n s_0,\omega_0,x_0)=S_x (u(s_0,\omega_0{\cdot}(n-1)s_0,u((n-1)s_0,\omega_0,x_0)))\\
&\quad\quad =S_x u(s_0,\omega_0{\cdot}(n-1)s_0,u((n-1)s_0,\omega_0,x_0))\cdot(u_x((n-1)s_0,\omega_0,x_0))^2\\ &
\quad\quad\quad+S_x u((n-1)s_0,\omega_0,x_0)
\end{split}
\end{equation}
for $n\in\{1,2,\dots,n_0\}$. Let us show by induction that $S_xu(ns_0,\omega_0,x_0)<0$ for every $n\in\{1,2,\dots,n_0\}$. Equality \eqref{eq:free1} shows it for $n=1$; and, since $u((n-1)s_0,\omega_0,x_0)\in[-k,k]$, \eqref{eq:free1} (resp.~the induction hypothesis) ensures that the first (resp.~second) term in the sum is strictly negative. In particular, $S_xu(t_0,\omega_0,x_0)=S_xu(n_0s_0,\omega_0,x_0)<0$, as asserted.
\end{proof}
\begin{proposition} Let $f\in C^{0,3}(\Omega\times\mathbb{R},\mathbb{R})$. Let us suppose that for all $(\omega,x)\in\Omega\times\mathbb{R}$ there exists a sequence $\{t_n\}_{n\in\mathbb{N}}$ of positive numbers with limit $0$ such that $(t_n,\omega,x)\in\mathcal{U}$ and $S_x u(t_n,\omega,x)\leq 0$ for every $n\in\mathbb{N}$. Then, $f$ is $\mathrm{(DC)}$.
\end{proposition}
\begin{proof} Let $(\omega,x)\in\Omega\times\mathbb{R}$ be fixed. It follows from Proposition \ref{prop:sch1} that $S_xu(t,\omega,x)>0$ for $t>0$ small enough if $(S_x u)_t(0,\omega,x)=f_{xxx}(\omega,x)>0$. Since our hypotheses ensure that this is not the case, $f_{xxx}(\omega,x)\le 0$ and hence $f$ is (DC).
\end{proof}
\section{A First One-Parametric Bifurcation Problem}\label{sec:bifurcationproblem}
This section deals with a parametric family of scalar ODEs
\begin{equation}\label{eq:parametricgeneraleq}
x'=f(\omega{\cdot}t,x)+\lambda\,,\quad\omega\in\Omega\,,
\end{equation}
where $f\in C^{0,1}(\Omega\times\mathbb{R},\mathbb{R})$ and $\lambda\in\mathbb{R}$. We will write $\eqref{eq:parametricgeneraleq}_\lambda$ to make reference to a particular value of the parameter. Let $\mathcal{I}^\lambda_{\omega,x_0}\rightarrow\mathbb{R}$, $t\mapsto u_\lambda(t,\omega,x_0)$ be the maximal solution of $\eqref{eq:parametricgeneraleq}_\lambda$ with $u_\lambda(0,\omega,x_0)=x_0$, and let $(\Omega\times\mathbb{R},\tau_\lambda)$ be the corresponding local skewproduct flow. We will also assume a coercivity property on $f$, which will guarantee the existence of a global attractor for all the flows $\tau_\lambda$.
\subsection{Coercivity and global attractor}\label{subsec:coerciveproperty}
A set $\mathcal{A}\subset\Omega\times\mathbb{R}$ is said to be the \emph{global attractor of the flow} $\tau$ if it is a compact $\tau$-invariant set and if it attracts every bounded set $\mathcal{C}\subset\Omega\times\mathbb{R}$; that is, if $\tau_t(\mathcal{C})$ is defined for any $t\geq 0$ and also $\lim_{t\rightarrow\infty} \text{dist}(\tau_t(\mathcal{C}),\mathcal{A})=0$, where
\begin{equation*}
\text{dist}(\mathcal{C}_1,\mathcal{C}_2)=\sup_{(\omega_1,x_1)\in\mathcal{C}_1}\left(\inf_{(\omega_2,x_2)\in\mathcal{C}_2}\big(\mathrm{dist}_{\Omega\times\mathbb{R}}((\omega_1,x_1),(\omega_2,x_2))\big)\right)
\end{equation*}
is the \emph{Hausdorff semidistance} from $\mathcal{C}_1$ to $\mathcal{C}_2$. The \emph{Hausdorff distance} between two compact sets $\mathcal{C}_1,\mathcal{C}_2\subset\Omega\times\mathbb{R}$ is defined by
\begin{equation*}
\text{dist}_\mathcal{H}(\mathcal{C}_1,\mathcal{C}_2)=\max\{\text{dist}(\mathcal{C}_1,\mathcal{C}_2),\text{dist}(\mathcal{C}_2,\mathcal{C}_1)\}\,.
\end{equation*}
 A function $f\colon\Omega\times\mathbb{R}\rightarrow\mathbb{R}$ is said to be \emph{coercive} ($\mathrm{(Co)}$ for short) if
\begin{equation*}
\lim_{|x|\rightarrow\infty}\frac{f(\omega,x)}{x}=-\infty
\end{equation*}
uniformly on $\omega\in\Omega$.
\begin{theorem} \label{th:attractorexistence} Let $f\in C^{0,1}(\Omega\times\mathbb{R},\mathbb{R})$ be $\mathrm{(Co)}$ and let $\tau$ be the flow induced by the family $\eqref{eq:parametricgeneraleq}_0$. Then,
\begin{enumerate}[label=\rm{(\roman*)}]
\item the flow $\tau$ admits global attractor
\begin{equation}\label{eq:attractorshape}
\mathcal{A}=\bigcup_{\omega\in\Omega}\big(\{\omega\}\times[\alpha_\mathcal{A}(\omega),\beta_\mathcal{A}(\omega)]\big)\,.
\end{equation}
In particular, any forward $\tau$-semiorbit is globally defined and bounded. Moreover, $\alpha_\mathcal{A}$ and $\beta_\mathcal{A}$ are, respectively, lower and upper semicontinuous $\tau$-equilibria, and can be obtained as the pullback limits
\begin{equation}\label{eq:limitstapasdef}
\begin{split}
\alpha_\mathcal{A}(\omega)&=\lim_{t\rightarrow\infty} u(t,\omega{\cdot}(-t),\rho_1)\ge\rho_1\,,\\
\beta_\mathcal{A}(\omega)&=\;\lim_{t\rightarrow\infty} \;u(t,\omega{\cdot}(-t),\rho_2)\le\rho_2\,,
\end{split}
\end{equation}
where the constants $\rho_1$ and $\rho_2$ satisfy $f(\omega,x)>0$ if $x\le\rho_1$ and
$f(\omega,x)<0$ if $x\ge\rho_2$ for all $\omega\in\Omega$.
\item $\mathcal{A}$ is the union of all the globally defined and bounded $\tau$-orbits.
\item If $\kappa\colon\Omega\rightarrow\mathbb{R}$ is a bounded global lower (resp.~upper) solution, then $\kappa(\omega)\leq \beta_\mathcal{A}(\omega)$ (resp.~$\alpha_\mathcal{A}(\omega)\leq\kappa(\omega)$) for every $\omega\in\Omega$; and if, in addition, it is strict, then $\kappa(\omega)<\beta_\mathcal{A}(\omega)$ (resp.~$\alpha_\mathcal{A}(\omega)<\kappa(\omega)$) for every $\omega\in\Omega$.
\end{enumerate}
\end{theorem}
\begin{proof}
(i) and (ii) These properties are proved by repeating the arguments leading to Theorem~16 of \cite{pullbackforwards} (see also Section~1.2 of \cite{carvalho1}). The existence of the constant $\rho_0$ is ensured by the coercivity property.

(iii) Let us work in the lower case. It is easy to deduce from the boundedness of any forward $\tau$-semiorbit that $x_0\leq\beta_\mathcal{A}(\omega)$ if and only if $u(t,\omega,x_0)$ is bounded from above as time decreases. The initial (non strict) conditions on $\kappa$ and a standard comparison argument ensures that $u(t,\omega,\kappa(\omega))\leq\kappa(\omega{\cdot}t)$ for any $t\leq 0$, and hence $u(t,\omega,\kappa(\omega))$ remains bounded from above as time decreases. Therefore, $\kappa\leq\beta_\mathcal{A}$. Now, we assume also that $\kappa$ is strict, and, for contradiction, that $\kappa(\omega_0)=\beta_\mathcal{A}(\omega_0)$ for some $\omega_0\in\Omega$. Then, $\kappa(\omega_0{\cdot}t)>u(t,\omega_0,\kappa(\omega_0))=u(t,\omega_0,\beta_\mathcal{A}(\omega_0))=\beta_\mathcal{A}(\omega_0{\cdot}t)$ if $t<0$, which is not possible. We proceed analogously in the upper case.
\end{proof}
Our next result guarantees that any compact $\tau$-invariant set which contains the graph of either $\alpha_\mathcal{A}$ or $\beta_\mathcal{A}$ has a nonpositive Lyapunov exponent corresponding to the ergodic measure defined by such equilibrium by \eqref{eq:new} for $m\in\mathfrak{M}_\mathrm{erg}(\Omega,\sigma)$.
\begin{proposition} \label{prop:jagerinmeasure} Let $f\in C^{0,2}(\Omega\times\mathbb{R},\mathbb{R})$ be $\mathrm{(Co)}$, and let $\alpha_\mathcal{A}$ and $\beta_\mathcal{A}$ be defined by \eqref{eq:attractorshape}. Given any $m\in\mathfrak{M}_\mathrm{erg}(\Omega,\sigma)$,
\begin{equation*}
\int_\Omega f_x(\omega,\alpha_\mathcal{A}(\omega))\; dm\leq 0\quad\text{and}\quad \int_\Omega f_x(\omega,\beta_\mathcal{A}(\omega))\; dm\leq 0\,.
\end{equation*}
\end{proposition}
\begin{proof}
We reason with $\beta_{\mathcal A}$, assuming for contradiction that $\int_\Omega f_x(\omega,\beta_\mathcal{A}(\omega))\; dm=\rho>0$. Birkhoff's Ergodic Theorem provides a nonempty $\sigma$-invariant subset $\Omega_0\subseteq\Omega$ and, for each $\omega\in\Omega_0$, a time $t_\omega>0$ such that $\int_0^t f_x(\omega{\cdot}s,\beta_{\mathcal A}(\omega{\cdot}s))\,ds\ge (\rho/2)\,t$ for all $t\ge t_\omega$. Let $L>0$ satisfy $|f_x(\omega,x)-f_x(\omega,\beta_{\mathcal A}(\omega))|\le L|x-\beta_{\mathcal A}(\omega)|$ for all $\omega\in\Omega$ and $x\in[\beta_{\mathcal A}(\omega),\beta_{\mathcal A}(\omega)+1]$. Let us take $k\in(0,1)$ with $Lk\le\rho/4$, fix $\omega_0\in\Omega_0$ and $x_0>\beta_{\mathcal A}(\omega_0)$, and use the definition of $\beta_{\mathcal A}$ to find $t_1>0$ such that $u(t,\omega_0,x_0)-\beta_{\mathcal A}(\omega_0{\cdot}t)\le k$ for all $t\ge t_1$. Then, for $\omega_1=\omega_0{\cdot}t_1$ (which belongs to $\Omega_0$) and $t>0$, there exists $\xi_t\in[\beta_{\mathcal A}(\omega_1),x_1]$ such that, if $x_1=u(t_1,\omega_0,x_0)$, then
\begin{equation}\label{distancia}
 u(t+t_1,\omega_0,x_0)-\beta_{\mathcal A}(\omega_0{\cdot}(t+t_1))=u_x(t,\omega_1,\xi_t)\,(x_1-\beta_{\mathcal A}(\omega_1))\,.
\end{equation}
It is easy to check that $|f_x(\omega_1{\cdot}s,u(s,\omega_1,\xi_t))-f_x(\omega_1{\cdot}s,\beta_{\mathcal A}(\omega_1{\cdot}s))|\le Lk\le \rho/4$ for all $s>0$, and to deduce that $u_x(t,\omega_1,\xi_t)=\exp\int_0^t f_x(\omega_1{\cdot}s,u(s,\omega_1,\xi_t))\,ds\ge e^{(\rho/4)t}$ if $t\ge t_{\omega_1}$. Hence, the left-hand term of \eqref{distancia} cannot converge to $0$ as $t\to\infty$, which is the sought-for contradiction. The argument is similar for $\alpha_{\mathcal A}$.
\end{proof}
Recall that any hyperbolic minimal set is a copy of the base: see Subsection~\ref{subsec:furstenbergremark}.
\begin{proposition} \label{prop:hyperbolicpositiveminimal} Let $(\Omega,\sigma)$ be minimal and let $f\in C^{0,2}(\Omega\times\mathbb{R},\mathbb{R})$ be $\mathrm{(Co)}$. Then,
\begin{enumerate}[label=\rm{(\roman*)}]
\item if $\tau$ admits a repulsive hyperbolic minimal set, then it admits at least three different minimal sets.
\item If $\tau$ admits two distinct hyperbolic minimal sets, then it admits at least three different minimal sets.
\end{enumerate}
\end{proposition}
\begin{proof}
(i) Let $\kappa:\Omega\rightarrow\mathbb{R}$ provide the repulsive hyperbolic copy of the base $\mathcal M$, and let $\mathcal{M}^u$ be the minimal set defined from the delimiter equilibrium $\beta_\mathcal{A}$ by \eqref{eq:precollision}. Assuming for contradiction that $\mathcal M\not<\mathcal M^u$ leads us to $\mathcal M=\mathcal M^u$, and hence to $\kappa(\w_0)=\beta_{\mathcal A}(\w_0)$ for a continuity point $\w_0$ of $\beta_{\mathcal A}$. So, $\lim_{t\to\infty}(u(t,\omega_0,x_0)-\kappa(\omega_0{\cdot}t))=0$ for $x_0>\kappa(\w_0)$. But this contradicts the repulsive hyperbolicity of $\mathcal M$: see e.g.~Proposition 2.8 of \cite{lineardissipativescalar}. An analogous argument shows the existence of $\mathcal{M}^l<\mathcal M$.


(ii) If there were exactly two minimal sets and they were hyperbolic attractive, then Theorem~3.4 of \cite{lineardissipativescalar} would guarantee that the global attractor is a hyperbolic copy of the base, which is a contradiction. Consequently, at least one of the two hyperbolic minimal sets is repulsive, and then (i) concludes the proof.
\end{proof}
Hereafter, we will consider the parametric problem $\eqref{eq:parametricgeneraleq}_\lambda$.
\begin{proposition} \label{prop:alphasuperequilibrium} Every global upper (resp.~lower) solution of $\eqref{eq:parametricgeneraleq}_\lambda$ is a strict global upper (resp.~lower) solution of $\eqref{eq:parametricgeneraleq}_\xi$ if $\xi<\lambda$ (resp.~$\lambda<\xi$). Particularly, any equilibrium for $\eqref{eq:parametricgeneraleq}_\lambda$ is a strong superequilibrium and a strong time-reversed subequilibrium for $\eqref{eq:parametricgeneraleq}_\xi$ if $\xi<\lambda$, as well as a strong subequilibrium and a strong time-reversed superequilibrium for $\eqref{eq:parametricgeneraleq}_\xi$ if $\lambda<\xi$.
\end{proposition}
\begin{proof} The results are easy consequences of the definitions of global upper and lower solutions and their relation with semiequilibria, described in Subsection~\ref{subsec:equilibria}.
\end{proof}
The information provided by Theorem~\ref{th:attractorexistence} plays a role on the next statement.
\begin{theorem} \label{th:attractorproperties} Let $f\in C^{0,1}(\Omega\times\mathbb{R},\mathbb{R})$ be $\mathrm{(Co)}$ and let
\begin{equation}\label{eq:parametricattractor}
\mathcal{A}_\lambda=\bigcup_{\omega\in\Omega}\big(\{\omega\}\times[\alpha_\lambda(\omega),\beta_\lambda(\omega)]\big)
\end{equation}
be the global attractor for the flow $\tau_\lambda$ given by $\eqref{eq:parametricgeneraleq}_\lambda$. Then,
\begin{enumerate}[label=\rm{(\roman*)}]
\item for every $\omega\in\Omega$, the maps $\lambda\mapsto \beta_\lambda(\omega)$ and $\lambda\mapsto \alpha_\lambda(\omega)$ are strictly increasing on $\mathbb{R}$ and they are, respectively, right- and left-continuous.
\item $\lim_{\lambda\rightarrow\pm\infty} \alpha_\lambda(\omega)=\lim_{\lambda\rightarrow\pm\infty}\beta_\lambda(\omega)=\pm\infty$ uniformly on $\Omega$.
\item If $f$ satisfies $\lim_{|x|\rightarrow\infty}f_x(\omega,x)=-\infty$ uniformly on $\Omega$, then there exists $\lambda_*>0$ such that $\mathcal{A}_\lambda$ is an attractive hyperbolic copy of the base if $|\lambda|\geq\lambda_*$.
\end{enumerate}
\end{theorem}
\begin{proof}(i) Let $\xi<\lambda$. Proposition~\ref{prop:alphasuperequilibrium} ensures that $\beta_\xi$ is a strict global lower solution of $\eqref{eq:parametricgeneraleq}_\lambda$. Hence, Theorem~\ref{th:attractorexistence}(iii) yields $\beta_\xi(\omega)<\beta_\lambda(\omega)$ for all $\omega\in\Omega$. Analogously, $\alpha_\lambda$ is a strict global upper solution for $\eqref{eq:parametricgeneraleq}_\xi$ and Theorem~\ref{th:attractorexistence}(iii) ensures $\alpha_\xi(\omega)<\alpha_\lambda(\omega)$ for all $\omega\in\Omega$.

As mentioned in the proof of Theorem~\ref{th:attractorexistence}(iii), $x_0\leq \beta_\lambda(\omega)$ if and only if $u_\lambda(t,\omega,x_0)$ is bounded from above as time decreases. Let $\omega_0\in\Omega$ be fixed. We take a decreasing sequence $\{\lambda_n\}_{n\in\mathbb{N}}$ with $\lambda_n\downarrow\lambda_0$ and call $y_0=\lim_{n\rightarrow\infty} \beta_{\lambda_n}(\omega_0)$, which (by the previous part) satisfies $y_0\geq\beta_{\lambda_0}(\omega_0)$. Then, for any $t\in\mathcal{I}_{\omega_0,y_0}^{\lambda_0}$, we have $u_{\lambda_0}(t,\omega_0,y_0)=\lim_{n\rightarrow\infty} u_{\lambda_n}(t,\omega_0,\beta_{\lambda_n}(\omega_0))=\lim_{n\rightarrow\infty} \beta_{\lambda_n}(\omega_0{\cdot}t)\leq \beta_{\lambda_1}(\omega_0{\cdot}t)$, so that $u_{\lambda_0}(t,\omega_0,y_0)$ remains bounded from above as time decreases. The previously mentioned property ensures that $y_0\leq\beta_{\lambda_0}(\omega_0)$, and hence that $\beta_{\lambda_0}(\omega_0)=\lim_{n\rightarrow\infty}\beta_{\lambda_n}(\omega_0)$. The proof is analogous for $\alpha_\lambda$.

(ii) Let us take $\rho>0$ and use the coercivity property to find
$\lambda_\rho>0$ large enough to guarantee that $f(\w,x)+\lambda_\rho>0$
whenever $x\le\rho$ and $f(\omega,x)-\lambda_\rho<0$ whenever $x\ge-\rho$.
Theorem \ref{th:attractorexistence}(i) ensures that, if $\lambda\ge\lambda_\rho$, then
$\alpha_\lambda(\w)
\ge\rho$ and $\beta_{-\lambda}(\w)
\le-\rho$ for all $\omega\in\Omega$, which imply the assertions in (ii).

(iii) The additional hypothesis on $f$ provides $\rho>0$ such that $f_x(\omega,x)<0$ if $|x|>\rho$ and $\omega\in\Omega$. By (ii), there exists $\lambda_*$ such that $\mathcal{A}_\lambda\subset\Omega\times[\rho,\infty)$ and $\mathcal{A}_{-\lambda}\subset\Omega\times(-\infty,-\rho]$ if $\lambda>\lambda_*$. We fix $\lambda$ with $|\lambda|>\lambda_*$, so all the Lyapunov exponents of any compact invariant subset of $\mathcal{A}_{\lambda}$ are strictly negative. Let $\mathcal{S}\subseteq\Omega$ be a $\sigma$-minimal set, and consider the family $\eqref{eq:parametricgeneraleq}_\lambda$ for $\omega\in\mathcal{S}$. The corresponding attractor is $\mathcal{A}_\lambda^\mathcal{S}=\{(\omega,x)\in\mathcal{A}_\lambda:\;\omega\in\mathcal{S}\}$, and Theorem~3.4 of \cite{lineardissipativescalar} guarantees that $\mathcal{A}_\lambda^\mathcal{S}$ is an attractive hyperbolic copy of the base of $(\mathcal{S}\times\mathbb{R},\tau)$. Let us check the same property for $\mathcal{A}_\lambda$, for which we proceed as in the second paragraph of the proof of Theorem~\ref{th:new}(ii). The only difference is that, given $(\omega,x)\in\mathcal{A}_\lambda$, we choose a point $\omega_0$ in a minimal subset $\mathcal{S}$ of the $\boldsymbol\upalpha$-limit of $\omega$ for $\sigma$, and a sequence $\{t_n\}_{n\in\mathbb{N}}\downarrow-\infty$ with $\omega_0=\lim_{n\rightarrow\infty} \omega{\cdot}t_n$ and such that there exists $\lim_{n\rightarrow\infty} u(t_n,\omega,x)$. Then, this limit is $\alpha_\lambda(\omega_0)$, the unique element of $(\mathcal{A}_\lambda^\mathcal{S})_{\omega_0}$. The rest of the argument is identical.
\end{proof}
The \emph{upper semicontinuity of the global attractors at $\lambda=\xi$} (see Chapter 3 of \cite{carvalho1}), defined by $\lim_{\lambda\rightarrow\xi}\text{dist}(\mathcal{A}_\lambda,\mathcal{A}_{\xi})=0$, is well-known for every $\xi\in\mathbb{R}$. In addition, the simultaneous continuity of $\lambda\mapsto\alpha_\lambda(\omega)$ and $\lambda\mapsto\beta_\lambda(\omega)$ at $\xi$ for all $\omega\in\Omega$ implies $\lim_{\lambda\rightarrow\xi}\text{dist}_\mathcal{H}(\mathcal{A}_\lambda,\mathcal{A}_{\xi})=0$, that is, the \emph{continuous variation of attractors at $\lambda=\xi$}. However, this is not a general property.
\begin{proposition} Let $f\in C^{0,1}(\Omega\times\mathbb{R},\mathbb{R})$ be $\mathrm{(Co)}$. Given $m\in\mathfrak{M}_\mathrm{erg}(\Omega,\sigma)$, there exists a countable set $\Sigma\subset\mathbb{R}$ such that $\lambda\mapsto\alpha_\lambda(\omega)$ and $\lambda\mapsto\beta_\lambda(\omega)$ are continuous at $\lambda=\xi$ for $m$-a.e. $\omega\in\Omega$ if $\xi\not\in\Sigma$.
\end{proposition}
\begin{proof} Recall that Theorem~\ref{th:attractorproperties}(i) ensures that $\lambda\mapsto\alpha_\lambda(\omega)$ is left-continuous and $\lambda\mapsto\beta_\lambda(\omega)$ is right-continuous for all $\omega\in\Omega$. It remains to show continuity in the other direction. Theorem~\ref{th:attractorproperties}(i) ensures that $\lambda\mapsto\int_\Omega \beta_\lambda(\omega)\, dm$ is an increasing function, so its set $\Sigma\subset\mathbb{R}$ of discontinuity points is at most countable. For $\xi\not\in\Sigma$, we have $\lim_{\lambda\uparrow\xi}\int_\Omega|\beta_\xi(\omega)-\beta_\lambda(\omega)|\, dm=\lim_{\lambda\uparrow\xi}\int_\Omega\big(\beta_\xi(\omega)-\beta_\lambda(\omega)\big)\, dm=0$, and hence $\lim_{\lambda\uparrow\xi}\beta_\lambda=\beta_\xi$ in $L^1(\Omega,m)$. Consequently, there exists $\{\lambda_n\}_{n\in\mathbb{N}}\uparrow\xi$ such that $\lim_{n\rightarrow\infty} \beta_{\lambda_n}(\omega)=\beta_\xi(\omega)$ for $m$-a.e. $\omega\in\Omega$, and the monotonicity ensures that $\lim_{\lambda\uparrow\xi}\beta_\lambda(\omega)=\beta_\xi(\omega)$ for $m$-a.e. $\omega\in\Omega$. We proceed analogously with $\alpha_\lambda$ and take the union of both sets.
\end{proof}
\subsection{Dealing with d-concavity} In this section, we will assume that $f$ is $\mathrm{(DC)}$ and $\mathrm{(Co)}$. We shall say that a function $f\colon\Omega\times\mathbb{R}\rightarrow\mathbb{R}$ is \emph{concave-convex} if there exist $x_l\leq x_u$ such that $f(\omega,\cdot)$ is concave in $[x_u,\infty)$ and convex in $(-\infty,x_l]$ for every $\omega\in\Omega$. The following proposition ensures that if $f$ is coercive and d-concave, it generates concave-convex dynamics. The uniformity on $\Omega$ of the coercivity property is not needed in its proof.
\begin{proposition} \label{prop:concaveconvex} Let $f\in C^{0,2}(\Omega\times\mathbb{R},\mathbb{R})$ be $\mathrm{(Co)}$. Then, for each $\omega\in\Omega$ there exist $\{x_{\omega,l}^n\}_{n\in\mathbb{N}}\downarrow-\infty$ and $\{x_{\omega,u}^n\}_{n\in\mathbb{N}}\uparrow\infty$ such that $f_{xx}(\omega,x_{\omega,l}^n)>0$ and $f_{xx}(\omega,x_{\omega,u}^n)<0$ for every $n\in\mathbb{N}$. Moreover, if $f$ is $\mathrm{(DC)}$, then there exist $x_l\leq x_u$ such that $f_{xx}(\omega,x)>0$ for all $x< x_l$ and all $\omega\in\Omega$ and such that $f_{xx}(\omega,x)<0$ for all $x> x_u$ and all $\omega\in\Omega$.
\end{proposition}
\begin{proof} We proceed with $\{x^n_{\omega,u}\}_{n\in\mathbb{N}}$. Let us fix $\omega\in\Omega$ and assume for contradiction the existence of $x_0$ such that $f_{xx}(\omega,\cdot)\geq0$ on $[x_0,\infty)$, in which case $f_x(\omega,\cdot)$ is nondecreasing on $[x_0,\infty)$. This fact and the mean value theorem ensure the existence of $\xi_x\in(x_0,x)$ for any $x>0$ such that
\begin{equation}\label{eq2}
\frac{f(\omega,x)-f(\omega,x_0)}{x-x_0}=f_x(\omega,\xi_x)\geq f_x(\omega,x_0)\,.
\end{equation}
From here, we get
\begin{equation*}
\lim_{x\rightarrow\infty}\frac{f(\omega,x)}{x-x_0}\geq \lim_{x\rightarrow\infty}\frac{f(\omega,x_0)}{x-x_0}+f_x(\omega,x_0)=f_x(\omega,x_0)>-\infty\,,
\end{equation*}
which contradicts the coercivity property. We proceed analogously to prove the existence of $\{x_{\omega,l}^n\}_{n\in\mathbb{N}}$.

Now, let us assume that $f$ is $\mathrm{(DC)}$. The previous properties and the nonincreasing character of $f_{xx}$ ensure that $I_\omega=\{x\colon f_{xx}(\omega,x)<0\}$ is a positive half-line. Let us define $g(\omega)=\inf I_\omega$, and let us check that $x_u=\sup_{\omega\in\Omega}g(\omega)$ is finite. For each $\omega_0\in\Omega$, we look for $y_0>g(\omega_0)$, so $f_{xx}(\omega_0,y_0)<0$ since $y_0\in I_{\omega_0}$. The continuity of $f_{xx}$ provides $f_{xx}(\omega,y_0)<0$ and hence $g(\omega)<y_0$ for all $\omega$ in an open neighborhood of $\omega_0$. The compactness of $\Omega$ proves that $x_u\in\mathbb{R}$. If $x>x_u$, then $x\in I_\omega$ for all $\omega\in\Omega$, and hence $f_{xx}(\omega,x)<0$, as asserted. We define $x_l$ analogously.
\end{proof}
\begin{proposition} \label{prop:510} Let $f\in C^{0,2}(\Omega\times\mathbb{R},\mathbb{R})$ be $\mathrm{(DC)}$ and $\mathrm{(Co)}$. Then, we have $\lim_{|x|\rightarrow\infty}f_x(\omega,x)=-\infty$ uniformly on $\Omega$.
\end{proposition}
\begin{proof} Let us take $a<x_l$ and $b>x_u$, where $x_l$ and $x_u$ are defined by Proposition~\ref{prop:concaveconvex}, which shows that $f_{xx}$ has a strictly negative upper bound in $[b,\infty)$ and a strictly positive lower bound in $(-\infty,a]$. Taylor's Theorem concludes the proof, taking into account the boundedness of $f_{x}(\cdot,a)$ and $f_{x}(\cdot,b)$ on $\Omega$.
\end{proof}
In particular, the last proposition shows that the hypotheses of Theorem~\ref{th:attractorproperties}(iii) are satisfied by a $C^{0,2}(\Omega\times\mathbb{R},\mathbb{R})$ function which is $\mathrm{(DC)}$ and $\mathrm{(Co)}$. The next result establishes a relation between two Lyapunov exponents of two compact sets which are $\tau_\lambda$-invariant for two different values of the parameter, which will be extremely useful in the proofs of the main results. Note that the $\mathrm{(Co)}$ condition is not required.
\begin{proposition} \label{prop:lyapexp2} Let $f\in C^{0,2}(\Omega\times\mathbb{R},\mathbb{R})$ be $\mathrm{(DC)}$, let us fix $m\in\mathfrak{M}_\mathrm{erg}(\Omega,\sigma)$ and $\lambda_\nu<\lambda_\mu$, and let $\beta_\nu,\beta_\mu\colon\Omega\rightarrow\mathbb{R}$ be two $m$-measurable equilibria for $\tau_{\lambda_\nu}$ and $\tau_{\lambda_\mu}$ respectively. If $\beta_\nu(\omega)<\beta_\mu(\omega)$ for all $m$-a.e. $\omega\in\Omega$, then
\begin{equation*}
\int_\Omega f_x(\omega,\beta_\nu(\omega))\, dm+\int_\Omega f_x(\omega,\beta_\mu(\omega))\; dm<0\,.
\end{equation*}
\end{proposition}
\begin{proof} The argument is similar to that of Proposition~\ref{prop:lyapexp1}. The function $k(\omega)=\beta_\mu(\omega)-\beta_\nu(\omega)$ satisfies
\begin{equation*}
k'(\omega{\cdot}t)/k(\omega{\cdot}t)=F(\omega{\cdot}t,k(\omega{\cdot}t))+(\lambda_\mu-\lambda_\nu)/k(\omega{\cdot}t)
\end{equation*}
for all $\omega$ at the $\sigma$-invariant set $\Omega_0$ at which $k(\omega)>0$ and $t\in\mathbb{R}$, where $F(\omega,y)=\int_0^1 f_x(\omega,sy+\beta_\nu(\omega))\, ds$. Then, the application of Birkhoff's Ergodic Theorem yields
\begin{equation*}
\int_\Omega F(\omega,k(\omega))\; dm=-(\lambda_\mu-\lambda_\nu)\int_\Omega\frac{1}{k(\omega)}\; dm<0\,,
\end{equation*}
which combined with \eqref{eq:326} completes the proof.
\end{proof}
\subsection{Bifurcation diagram with minimal base flow}\label{subsec:bifurcationdiagramminimal} Hereafter, we will always assume that $(\Omega,\sigma)$ is a minimal flow. The results in this section will describe several possible bifurcation diagrams of minimal sets for $\eqref{eq:parametricgeneraleq}_\lambda$. The bifurcations of the family of global attractors will be deduced. We point out that an (SDC)$_*$ function $f$ satisfies all the strict d-concavity conditions required in the results of this subsection (see Proposition \ref{prop:sdcimplications}), which are optimized to be less restrictive.

The next result establishes conditions under which the global bifurcation diagram presents two local saddle-node bifurcation points of minimal sets which are also points of discontinuity of the global attractor. The maps $\alpha_\lambda$ and $\beta_\lambda$ of its statement are those of $\eqref{eq:parametricattractor}_\lambda$. 
\begin{theorem}[Global bifurcation diagram with double saddle-node] \label{th:bifdiagram} Let $f\in C^{0,2}(\Omega\times\mathbb{R},\mathbb{R})$ be $(\mathrm{SDC})$ on $J$ and $\mathrm{(Co)}$, where $J$ is a compact interval which contains the interval of variation of $\mathcal{A}_\lambda$ for all $\lambda\in[-\lambda_*,\lambda_*]$ with $\lambda_*$ given by Theorem~$\mathrm{\ref{th:attractorproperties}(iii)}$ and Proposition~$\mathrm{\ref{prop:510}}$. Assume that there exists $\lambda_0\in\mathbb{R}$ such that $\tau_{\lambda_0}$ admits three different minimal sets. Then, there exists an interval $I=(\lambda_-,\lambda_+)$ with $\lambda_0\in I$ such that
\begin{enumerate}[label=\rm{(\roman*)}]
\item for every $\lambda\in I$, there exist exactly three $\tau_\lambda$-minimal sets $\mathcal{M}_\lambda^l<\mathcal{N}_\lambda<\mathcal{M}_\lambda^u$ which are hyperbolic copies of the base, given by the graphs of $\alpha_\lambda<\kappa_\lambda<\beta_\lambda$. In addition, $\mathcal{N}_\lambda$ is repulsive and $\mathcal{M}_\lambda^l$, $\mathcal{M}_\lambda^u$ are attractive, and $\lambda\mapsto\kappa_\lambda$ is strictly decreasing on $I$.
\item The graphs of $\kappa_\lambda$ and $\beta_\lambda$ (resp.~$\alpha_\lambda$ and $\kappa_\lambda$) collide on a $\sigma$-invariant residual set as $\lambda\downarrow\lambda_-$ (resp.~$\lambda\uparrow\lambda_+$), giving rise to a nonhyperbolic minimal set $\mathcal{M}_{\lambda_-}^u$ (resp.~$\mathcal{M}_{\lambda_+}^l$). In addition, there is exactly other minimal set for $\tau_{\lambda_-}$ (resp.~$\tau_{\lambda_+}$) and it is an attractive hyperbolic copy of the base given by the graph $\mathcal{M}_{\lambda_-}^l$ of $\alpha_{\lambda_-}$ (resp.~$\mathcal{M}_{\lambda_+}^u$ of $\beta_{\lambda_+}$).
\item For $\lambda\in(-\infty,\lambda_-)\cup(\lambda_+,\infty)$, $\mathcal{A}_\lambda$ is an attractive hyperbolic copy of the base, given by the graph of the map $\alpha_\lambda=\beta_\lambda$.
\end{enumerate}
In particular, local saddle-node bifurcations of minimal sets and discontinuous bifurcations of attractors occur at $\lambda_-$ and $\lambda_+$.
\end{theorem}
\begin{proof}
Since $(\Omega,\sigma)$ is minimal, Theorem~\ref{th:new} ensures that the three $\tau_{\lambda_0}$-minimal sets are the unique ones, and that they are hyperbolic copies of the base. The persistence under small $C^{0,1}$ perturbations of a hyperbolic set (see e.g. Theorem~2.8 of \cite{potz}) provides a maximal neighborhood $I$ of $\lambda_0$ such that, for any $\lambda\in I$, there are (exactly) three $\tau_\lambda$-minimal sets $\mathcal{M}_\lambda^l<\mathcal{N}_\lambda<\mathcal{M}_\lambda^u$ and they are hyperbolic copies of the base, given by the graphs of $\alpha_\lambda<\kappa_\lambda<\beta_\lambda$. This persistence also ensures the continuity of the maps $I\rightarrow C(\Omega,\mathbb{R})$, $\lambda\mapsto \alpha_\lambda,\kappa_\lambda,\beta_\lambda$ with respect to the uniform topology. Proposition~\ref{prop:precollision} ensures that $\alpha_\lambda$ and $\beta_\lambda$ are the delimiter equilibria of $\mathcal{A}_\lambda$, and hence Theorem~\ref{th:attractorproperties}(i) shows that $\lambda\mapsto\alpha_\lambda$ and $\lambda\mapsto\beta_\lambda$ are increasing on $I$. Let us check that $\lambda\mapsto\kappa_\lambda$ is strictly decreasing on $I$. We take $\xi<\lambda$ in $I$ close enough to guarantee $\alpha_\xi<\kappa_\lambda<\beta_\xi$. Since $\kappa_\lambda$ is a continuous strong time-reversed $\tau_\xi$-subequilibrium, there exist $e>0$ and $s^*<0$ such that $u_\xi(s,\omega,\kappa_\lambda(\omega))\geq\kappa_\lambda(\omega{\cdot}s)+e$ for all $\omega\in\Omega$ and $s\leq s^*$ (see Subsection~\ref{subsec:equilibria}). This implies that the graph of $\kappa_\lambda$ is strictly below any $\tau_\xi$-minimal subset $\mathcal{N}_\xi^*$ of the $\boldsymbol\upalpha$-limit set for $\tau_\xi$ of a point $(\omega,\kappa_\lambda(\omega))$. The repulsive properties of $\mathcal{M}^l_\lambda$ and $\mathcal{M}_\lambda^u$ as time decreases ensure that they are not $\boldsymbol\upalpha$-limit sets of points placed outside them. Therefore, $\mathcal{N}_\xi^*=\mathcal{N}_\xi$ and hence $\kappa_\xi>\kappa_\lambda$, as asserted.

Notice that $I\subseteq(-\lambda_*,\lambda_*)$: see Theorem~\ref{th:attractorproperties}(iii). We define $\lambda_-=\inf\{\lambda<\lambda_0\colon$ there exist three hyperbolic $\tau_\xi$-minimal sets $\forall\xi\in [\lambda,\lambda_0]\}=\inf I\in[-\lambda_*,\lambda_0)$. Since $I$ is open, $\lambda_-\not\in I$, and hence there exist at most two $\tau_{\lambda_-}$-minimal sets. We define $\beta_{\lambda_-}(\omega)=\lim_{\lambda\downarrow\lambda_-}\beta_\lambda(\omega)$, and analogously $\alpha_{\lambda_-}(\omega)=\lim_{\lambda\downarrow\lambda_-}\alpha_\lambda(\omega)$ and $\kappa_{\lambda_-}(\omega)=\lim_{\lambda\downarrow\lambda_-}\kappa_\lambda(\omega)$. As they are monotone limits of continuous functions, they are semicontinuous on $\Omega$. In particular, $\alpha_{\lambda_-}$ and $\beta_{\lambda_-}$ are upper semicontinuous and $\kappa_{\lambda_-}$ is lower semicontinuous. The continuous variation with respect to $\lambda$ ensures
\begin{equation*}
\beta_{\lambda_-}(\omega{\cdot}t)=\lim_{\lambda\downarrow\lambda_-}\beta_\lambda(\omega{\cdot}t)=\lim_{\lambda\downarrow\lambda_-} u_\lambda(t,\omega,\beta_\lambda(\omega))=u_{\lambda_-}(t,\omega,\beta_{\lambda_-}(\omega))\,,
\end{equation*}
that is, $\beta_{\lambda_-}$ is a $\tau_{\lambda_-}$-equilibrium. The same holds for $\alpha_{\lambda_-}$ and $\kappa_{\lambda_-}$. Theorem~\ref{th:attractorproperties}(i) shows that $\beta_{\lambda_-}$ is the same map that in $\eqref{eq:parametricattractor}_{\lambda_-}$, so the notation is coherent. Now, Proposition~\ref{prop:lyapexp2} applied to $\kappa_{\lambda_0}$ (whose graph is a repulsive hyperbolic copy of the base) and to any ergodic measure concentrated on the minimal set $\mathcal{M}_{\lambda_-}^l=\mathrm{cl}_{\Omega\times\mathbb{R}}\{(\omega_0{\cdot}t,\alpha_{\lambda_-}(\omega_0{\cdot}t))\colon\; t\in\mathbb{R}\}$, being $\omega_0$ a continuity point of $\alpha_{\lambda_-}$, ensures that every Lyapunov exponent of this minimal set is strictly negative. Consequently, $\mathcal{M}^l_{\lambda_-}$ is an attractive hyperbolic copy of the base, so $\alpha_{\lambda_-}$ is, in fact, continuous and the hyperbolic continuation of $\alpha_\lambda$ as $\lambda\rightarrow\lambda_-$.

Let $\mathcal{M}$ be the minimal set associated to $\kappa_{\lambda_-}$ by \eqref{eq:precollision}. The strict monotonicity properties of $\alpha_\lambda$ and $\kappa_\lambda$ ensure that $\mathcal{M}_{\lambda_-}^l<\mathcal{M}$. Therefore, $\mathcal{M}_{\lambda_-}^l$ and $\mathcal{M}_{\lambda_-}^u=\mathcal{M}$ are the two unique $\tau_{\lambda_-}$-minimal sets. This fact and Proposition~\ref{prop:precollision} ensure these properties: that $\alpha_{\lambda_-}$ is the lower delimiter of $\mathcal{A}_{\lambda_-}$, and that $\mathcal{M}_{\lambda_-}^u$ is also associated to $\beta_{\lambda_-}$ by \eqref{eq:precollision}. In particular, $\beta_{\lambda_-}(\omega)=\kappa_{\lambda_-}(\omega)$ for all $\omega$ in the residual set $\mathcal{R}$ of their common continuity points, and hence $(\mathcal{M}_{\lambda_-}^u)_\omega$ is a singleton for all $\omega\in\mathcal{R}$. Note also that $\mathcal{M}_{\lambda_-}^u$ is contained in the set $\bigcup_{\omega\in\Omega} \big(\{\omega\}\times [\kappa_{\lambda_-}(\omega),\beta_{\lambda_-}(\omega)]\big)$, which is a compact $\tau_{\lambda_-}$-invariant pinched set.

The hyperbolic attractive character of $\alpha_{\lambda_-}$ yields the first inequality of
\begin{equation}\label{eq:formulaconnumero}
\begin{split}
&\int_\Omega f_x(\omega,\alpha_{\lambda_-}(\omega))\, dm<0,\quad\int_\Omega f_x(\omega,\kappa_{\lambda_-}(\omega))\, dm\geq0\\ &\qquad\qquad\qquad\text{and}\quad \int_\Omega f_x(\omega,\beta_{\lambda_-}(\omega))\, dm\leq0\,.
\end{split}
\end{equation}
The other ones follow from Lebesgue's Monotone Convergence Theorem. Let us check that $\mathcal{A}_\lambda$ is an attractive hyperbolic copy of the base for all $\lambda<\lambda_-$. If $\lambda<\lambda_-$, Theorem~\ref{th:attractorproperties}(i) and Proposition~\ref{prop:atravesar} provide $e>0$ and $s>0$ such that
\begin{equation*}
\kappa_{\lambda_-}(\omega)-e>u_\lambda(s,\omega{\cdot}(-s),\beta_{\lambda_-}(\omega{\cdot}(-s)))>u_\lambda(s,\omega{\cdot}(-s),\beta_\lambda(\omega{\cdot}(-s)))=\beta_{\lambda}(\omega)
\end{equation*}
for all $\omega\in\Omega$. In particular, any $\tau_\lambda$-equilibrium $\kappa:\Omega\rightarrow\mathbb{R}$ is strictly below $\kappa_{\lambda_-}$. Proposition~\ref{prop:lyapexp2} and the second inequality of (\ref{eq:formulaconnumero}) ensure that $\int_\Omega f_x(\omega,\kappa(\omega))\, dm<0$ for all $m\in\mathfrak{M}_\mathrm{erg}(\Omega,\sigma)$, which combined with (\ref{eq:new}) ensures that all the Lyapunov exponents of $\mathcal{A}_\lambda$ are strictly negative. This property and Theorem~3.4 of \cite{lineardissipativescalar} prove the assertion. Note that $\mathcal{M}_{\lambda_-}^u$ is nonhyperbolic, since it does not persist for $\lambda<\lambda_-$.

The same arguments for $\lambda_+=\sup I$ complete the proof of (i), (ii) and (iii). Note that a saddle-node of minimal sets takes place at $\lambda_-$ (resp.~$\lambda_+$), as the two minimal sets which collide at that value of the parameter actually disappear for $\lambda<\lambda_-$ (resp.~$\lambda>\lambda_+$). In addition, $\lambda\mapsto\beta_\lambda(\omega)$ and $\lambda\mapsto \alpha_\lambda(\omega)$ are respectively discontinuous at $\lambda_-$ and $\lambda_+$ for all $\omega\in\Omega$, which shows the ``strong'' discontinuity of the attractor at both bifurcation points.
\end{proof}
It is interesting to note that the dynamics at the nonhyperbolic minimal set occurring at any of the bifurcation points of the previous theorem can be highly complicated. The next result contributes to understand the possibilities for this complex dynamics. In particular,
strictly positive and strictly negative Lyapunov exponents can coexist on that set.
\begin{proposition} \label{prop:alternatives} Let us add to the hypotheses of Theorem~$\mathrm{\ref{th:bifdiagram}}$ that $f$ is $\mathrm{(SDC)}_m$ on the variation interval $J_{\mathcal{A}_{\lambda_-}}$ of $\mathcal{A}_{\lambda_-}$ for a measure $m\in\mathfrak{M}_\mathrm{erg}(\Omega,\sigma)$. Then,
\begin{enumerate}[label=\rm{(\roman*)}]
\item one of the following situations holds:
\begin{enumerate}[label=\rm{(\alph*)}]
\item $m(\{\omega\colon\;\kappa_{\lambda_-}(\omega)=\beta_{\lambda_-}(\omega)\})=1$ and
\begin{equation*}
\int_\Omega f_x(\omega,\beta_{\lambda_-}(\omega))\; dm=\int_\Omega f_x(\omega,\kappa_{\lambda_-}(\omega))\; dm=0\,,
\end{equation*}
\item $m(\{\omega\colon\;\kappa_{\lambda_-}(\omega)=\beta_{\lambda_-}(\omega)\})=0$,
\begin{equation*}
\int_\Omega f_x(\omega,\beta_{\lambda_-}(\omega))\; dm<0\quad\text{and}\quad \int_\Omega f_x(\omega,\kappa_{\lambda_-}(\omega))\; dm>0\,.
\end{equation*}
\end{enumerate}
Moreover, situation $\mathrm{(a)}$ holds if and only if all the Lyapunov exponents of the compact $\tau_{\lambda_-}$-invariant set $\bigcup_{\omega\in\Omega} \big(\{\omega\}\times [\kappa_{\lambda_-}(\omega),\beta_{\lambda_-}(\omega)]\big)$ are zero.
\item Let $\alpha_\mathcal{M}$ and $\beta_\mathcal{M}$ be the delimiting equilibria of $\mathcal{M}=\mathcal{M}_{\lambda_-}^u$. Then, $\alpha_\mathcal{M}$ (resp.~$\beta_\mathcal{M}$) coincides with $\kappa_{\lambda_-}$ or $\beta_{\lambda_-}$ $m$-a.e. If, in addition, $\mathfrak{M}_\mathrm{erg}(\Omega,\sigma)=\{m\}$, then $\alpha_\mathcal{M}=\kappa_{\lambda_-}$ and $\beta_\mathcal{M}=\beta_{\lambda_-}$ $m$-a.e.
\end{enumerate}
Analogous properties hold for $\lambda_+$.
\end{proposition}
\begin{proof} (i) As $\{\omega\colon\;\kappa_{\lambda_-}(\omega)=\beta_{\lambda_-}(\omega)\}$ is $\sigma$-invariant and $m\in\mathfrak{M}_\mathrm{erg}(\Omega,\sigma)$, then $m(\{\omega\colon\;\kappa_{\lambda_-}(\omega)=\beta_{\lambda_-}(\omega)\})\in\{0,1\}$. If $m(\{\omega\colon\;\kappa_{\lambda_-}(\omega)=\beta_{\lambda_-}(\omega)\})=1$, then \eqref{eq:formulaconnumero} yields $0\leq\int_\Omega f_x(\omega,\beta_{\lambda_-}(\omega))\; dm=\int_\Omega f_x(\omega,\kappa_{\lambda_-}(\omega))\; dm\leq0$, and hence (a) holds. If $m(\{\omega\colon\;\kappa_{\lambda_-}(\omega)=\beta_{\lambda_-}(\omega)\})=0$, then $\alpha_{\lambda_-}$, $\kappa_{\lambda_-}$ and $\beta_{\lambda_-}$ define three different elements of $\mathfrak{M}_\mathrm{erg}(\mathcal{A}_{\lambda_-},\tau_{\lambda_-})$ projecting onto $m$. By reviewing the proof of Theorem~\ref{th:numergmeas}, we see that $\int_\Omega f_x(\omega,\beta_{\lambda_-}(\omega))\; dm<0$ and $\int_\Omega f_x(\omega,\kappa_{\lambda_-}(\omega))\; dm>0$, as asserted. Let us prove the last assertion. If situation (a) holds, then every equilibria with graph contained in $\bigcup_{\omega\in\Omega} \big(\{\omega\}\times [\kappa_{\lambda_-}(\omega),\beta_{\lambda_-}(\omega)]\big)$ must coincide $m$-a.e. with $\kappa_{\lambda_-}$ and $\beta_{\lambda_-}$, so it defines a zero exponent. Conversely, if situation (b) holds, then there exist nonzero exponents of $\bigcup_{\omega\in\Omega} \big(\{\omega\}\times [\kappa_{\lambda_-}(\omega),\beta_{\lambda_-}(\omega)]\big)$.

(ii) According to Theorem~\ref{th:numergmeas}, $\mathcal{A}_{\lambda_-}$ concentrates at most three ergodic measures projecting onto $m$. Since $\mathcal{M}_{\lambda_-}^l$ concentrates one, then the $\tau_{\lambda_-}$-equilibria $\kappa_{\lambda_-}\leq\alpha_\mathcal{M}\leq\beta_\mathcal{M}\leq\beta_{\lambda_-}$ can define at most two (by \eqref{eq:new}), which ensures the first assertion in (ii). Now, let $(\Omega,\sigma)$ be uniquely ergodic. If $m(\{\omega\colon\; \kappa_{\lambda_-}(\omega)=\beta_{\lambda_-}(\omega)\})=1$, then the four equilibria define the same measure, so that they coincide $m$-a.e. If $m(\{\omega\colon\;\kappa_{\lambda_-}(\omega)=\beta_{\lambda_-}(\omega)\})=0$, then $\kappa_{\lambda_-}$ and $\beta_{\lambda_-}$ define the two remaining measures. By contradiction, assume that $\alpha_\mathcal{M}=\beta_\mathcal{M}$ $m$-a.e. and that they coincide $m$-a.e. with $\beta_{\lambda_-}$ (resp.~with $\kappa_{\lambda_-}$). Then, (i)(b) ensures that the unique Lyapunov exponent of $\mathcal{M}$ is negative (resp.~positive), and hence $\mathcal{M}$ is a hyperbolic copy of the base which contradicts Theorem~\ref{th:bifdiagram}(ii). In consequence, $\alpha_\mathcal{M}=\kappa_{\lambda_-}$ and $\beta_\mathcal{M}=\beta_{\lambda_-}$ $m$-a.e.
\end{proof}
The next two results show an alternative hypothesis to get the bifurcation diagram of Theorem 5.12, and conditions under which this alternative hypothesis holds unless the flow $\tau_\lambda$ admits a unique minimal set for all the values of the parameter; that is, conditions under which there are only two possible global bifurcation diagrams: that of Theorem~\ref{th:bifdiagram}, and that of Theorems~\ref{prop:orderedattractors} and \ref{th:onlyoneminimal} below. The definition of the Sacker and Sell spectrum is given in Subsection~\ref{subsec:furstenbergremark}.
\begin{theorem} \label{th:Sformation} Let $f\in C^{0,2}(\Omega\times\mathbb{R},\mathbb{R})$ be $(\mathrm{SDC})$ on $J$ and $\mathrm{(Co)}$, where $J$ is a compact interval which contains the interval of variation of $\mathcal{A}_\lambda$ for all $\lambda\in[-\lambda_*,\lambda_*]$, with $\lambda_*$ given by Theorem~$\mathrm{\ref{th:attractorproperties}(iii)}$. Assume the existence of $\xi\in\mathbb{R}$ such that $\tau_\xi$ admits exactly two different minimal sets $\mathcal{M}_1<\mathcal{M}_2$ with $\mathcal{M}_1$ (resp.~$\mathcal{M}_2$) hyperbolic. Then, the bifurcation diagram of $\eqref{eq:parametricgeneraleq}_\lambda$ is that described by Theorem~$\mathrm{\ref{th:bifdiagram}}$, with $\lambda_-=\xi$ (resp.~$\lambda_+=\xi$).
\end{theorem}
\begin{proof} Let us suppose that the lower $\tau_\xi$-minimal set $\mathcal{M}_1\subseteq\mathcal{A}_\xi$ is the hyperbolic one and let $\mathcal{M}_2$ be the upper $\tau_\xi$-minimal set, with $\alpha(\omega)=\inf(\mathcal{M}_2)_\omega$ and $\beta(\omega)=\sup(\mathcal{M}_2)_\omega$. The other case follows analogously. Note that Proposition~\ref{prop:hyperbolicpositiveminimal} guarantees that $\mathcal{M}_1$ is hyperbolic attractive, as otherwise there would exist three $\tau_\xi$-minimal sets. Let $\lambda_0>\xi$ be close enough to $\xi$ to guarantee the existence of an attractive hyperbolic minimal set $\mathcal{M}^l_{\lambda_0}<\mathcal{M}_2$, obtained by hyperbolic continuation of $\mathcal{M}_1$. Proposition~\ref{prop:alphasuperequilibrium} ensures that $\alpha$ and $\beta$ are semicontinuous strong $\tau_{\lambda_0}$-subequilibria and Proposition~\ref{prop:atravesar} ensures that there exist $s_1>0$ and $e_1>0$ such that $\beta(\omega)+e_1<(\alpha)_{s_1}^{\lambda_0}(\omega)=u_{\lambda_0}(s_1,\omega{\cdot}(-s_1),\alpha(\omega{\cdot}(-s_1)))$ for all $\omega\in\Omega$. Let $\mathcal{M}^u_{\lambda_0}$ be the $\boldsymbol\upomega$-limit set for $\tau_{\lambda_0}$ of a point $(\omega_0,\alpha(\omega_0))$ of the graph of $\alpha$, and note that it coincides with the $\boldsymbol\upomega$-limit set for $\tau_{\lambda_0}$ of the point $(\omega_0{\cdot}s_1,(\alpha)_{s_1}^{\lambda_0}(\omega_0{\cdot}s_1))$ of the graph of $(\alpha)_{s_1}^{\lambda_0}$. As $(\alpha)_{s_1}^{\lambda_0}$ is a lower semicontinuous subequilibria, taking into account the inequality provided by Proposition~\ref{prop:atravesar}, we get $\mathcal{M}_2<\mathcal{M}^u_{\lambda_0}$. Using now that $\alpha$ and $\beta$ are semicontinuous strong time-reversed $\tau_{\lambda_0}$-superequilibria, and the analogous of Proposition~\ref{prop:atravesar} for the time-reversed case, we check that a minimal subset $\mathcal{N}_{\lambda_0}$ contained in the $\boldsymbol\upalpha$-limit set for $\tau_{\lambda_0}$ of some point $(\omega_0,\beta(\omega_0))$ of the graph of $\beta$ satisfies $\mathcal{N}_{\lambda_0}<\mathcal{M}_2$. In addition, $\mathcal{N}_{\lambda_0}\neq\mathcal{M}_{\lambda_0}^l$, due to the repulsive properties of $\mathcal{M}_{\lambda_0}^l$ as time decreases (see e.g. Proposition 2.8 of~\cite{lineardissipativescalar}). Hence, there exist three different $\tau_{\lambda_0}$-minimal sets: all the hypotheses of Theorem~\ref{th:bifdiagram} hold.
\end{proof}
\begin{theorem} \label{th:twominimalsets} Let $f\in C^{0,2}(\Omega\times\mathbb{R},\mathbb{R})$ be $\mathrm{(SDC)}$ and $\mathrm{(Co)}$. Assume that the flow $\tau$ defined by \eqref{eq:flujotau} admits at least two minimal sets $\mathcal{M}_1$ and $\mathcal{M}_2$.
\begin{enumerate}[label=\rm{(\roman*)}]
\item If the Sacker and Sell spectrum of $f_x$ on $\mathcal{M}_1$ is contained in $[0,\infty)$, then $\mathcal{M}_2$ is hyperbolic attractive.
\item If the Sacker and Sell spectrum of $f_x$ on one of the minimal sets reduces to a point, then either $\mathcal{M}_1$ or $\mathcal{M}_2$ is hyperbolic attractive.
\item If $(\Omega,\sigma)$ is uniquely ergodic, then either $\mathcal{M}_1$ or $\mathcal{M}_2$ is hyperbolic attractive.
\end{enumerate}
\end{theorem}
\begin{proof} (i) We will prove that every Lyapunov exponent of $\mathcal{M}_2$ is strictly negative (see Subsection~\ref{subsec:furstenbergremark}). Let $\nu\in\mathfrak{M}_\mathrm{erg}(\mathcal{M}_2,\tau)$ projecting onto $m\in\mathfrak{M}_\mathrm{erg}(\Omega,\sigma)$, let $\beta_\nu:\Omega\rightarrow\mathbb{R}$ be an $m$-measurable equilibrium with graph contained in $\mathcal{M}_2$ satisfying \eqref{eq:new} for $\nu$, and let $\mu\in\mathfrak{M}_\mathrm{erg}(\mathcal{M}_1,\tau)$ be defined by the upper delimiter $\beta_\mu$ of $\mathcal{M}_1$, $m$, and \eqref{eq:new}. As the hypothesis of (i) ensures that $\gamma(\mu,\mathcal{M}_1)\geq0$, Proposition~\ref{prop:lyapexp1} applied to $\beta_\nu$ and $\beta_\mu$ ensures that $\gamma(\nu,\mathcal{M}_2)<0$, as asserted.

(ii) Suppose that the Sacker and Sell spectrum of $f_x$ on $\mathcal{M}_1$ reduces to a point $\{a\}$. If $a<0$, then $\mathcal{M}_1$ is hyperbolic attractive. If not, then $\{a\}\subset[0,\infty)$, and (i) proves the assertion.

(iii) In the uniquely ergodic case, if $\mathcal{M}_1$ is not hyperbolic attractive, then there exists $\mu\in\mathfrak{M}_\mathrm{erg}(\mathcal{M}_1,\tau)$ projecting onto the unique measure $m\in\mathfrak{M}_\mathrm{erg}(\Omega,\sigma)$ such that $\gamma(\mu,\mathcal{M}_1)\geq0$. Since every ergodic measure on $\mathcal{M}_2$ projects onto $m$, Proposition~\ref{prop:lyapexp1} guarantees that $\gamma(\nu,\mathcal{M}_2)<0$ for every $\nu\in\mathfrak{M}_\mathrm{erg}(\mathcal{M}_2,\tau)$.
\end{proof}
The simplest global bifurcation diagram of minimal sets (without bifurcation points) occurs when the flow $\tau_\lambda$ admits only one minimal set for every value of the parameter. The next two results analyze this situation.
\begin{theorem}\label{prop:orderedattractors} Let $f\in C^{0,1}(\Omega\times\mathbb{R},\mathbb{R})$ be $\mathrm{(Co)}$ such that there exists only one $\tau_\lambda$-minimal set for all $\lambda\in\mathbb{R}$. Then, $\mathcal{A}_\lambda$ is pinched for all $\lambda\in\mathbb{R}$, and $\mathcal{A}_\lambda<\mathcal{A}_\xi$ (i.e.,~$\beta_\lambda<\alpha_\xi$) if $\lambda<\xi$.
\end{theorem}
\begin{proof} Propositions~\ref{prop:precollision} and \ref{prop:collision} applied to the delimiter equilibria $\alpha_\lambda$ and $\beta_\lambda$ of $\mathcal{A}_\lambda$ show that they coincide in a residual set of common continuity points. Proposition~\ref{prop:alphasuperequilibrium} shows that $\alpha_\lambda$ and $\beta_\lambda$ are semicontinuous strong $\tau_\xi$-subequilibria if $\xi>\lambda$, so Proposition~\ref{prop:atravesar} provides $e>0$ such that $\beta_\lambda(\omega)+e<\alpha_\xi(\omega)$ for all $\omega\in\Omega$.
\end{proof}
\begin{theorem} \label{th:onlyoneminimal} Let $f\in C^{0,2}(\Omega\times\mathbb{R},\mathbb{R})$ be $\mathrm{(DC)}$ and $\mathrm{(Co)}$ such that there exists only one $\tau_\lambda$-minimal set for all $\lambda\in\mathbb{R}$.
\begin{enumerate}[label=\rm{(\roman*)}]
\item If $(\Omega,\sigma)$ is uniquely ergodic (resp.~finitely ergodic), then there exists at most a value (resp.~a finite number of values) of the parameter at which the minimal set is a nonhyperbolic copy of the base.
\item If there exists $\lambda_0\in\mathbb{R}$ such that the Sacker and Sell spectrum of $f_x$ on $\mathcal{A}_{\lambda_0}$ is $\{0\}$, then $\lambda_0$ is the only value of the parameter at which the minimal set is nonhyperbolic.
\end{enumerate}
\end{theorem}
\begin{proof} (i) Assume that $\mathcal{M}_i$ is a nonhyperbolic $\tau_{\lambda_i}$-minimal set for $i\in\{1,2\}$, with $\lambda_1<\lambda_2$. Then, $\mathcal{M}_i$ concentrates a $\tau_i$-ergodic measure giving rise to a nonnegative Lyapunov exponent. Let $m_i\in\mathfrak{M}_\mathrm{erg}(\Omega,\sigma)$ be the projection of this measure. Theorem~\ref{prop:orderedattractors} shows that $\mathcal{M}_{\lambda_1}<\mathcal{M}_{\lambda_2}$, and hence the previous property and Proposition~\ref{prop:lyapexp2} show that $m_1\neq m_2$. The assertions follow easily.

(ii) The hypothesis ensures that all the Lyapunov exponents of the only $\tau_{\lambda_0}$-minimal set are zero, so it is nonhyperbolic and we can reason as in (i).
\end{proof}
The last result of this section is a local bifurcation theorem: at least a local saddle-node bifurcation of minimal sets occurs under the $\mathrm{(SDC)}$ hypothesis on an interval which contains the variation interval of two minimal sets with some hyperbolicity properties. Recall that $(\Omega,\sigma)$ is assumed to be minimal.
\begin{theorem}[Local saddle-node bifurcations] Let $f\in C^{0,2}(\Omega\times\mathbb{R},\mathbb{R})$ be $\mathrm{(SDC)}$ on a compact interval $J$ and let us suppose that the flow $\tau$ given by $\eqref{eq:parametricgeneraleq}_0$ admits two minimal sets $\mathcal{M}_1<\mathcal{M}_2$ contained on $\Omega\times\mathrm{int}\,J$. Then,
\begin{enumerate}[label=\rm{(\roman*)}]
\item if $\mathcal{M}_1$ is hyperbolic attractive or $\mathcal{M}_2$ is hyperbolic repulsive, then there exists $\lambda_+>0$ such that $\eqref{eq:parametricgeneraleq}_\lambda$ exhibits a local saddle-node bifurcation of minimal sets at $\lambda_+$.
\item If $\mathcal{M}_2$ is hyperbolic attractive or $\mathcal{M}_1$ is hyperbolic repulsive, then there exists $\lambda_-<0$ such that $\eqref{eq:parametricgeneraleq}_\lambda$ exhibits a local saddle-node bifurcation of minimal sets at $\lambda_-$.
\item If both $\mathcal{M}_1$ and $\mathcal{M}_2$ are hyperbolic attractive, then there exists an intermediate repulsive hyperbolic minimal set $\mathcal{M}$, and two local saddle-node bifurcations of minimal sets take place at $\lambda_-$ and $\lambda_+$, with $\lambda_-<0<\lambda_+$.
\end{enumerate}
\end{theorem}
\begin{proof} We define $\widetilde{f}$ by extending $f$ outside of $\Omega\times J$ as in \eqref{eq:foutsideextension}. Then, $\widetilde{f}\in C^{0,2}(\Omega\times\mathbb{R},\mathbb{R})$ is $\mathrm{(Co)}$ and $\mathrm{(SDC)}$. In addition, the minimal sets for the flow $\widetilde\tau_\lambda$ provided by $\widetilde f+\lambda$ and contained in $\Omega\times J$ are also minimal sets for the flow $\tau_\lambda$ given by $f+\lambda$, and hence a local saddle-node bifurcation of minimal sets (\emph{lsnb} for short) for $\widetilde\tau_\lambda$ at some $\lambda_0\in\mathbb{R}$ taking place on $\Omega\times\mathrm{int}\, J$ ensures an lsnb for $\tau_\lambda$ at $\lambda_0$.

Assume that $\widetilde\tau_0$ admits three minimal sets $\mathcal{N}_1<\mathcal{N}_2<\mathcal{N}_3$. Theorem~\ref{th:bifdiagram} provides an lsnb at $\lambda_-<0$ on the open band of $\Omega\times\mathbb{R}$ delimited by $\mathcal{N}_2$ and $\mathcal{N}_3$, as well as an lsnb at $\lambda_+>0$ on the open band delimited by $\mathcal{N}_1$ and $\mathcal{N}_2$. Under the hypotheses of (iii), Proposition~\ref{prop:hyperbolicpositiveminimal} provides a $\widetilde\tau_0$-minimal set $\mathcal{M}$ with $\mathcal{M}_1<\mathcal{M}<\mathcal{M}_2$, and hence we have two lsnb on $\Omega\times\mathrm{int}\, J$ at $\lambda_-<0$ and $\lambda_+>0$. If $\mathcal{M}_1$ (resp.~$\mathcal{M}_2$) is hyperbolic repulsive, Proposition~\ref{prop:hyperbolicpositiveminimal} provides $\mathcal{M}$ with $\mathcal{M}<\mathcal{M}_1<\mathcal{M}_2$ (resp.~$\mathcal{M}_1<\mathcal{M}_2<\mathcal{M}$), and hence we have at least an lsnb on $\Omega\times\mathrm{int}\, J$ at $\lambda_-<0$ (resp.~$\lambda_+>0$). If $\mathcal{M}_1$ (resp.~$\mathcal{M}_2$) is hyperbolic attractive but $\mathcal{M}_2$ (resp.~$\mathcal{M}_1$) is nonhyperbolic, Theorems~\ref{th:Sformation} and \ref{th:bifdiagram} show that $\widetilde\tau_\lambda$ has an lsnb on $\Omega\times\mathrm{int}\, J$ at $\lambda_+>0$ (resp.~$\lambda_-<0$). These properties prove (i) and (ii).
\end{proof}
\section{A Second One-Parametric Bifurcation Problem}\label{sec:secondoneparametric}
In this section we shall bring forward a technical procedure, which is the change of skewproduct base. This procedure will allow us to study the different one-parametric bifurcation problems which we will describe in Subsection~\ref{subsec:biflambdax}.
\subsection{Procedure of change of skewproduct base}\label{subsec:sumslyap}
The main purposes of the technical procedure described in this subsection are to transform a general minimal set for a given flow in a copy of the base for an extended flow and to find the relations between the Lyapunov exponents for these two flows.

Let $(\Omega\times\mathbb{R},\tau_1)$ and $(\Omega\times\mathbb{R},\tau_2)$ be local skewproduct flows on $\Omega\times\mathbb{R}$, given by
\begin{equation*}
\tau_i\colon\mathcal{U}_i\subseteq\mathbb{R}\times\Omega\times\mathbb{R}\rightarrow\Omega\times\mathbb{R}\,,\;(t,\omega,x_0)\rightarrow(\omega{\cdot}t,u_i(t,\omega,x_0))\,,
\end{equation*}
for $i\in\{1,2\}$, where $\mathcal{I}^i_{\omega,x_0}\rightarrow\mathbb{R}$, $t\mapsto u_i(t,\omega,x_0)$ represents the maximal solution of $x'=f_i(\omega{\cdot}t,x)$ satisfying $u_i(0,\omega,x_0)=x_0$, with $f_i\in C^{0,1}(\Omega\times\mathbb{R},\mathbb{R})$, and $\mathcal{U}_i=\bigcup_{(\omega,x_0)\in\Omega\times\mathbb{R}} (\mathcal{I}^i_{\omega,x_0}\times\{(\omega,x_0)\})$ for $i\in\{1,2\}$. We allow $f_2$ to be different from $f_1$ to obtain a more general framework, although frequently the assumption $\tau_1=\tau_2$ simplifies the scenario.

Let $\Upsilon\subset\Omega\times\mathbb{R}$ be a compact $\tau_1$-invariant set projecting onto $\Omega$. As $\Upsilon$ is composed by global $\tau_1$-orbits, $(\Upsilon,\tau_1)$ is a global flow, which will provide the base of the new local skewproduct flow. We represent $\upsilon{\cdot}t=\tau_1(t,\upsilon)$ for $\upsilon=(\omega,z_0)\in\Upsilon$ and consider the local skewproduct flow defined on $\Upsilon\times\mathbb{R}$ by
\begin{equation*}
\phi_{\tau_1,\tau_2}\colon\mathcal{U}_{\tau_1,\tau_2}\subseteq\mathbb{R}\times \Upsilon\times\mathbb{R}\rightarrow \Upsilon\times\mathbb{R},\;(t,\upsilon,x_0)\mapsto \big(\upsilon{\cdot}t,u_2(t,\pi(\upsilon),x_0)\big)\,,
\end{equation*}
where $\pi\colon\Upsilon\rightarrow\Omega$, $\upsilon=(\omega,z)\mapsto\omega$ and $\mathcal{U}_{\tau_1,\tau_2}=\bigcup_{(\upsilon,x_0)\in\Upsilon\times\mathbb{R}} (\mathcal{I}^2_{\pi(\upsilon),x_0}\times\{(\upsilon,x_0)\})$. We shall say that $(\Upsilon\times\mathbb{R},\phi_{\tau_1,\tau_2})$ is the \emph{local skewproduct flow obtained from $(\Omega\times\mathbb{R},\tau_1)$ and $(\Omega\times\mathbb{R},\tau_2)$ by a change of base}. The second component of this flow represents the solutions of the scalar ordinary differential equation
\begin{equation}\label{eq:ly2}
x'=\widetilde f_2(\upsilon{\cdot} t,x)\,,
\end{equation}
where $\widetilde f_2(\upsilon,x)= f_2(\pi (\upsilon),x)$. Note that, if $\tau_1=\tau_2$, then the map $\Upsilon\rightarrow\mathbb{R}$, $(\omega,z)\mapsto z$ is a continuous $\phi_{\tau_1,\tau_1}$-equilibrium, that is, its graph is a copy of the base.

Now, we will discuss how ergodic measures of $(\Upsilon\times\mathbb{R},\phi_{\tau_1,\tau_2})$ are related to ergodic measures of $(\Omega\times\mathbb{R},\tau_2)$. Let $\mathcal{K}\subset\Omega\times\mathbb{R}$ be a compact $\tau_2$-invariant set. Then, we define $\widetilde{\mathcal{K}}=\{(\omega,z,x)\colon\, (\omega,z)\in\Upsilon,\,(\omega,x)\in\mathcal{K}\}\subset\Upsilon\times\mathbb{R}$. It is clear that $\widetilde{\mathcal{K}}$ is a compact $\phi_{\tau_1,\tau_2}$-invariant set and that $\widetilde\pi(\widetilde{\mathcal{K}})=\mathcal{K}$, where $\widetilde\pi\colon\Upsilon\times\mathbb{R}\rightarrow\Omega\times\mathbb{R}$, $(\omega,z,x)\mapsto (\omega,x)$ is a flow epimorphism. Let $\nu\in\mathfrak{M}_\mathrm{erg}(\Upsilon,\tau_1)$ be a fixed measure which projects onto $m\in\mathfrak{M}_\mathrm{erg}(\Omega,\sigma)$. We will establish a one-to-one correspondence between the $\phi_{\tau_1,\tau_2}$-ergodic measures on $\widetilde{\mathcal{K}}$ which project onto $\nu$ and the $\tau_2$-ergodic measures on $\mathcal{K}$ which project onto $m$. To this end, we fix an $m$-measurable equilibrium $\beta_\nu:\Omega\rightarrow\mathbb{R}$ satisfying \eqref{eq:new} for $\nu$.

Let $\widetilde\mu\in\mathfrak{M}_\mathrm{erg}(\widetilde{\mathcal{K}},\phi_{\tau_1,\tau_2})$ project onto $\nu$, and let $\widetilde\beta_{\widetilde\mu}\colon\Upsilon\rightarrow\mathbb{R}$ be a $\nu$-measurable $\phi_{\tau_1,\tau_2}$-equilibrium satisfying \eqref{eq:new} for $\widetilde\mu$. Then, for every continuous $\widetilde h\colon\widetilde{\mathcal{K}}\rightarrow\mathbb{R}$,
\begin{equation}\label{eq:3.36}
\int_{\widetilde{\mathcal{K}}} \widetilde h(\upsilon,x) \,d\widetilde\mu=\int_{\Upsilon} \widetilde h(\upsilon,\widetilde\beta_{\widetilde\mu}(\upsilon))\; d\nu=\int_\Omega \widetilde h(\omega,\beta_\nu(\omega),\widetilde\beta_{\widetilde\mu}(\omega,\beta_\nu(\omega)))\; dm\,.
\end{equation}
Notice that, as $\widetilde h(\cdot,\widetilde\beta_{\widetilde\mu}(\cdot))$ is not necessarily continuous, the second equality in \eqref{eq:3.36} is not immediate: it follows from applying Lusin's Theorem to $g(\upsilon)=\widetilde h(\upsilon,\widetilde\beta_{\widetilde\mu}(\upsilon))$. Notice also that, as $\beta_\nu$ is a $\tau_1$-equilibrium, then $(\omega{\cdot}t,\beta_\nu(\omega{\cdot}t))=(\omega,\beta_\nu(\omega)){\cdot}t$. So, if we define $\beta_\mu(\omega)=\widetilde\beta_{\widetilde\mu}(\omega,\beta_\nu(\omega))$,
\begin{equation*}
\beta_\mu(\omega{\cdot} t)=\widetilde\beta_{\widetilde\mu}((\omega,\beta_\nu(\omega)){\cdot}t)=u_2(t,\pi(\omega,\beta_\nu(\omega)),\widetilde\beta_{\widetilde\mu}(\omega,\beta_\nu(\omega))=u_2(t,\omega,\beta_\mu(\omega))\,.
\end{equation*}
Consequently, $\beta_\mu\colon\Omega\rightarrow\mathbb{R}$ is an $m$-measurable $\tau_2$-equilibrium with graph contained in $\mathcal{K}$ (as that of $\widetilde\beta_{\widetilde\mu}$ is in $\widetilde{\mathcal{K}}$), and hence $\beta_\mu$ defines a $\tau_2$-ergodic measure $\mu$ on $\mathcal{K}$ which projects onto $m$, by \eqref{eq:new}: given a continuous function $h\colon\mathcal{K}\rightarrow\mathbb{R}$, we have
\begin{equation}\label{eq:renew}
\begin{split}
\int_\mathcal{K} h(\omega,x)\, d\mu &=\int_\Omega h(\omega,\widetilde\beta_{\widetilde\mu}(\omega,\beta_\nu(\omega))\, dm\\ &=\int_\mathcal{K}\widetilde h(\omega,z,\widetilde\beta_{\widetilde\mu}(\omega,z))\, d\nu=\int_{\widetilde{\mathcal{K}}}\widetilde h(\omega,z,x)\, d\widetilde\mu
\end{split}
\end{equation}
for $\widetilde h\colon\widetilde{\mathcal{K}}\rightarrow\mathbb{R}$, $(\omega,z,x)\mapsto h(\omega,x)$. Relation \eqref{eq:3.36} has been used in the second equality. In particular, \eqref{eq:renew} shows that the definition of $\mu$ is independent of the choice of $\widetilde\beta_{\widetilde\mu}$. Notice that \eqref{eq:renew} ensures that the Lyapunov exponent of $\widetilde{\mathcal{K}}$ for \eqref{eq:ly2} with respect to $\widetilde\mu$ coincides with that of $\mathcal{K}$ for $z'=f_2(\omega{\cdot}t,z)$ with respect to $\mu$.

On the other hand, \eqref{eq:new} associates an $m$-measurable $\tau_2$-equilibrium $\beta_\mu\colon\Omega\rightarrow\mathbb{R}$ with graph in $\mathcal{K}$ to a $\tau_2$-ergodic measure $\mu$ on $\mathcal{K}$ which projects onto $m$. We define
\begin{equation*}
\int_{\widetilde{\mathcal{K}}}\widetilde h (\upsilon,x)\,d\widetilde\mu=\int_\Omega\widetilde h(\omega,\beta_\nu(\omega),\beta_\mu(\omega))\; dm
\end{equation*}
for any continuous map $\widetilde h\colon\mathcal{K}\rightarrow\mathbb{R}$, and obtain an ergodic measure $\widetilde\mu$ on $\widetilde{\mathcal{K}}$ which projects onto $\nu$. As in the previous case, it can be seen that this process is well defined, that is, it is independent of the choice of $\beta_\mu$. It is not hard to check that this process is the inverse of the aforementioned.

Now, let us suppose that both $\Omega$ and $\Upsilon\subset\Omega\times\mathbb{R}$ are minimal for their respective flows. We will construct a one-to-one correspondence between the minimal sets of $(\Upsilon\times\mathbb{R},\phi_{\tau_1,\tau_2})$ and of $(\Omega\times\mathbb{R},\tau_2)$. Recall that the sections of a minimal set on a scalar skewproduct flow with minimal base are singletons for all the points on a residual subset of the base (see Proposition~\ref{prop:precollision}). Given a $\phi_{\tau_1,\tau_2}$-minimal set $\hat{\mathcal{M}}\subset\Upsilon\times\mathbb{R}$, we define $\mathcal{M}=\widetilde\pi(\hat{\mathcal{M}})=\{(\omega,x)\colon\;(\omega,z,x)\in\hat{\mathcal{M}}\}$. It is easy to check that $\mathcal{M}\subset\Omega\times\mathbb{R}$ is $\tau_2$-minimal. Conversely, given a $\tau_2$-minimal set $\mathcal{M}\subset\Omega\times\mathbb{R}$, we look for $\omega_0\in\Omega$ such that $(\Upsilon)_{\omega_0}=\{z_0\}$ and $(\mathcal{M})_{\omega_0}=\{x_0\}$, and take $\hat{\mathcal{M}}=\mathrm{cl}_{\Upsilon\times\mathbb{R}}\{\phi_{\tau_1,\tau_2}(t,(\omega_0,z_0),x_0):\;t\in\mathbb{R}\}$. It is easy to check that $(\hat{\mathcal{M}})_{(\omega_0,z_0)}=\{x_0\}$. As $\hat{\mathcal{M}}$ is a compact $\phi_{\tau_1,\tau_2}$-invariant set, it contains a minimal set $\hat{\mathcal{M}}'\subseteq\hat{\mathcal{M}}$, and necessarily $(\hat{\mathcal{M}}')_{(\omega_0,z_0)}=\{x_0\}$. Consequently, $\hat{\mathcal{M}}\subseteq\hat{\mathcal{M}}'$, so $\hat{\mathcal{M}}$ is itself minimal. Moreover, it is easy to check that $\hat{\mathcal{M}}$ is independent of the choice of $\omega_0$ in the residual sets of points $\omega$ for which both $(\Upsilon)_\omega$ and $(\mathcal{M})_\omega$ are singletons, and also that $\widetilde\pi(\hat{\mathcal{M}})=\mathcal{M}$, since $(\omega_0,x_0)\in\widetilde\pi(\hat{\mathcal{M}})$. Finally, $\hat{\mathcal{M}}$ is the only $\phi_{\tau_1,\tau_2}$-minimal set which verifies $\widetilde\pi(\hat{\mathcal{M}})=\mathcal{M}$, since $\widetilde\pi(\hat{\mathcal{M}}')=\mathcal{M}$ would imply $(\omega_0,z_0,x_0)\in \hat{\mathcal{M}}'$. This is the one-to-one correspondence we referred to.
\subsection{A new bifurcation problem}\label{subsec:biflambdax}
Throughout this section, $(\Omega,\sigma)$ will be minimal and $f\in C^{0,1}(\Omega\times\mathbb{R},\mathbb{R})$. As explained in Subsection~\ref{subsec:skewproductflow}, the solutions of
\begin{equation}\label{eq:62primera}
x'=f(\omega{\cdot}t,x)\,,\quad\omega\in\Omega\,,
\end{equation}
induce the local skewproduct flow
\begin{equation}\label{eq:62flujo}
\tau\colon\mathcal{U}\subseteq\mathbb{R}\times\Omega\times\mathbb{R}\rightarrow\Omega\times\mathbb{R},\quad (t,\omega,x_0)\mapsto (\omega{\cdot}t,u(t,\omega,x_0))\,.
\end{equation}
A $\tau$-minimal set $\mathcal{M}$ can be understood as the closure of the graph of a recurrent solution $\widetilde x(t)=u(t,\omega_0,\widetilde x(0))$ of one of the equations. Therefore, the classical problem of bifurcation of recurrent solutions ``around'' $\widetilde x(t)$ can be included in the analysis of bifurcation patterns for the family of equations
\begin{equation}\label{eq:62}
x'=f(\omega{\cdot}t,x)+\lambda(x-u(t,\omega,z))\,,
\end{equation}
for $(\omega,z)\in\mathcal{M}$. As explained in Subsection~\ref{subsec:sumslyap}, the families $\eqref{eq:62}_\lambda$ induce local skewproduct flows on the bundle $\mathcal{M}\times\mathbb{R}$, with base given by the global restricted flow $(\mathcal{M},\tau)$: we define $\hat f:\mathcal{M}\times\mathbb{R}\rightarrow\mathbb{R}$, $(\omega,z,x)\mapsto f(\omega,x)$ and $\beta:\mathcal{M}\rightarrow\mathbb{R}$, $(\omega,z)\mapsto z$ and rewrite $\eqref{eq:62}_\lambda$ as
$ x'=\hat f(\tau(t,\omega,z),x)+\lambda(x-\beta(\tau(t,\omega,z)))$ for $(\omega,z)\in\mathcal{M}$. The change of variables $y=x-\beta(\tau(t,\omega,z))$ takes this equation to
\begin{equation}\label{eq:new63}
y'=\hat g(\tau(t,\omega,z),y)+\lambda y\,,
\end{equation}
for $(\omega,z)\in\mathcal{M}$, where $\hat g(\omega,z,y)=\hat f(\omega,z,y+z)-\hat f(\omega,z,z)=f(\omega,y+z)-f(\omega,z)$. Clearly $\hat{\mathcal{M}}_0=\Omega\times\{0\}$ is a minimal set for the flow induced by $\eqref{eq:new63}_\lambda$ for all $\lambda\in\mathbb{R}$, coming from $\widetilde{\mathcal{M}}=\{(\omega,z,z)\colon\, (\omega,z)\in\mathcal{M}\}$, which is minimal for $\eqref{eq:62}_\lambda$.

In this way, studying bifurcations for $\eqref{eq:62}_\lambda$ ``around'' $\widetilde x(t)$ has been transformed in studying bifurcations for $\eqref{eq:new63}_\lambda$ ``around'' zero. To simplify the forthcoming discussion, we recover the standard notation of this work: the analysis of $\eqref{eq:new63}_\lambda$ is included in that of
\begin{equation}\label{eq:65}
x'=f(\omega{\cdot}t,x)+\lambda x\,,
\end{equation}
where $(\Omega,\sigma)$ is a minimal flow, $f\in C^{0,1}(\Omega\times\mathbb{R},\mathbb{R})$ satisfies $f(\omega,0)=0$ for every $\omega\in\Omega$, and $\lambda\in\mathbb{R}$. We will impose appropriate d-concavity and coercivity hypotheses when needed. We will analyze the bifurcation patterns of minimal sets and attractors which may arise for $\eqref{eq:65}_\lambda$ and, afterwards, translate the information to understand the possible bifurcation patterns for $\eqref{eq:62}_\lambda$: see Subsection~\ref{subsec:porponer}.

The solutions of $\eqref{eq:65}_\lambda$ define a local skewproduct flow
\begin{equation*}
\hat \tau_\lambda\colon\hat{\mathcal{U}}_\lambda\subseteq\mathbb{R}\times\Omega\times\mathbb{R}\rightarrow\Omega\times\mathbb{R},\quad (t,\omega,x_0)\mapsto (\omega{\cdot}t,\hat u_\lambda(t,\omega,x_0))\,.
\end{equation*}
Notice that $\hat{\mathcal{M}}_0=\Omega\times\{0\}$ is always a $\hat\tau_\lambda$-minimal set.
\begin{proposition}\label{prop:alphasubeq2bifurcation}The following statements hold:
\begin{enumerate}[label=\rm{(\roman*)}]
\item Any strictly positive (resp.~negative) global upper solution of $\eqref{eq:65}_\lambda$ is a strict global upper solution of $\eqref{eq:65}_\xi$ whenever $\xi<\lambda$ (resp.~$\lambda<\xi$). In particular, any strictly positive (resp.~negative) equilibrium for $\eqref{eq:65}_\lambda$ is a strong superequilibrium and a strong time-reversed subequilibrium for $\eqref{eq:65}_\xi$ whenever $\xi<\lambda$ (resp.~$\lambda<\xi$).
\item Any strictly positive (resp.~negative) global lower solution of $\eqref{eq:65}_\lambda$ is a strict global lower solution of $\eqref{eq:65}_\xi$ whenever $\lambda<\xi$ (resp.~$\xi<\lambda$). In particular, any strictly positive (resp.~negative) equilibrium for $\eqref{eq:65}_\lambda$ is a strong subequilibrium and a strong time-reversed superequilibrium for $\eqref{eq:65}_\xi$ whenever $\lambda<\xi$ (resp.~$\xi<\lambda$).
\end{enumerate}
\end{proposition}
\begin{proof}
The statements follow from the properties described in Subsection~\ref{subsec:equilibria}.
\end{proof}
The information provided by Theorem~\ref{th:attractorexistence} plays a role in the next statement.
\begin{proposition} \label{prop:alphas} Assume that $f\in C^{0,1}(\Omega\times\mathbb{R},\mathbb{R})$ is also $\mathrm{(Co)}$ and let
\begin{equation}\label{eq:newattractor}
\hat{\mathcal{A}}_\lambda=\bigcup_{\omega\in\Omega}\big(\{\omega\}\times[\hat{\alpha}_\lambda(\omega),\hat{\beta}_\lambda(\omega)]\big)
\end{equation}
be the global attractor of $\eqref{eq:65}_\lambda$. Then,
\begin{enumerate}[label=\rm{(\roman*)}]
\item $\hat\alpha_\lambda(\omega)\leq 0\leq\hat\beta_\lambda(\omega)$ for every $\omega\in\Omega$ and $\lambda\in\mathbb{R}$.
\item For every $\omega\in\Omega$, the maps $\lambda\mapsto \hat\beta_\lambda(\omega)$ and $\lambda\mapsto\hat\alpha_\lambda(\omega)$ are respectively nondecreasing and nonincreasing on $\mathbb{R}$ and both are right-continuous.
\item $\lim_{\lambda\rightarrow\infty} \hat\alpha_\lambda(\omega)=-\infty$ and $\lim_{\lambda\rightarrow\infty} \hat\beta_\lambda(\omega)=\infty$ uniformly on $\Omega$. In particular, $\hat{\tau}_\lambda$ admits at least three minimal sets for $\lambda$ large enough.
\item There exists $\lambda_0\in\mathbb{R}$ such that $\hat{\mathcal{A}}_\lambda=\hat{\mathcal{M}}_0=\Omega\times\{0\}$ for every $\lambda<\lambda_0$ and it is hyperbolic attractive.
\end{enumerate}
\end{proposition}
\begin{proof} Property (i) follows from $\hat{\mathcal{M}}_0=\Omega\times\{0\}\subseteq\hat{\mathcal{A}}_\lambda$. The arguments used to prove Theorem~\ref{th:attractorproperties}(i)-(ii) can be used to check (ii) and the first assertion in (iii). The second one follows from the first one, Proposition~\ref{prop:precollision} and the minimality of $\hat{\mathcal{M}}_0$. To prove (iv), note that (ii) ensures that $\mathcal{A}_\lambda\subseteq\mathcal{A}_\xi$ if $\lambda<\xi$. Let us fix $\xi\in\mathbb{R}$, take $r>0$ such that $\mathcal{A}_\xi\subseteq\Omega\times[-r,r]$ and define $\lambda_0=\min\{\xi,\,\inf\{-f_x(\omega,x)\colon\,(\omega,x)\in\Omega\times[-r,r]\}\}$. Then, $f_x(\omega,x)+\lambda<0$ for all $\lambda<\lambda_0$ and $(\omega,x)\in\Omega\times[-r,r]$, so $\gamma(\mathcal{A}_\lambda,\nu)<0$ for all $\nu\in\mathfrak{M}_\mathrm{erg}(\mathcal{A}_\lambda,\tau)$. Theorem~3.4 of \cite{lineardissipativescalar} guarantees that $\mathcal{A}_\lambda$ is hyperbolic attractive and coincides with $\Omega\times\{0\}$.
\end{proof}
Now, we will study the bifurcation problem $\eqref{eq:65}_\lambda$ in terms of two parameters:
\begin{equation*}
\begin{split}
\mu_-=\inf\{\lambda\colon\text{ the graph of }\hat\alpha_\xi\text{ is a hyperbolic minimal set }\hat{\mathcal{M}}_\xi^l<\hat{\mathcal{M}}_0,\;\forall\xi>\lambda\}\,,\\
\mu_+=\inf\{\lambda\colon\text{ the graph of }\hat\beta_\xi\text{ is a hyperbolic minimal set }\hat{\mathcal{M}}_\xi^u>\hat{\mathcal{M}}_0,\;\forall\xi>\lambda\}\,.\\
\end{split}
\end{equation*}
\begin{theorem}\label{th:pitchfork} Let $f\in C^{0,2}(\Omega\times\mathbb{R},\mathbb{R})$ be $\mathrm{(Co)}$ and $\mathrm{(SDC)}$, and let $[-\lambda_+,-\lambda_-]$ be the Sacker and Sell spectrum of $f_x$ on $\hat{\mathcal{M}}_0=\Omega\times\{0\}$, with $\lambda_-\leq\lambda_+$. Then,
\begin{enumerate}[label=\rm{(\roman*)}]
\item $\hat{\mathcal{M}}_0$ is hyperbolic attractive (resp.~repulsive) if $\lambda<\lambda_-$ (resp.~$\lambda>\lambda_+$) and nonhyperbolic if $\lambda\in[\lambda_-,\lambda_+]$. Both $\hat{\mathcal{M}}_\lambda^l$ for $\lambda>\mu_-$ and $\hat{\mathcal{M}}_\lambda^u$ for $\lambda>\mu_+$ are hyperbolic attractive.
\item $\mu_-,\mu_+\in(-\infty,\lambda_+]$, and either $\mu_+=\lambda_+$ (which happens if $\hat\beta_\lambda$ collides with $0$ on a residual $\sigma$-invariant set as $\lambda\downarrow\lambda_+$) or $\mu_-=\lambda_+$ (which happens if $\hat\alpha_\lambda$ collides with $0$ on a residual $\sigma$-invariant set as $\lambda\downarrow\lambda_+$).
\item {\rm(Global pitchfork bifurcation).} If $\mu_-=\mu_+=\lambda_+$, then a global pitchfork bifurcation pattern of minimal sets arises around $\hat{\mathcal{M}}_0$ at $\lambda_+$: $\hat\tau_\lambda$ admits the three different hyperbolic minimal sets $\hat{\mathcal{M}}_\lambda^l<\hat{\mathcal{M}}_0<\hat{\mathcal{M}}_\lambda^u$ for $\lambda>\lambda_+$; both $\hat\alpha_\lambda$ and $\hat\beta_\lambda$ collide with $0$ on a residual $\sigma$-invariant set as $\lambda\downarrow\lambda_+$; $\hat{\mathcal{M}}_0$ is the unique $\hat\tau_\lambda$-minimal set if $\lambda\leq\lambda_+$; and $\hat{\mathcal{A}}_\lambda=\hat{\mathcal{M}}_0$ if $\lambda<\lambda_-$.
\end{enumerate}
\end{theorem}
\begin{proof}(i) The spectrum of the linearized equation $f_x+\lambda$ on $\hat{\mathcal{M}}_0$ is $[\lambda-\lambda_+,\lambda-\lambda_-]$, which yields the stated hyperbolicity properties for $\hat{\mathcal{M}}_0$ (see Section~\ref{subsec:furstenbergremark}). Proposition~\ref{prop:jagerinmeasure} ensures the attractiveness of the hyperbolic sets $\hat{\mathcal{M}}_\lambda^l$ and $\hat{\mathcal{M}}_\lambda^u$.

(ii) Since $\hat{\mathcal{M}}_0$ is hyperbolic repulsive if and only if $\lambda>\lambda_+$, Proposition~\ref{prop:hyperbolicpositiveminimal}(i) and Theorem~\ref{th:new}(ii) ensure that $\lambda_+=\inf\{\lambda\in\mathbb{R}\colon\,\hat\tau_\xi$ admits three different hyperbolic minimal sets for every $\xi\in(\lambda,\infty)\}$. Proposition~\ref{prop:precollision} shows that these three minimal sets are $\hat{\mathcal{M}}_\lambda^l<\hat{\mathcal{M}}_0<\hat{\mathcal{M}}_\lambda^u$, respectively given by the graphs of $\hat\alpha_\lambda<0<\hat\beta_\lambda$. Consequently, $\mu_+,\mu_-\leq\lambda_+$. Proposition~\ref{prop:alphas}(iv) ensures that $\mu_+$ and $\mu_-$ are finite. Proposition~\ref{prop:alphas}(i)-(ii) ensure that $\hat\alpha_{\lambda_+}(\omega)=\lim_{\lambda\downarrow\lambda_+}\hat\alpha_\lambda(\omega)\leq 0$ and $\hat\beta_{\lambda_+}(\omega)=\lim_{\lambda\downarrow\lambda_+}\hat\beta_\lambda(\omega)\geq0$. Let $\mathcal{R}$ be the residual set of common continuity points for $\hat\alpha_{\lambda_+}$ and $\hat\beta_{\lambda_+}$. Assume first that $\hat\alpha_{\lambda_+}(\omega)<0<\hat\beta_{\lambda_+}(\omega)$ for all $\omega\in\Omega$. Proposition~\ref{prop:collision} ensures that, in this case, there exist three different minimal sets at $\lambda_+$ and Theorem~\ref{th:new} guarantees that they are hyperbolic. Hence, the nonhyperbolicity of $\hat{\mathcal{M}}_0$ at $\lambda_+$ precludes this case. Consequently, there exists $\omega_0\in\Omega$ such that either $\hat\alpha_{\lambda_+}(\omega_0)=0$ or $\hat\beta_{\lambda_+}(\omega_0)=0$, which according to Propositions~\ref{prop:postcollision} and \ref{prop:collision} ensures that either $\hat\alpha_{\lambda_+}$ or $\hat\beta_{\lambda_+}$ coincides with $0$ on $\mathcal{R}$ and hence that either $\mu_-=\lambda_+$ or $\mu_+=\lambda_+$.

(iii) Let us see that $\hat\alpha_+(\omega)=\hat\beta_+(\omega)=0$ for all $\omega\in\mathcal{R}$ if $\mu_-=\mu_+=\lambda_+$. We proceed by contradiction: Proposition~\ref{prop:collision} ensures that either $\hat\alpha_{\lambda_+}(\omega)<0$ for all $\omega\in\mathcal{R}$ or $\hat\beta_{\lambda_+}(\omega)>0$ for all $\omega\in\mathcal{R}$. Proposition~\ref{prop:precollision} ensures that there exists another $\hat\tau_{\lambda_+}$-minimal set distinct from $\hat{\mathcal{M}}_0$ and Theorem~\ref{th:twominimalsets}(i) ensures that it is hyperbolic attractive (the spectrum of $f_x$ on $\hat{\mathcal{M}}_0$ for $\hat\tau_{\lambda_+}$ is contained in $[0,\infty)$). Therefore, it has a hyperbolic continuation for $\lambda<\lambda_+$ close enough, contradicting the definition of either $\mu_-$ or $\mu_+$. The assertion follows. Now, Proposition~\ref{prop:alphas}(ii) ensures that $\hat{\mathcal{A}}_\lambda$ is pinched for every $\lambda<\lambda_+$, as $\hat{\mathcal{A}}_\lambda\subseteq\hat{\mathcal{A}}_{\lambda_+}$, and hence $\hat{\mathcal{M}}_0$ is the unique $\hat\tau_\lambda$-minimal set; since $\hat{\mathcal{M}}_0$ is hyperbolic attractive for $\lambda<\lambda_+$, Theorem 3.4 of \cite{lineardissipativescalar} ensures that $\hat{\mathcal{A}}_\lambda=\hat{\mathcal{M}}_0$ for $\lambda<\lambda_+$.
\end{proof}
Recall that Theorem~\ref{th:pitchfork}(ii) ensures $\mu_-=\lambda_+$ or $\mu_+=\lambda_+$ (or both).
\begin{theorem}[Local saddle-node and transcritical bifurcations]\label{th:transcritical} Let $f\in C^{0,2}(\Omega\times\mathbb{R},\mathbb{R})$ be $\mathrm{(Co)}$ and $\mathrm{(SDC)}$, and let $[-\lambda_+,-\lambda_-]$ be the Sacker and Sell spectrum of $f_x$ on $\hat{\mathcal{M}}_0=\Omega\times\{0\}$, with $\lambda_-\leq\lambda_+$. Assume that $\mu_+=\lambda_+$ and $\mu_-<\lambda_-$. Then,
\begin{enumerate}[label=\rm{(\roman*)}]
\item $\hat\tau_\lambda$ admits exactly three minimal sets $\hat{\mathcal{M}}_\lambda^l<\hat{\mathcal{M}}_0<\hat{\mathcal{M}}_\lambda^u$ for $\lambda>\lambda_+$, and $\hat\beta_\lambda$ collides with $0$ on a residual $\sigma$-invariant set as $\lambda\downarrow\lambda_+$,
\item $\hat\tau_\lambda$ admits exactly two minimal sets $\hat{\mathcal{M}}_\lambda^l<\hat{\mathcal{M}}_0$ for $\lambda\in[\lambda_-,\lambda_+]$,
\item $\hat\tau_\lambda$ admits three hyperbolic minimal sets $\hat{\mathcal{M}}_\lambda^l<\hat{\mathcal{N}}_\lambda<\hat{\mathcal{M}}_0$ if $\lambda\in(\mu_-,\lambda_-)$, where $\hat{\mathcal{N}}_\lambda$ is hyperbolic repulsive and given by the graph of a continuous map $\hat\kappa_\lambda\colon\Omega\rightarrow\mathbb{R}$ which increases strictly as $\lambda$ increases in $(\mu_-,\lambda_-)$, and which collides with $\hat\alpha_\lambda$ (resp.~with $0$) on a residual set as $\lambda\downarrow\mu_-$ (resp.~$\lambda\uparrow\lambda_-$),
\item $\hat\tau_{\mu_-}$ admits exactly two different minimal sets $\hat{\mathcal{M}}_{\mu_-}^l<\hat{\mathcal{M}}_0$, with $\hat{\mathcal{M}}_{\mu_-}^l$ nonhyperbolic,
\item $\hat{\mathcal{A}}_\lambda=\hat{\mathcal{M}}_0$ for $\lambda<\mu_-$.
\end{enumerate}
In particular, $\mu_-$, $\lambda_-$ and $\lambda_+$ are the unique bifurcation points: a local saddle-node bifurcation of minimal sets occurs around $\hat{\mathcal{M}}_{\mu_-}$ at $\mu_-$, as well as a discontinuous bifurcation of attractors; and the bifurcation pattern which arises around $\hat{\mathcal{M}}_0$ on $[\lambda_-,\lambda_+]$ can be understood as a generalized local transcritical bifurcation of minimal sets around $\hat{\mathcal{M}}_0$ through the interval $[\lambda_-,\lambda_+]$ (which is classic for $\lambda_-=\lambda_+$).

The possibilities for the global bifurcation diagram are symmetric to the described ones with respect to the horizontal axis if $\mu_-=\lambda_+$ and $\mu_+<\lambda_-$.
\end{theorem}
\begin{proof} (i) The results in (i) follow from Theorem~\ref{th:pitchfork}(i)-(ii).

(ii) The nonhyperbolicity of $\hat{\mathcal{M}}_0$ ensured by Theorem~\ref{th:pitchfork}(i) and Theorem~\ref{th:new} preclude the existence of other minimal set apart from $\hat{\mathcal{M}}_0$ and $\hat{\mathcal{M}}_\lambda^l$ for $\lambda\in[\lambda_-,\lambda_+]$.

(iii) As both $\hat{\mathcal{M}}_0$ and $\hat{\mathcal{M}}^l_\lambda$ are hyperbolic attractive for $\lambda\in(\mu_-,\lambda_-)$, Proposition~\ref{prop:hyperbolicpositiveminimal} and Theorem~\ref{th:new} ensure the existence of an intermediate repulsive hyperbolic minimal set $\hat{\mathcal{N}}_{\lambda}$ for $\lambda\in(\mu_-,\lambda_-)$ given by the graph of a continuous equilibrium $\hat\kappa_{\lambda}\colon\Omega\rightarrow\mathbb{R}$. Reasoning as in Theorem~\ref{th:bifdiagram}, but taking into account that in this case the kind of monotonicity depends on the half-plane in which we are working (and we are in the negative one), we can check that $\lambda\mapsto \hat\kappa_\lambda(\omega)$ is strictly increasing on $(\mu_-,\lambda_-)$ for all $\omega\in\Omega$. In consequence, $\hat\kappa_{\lambda_-}(\omega)=\lim_{\lambda\uparrow\lambda_-}\hat\kappa_\lambda(\omega)$ defines a lower semicontinuous $\hat\tau_{\lambda_-}$-equilibrium. Since $\hat\kappa_{\lambda_-}>\hat\kappa_\lambda>\hat\alpha_\lambda>\hat\alpha_{\lambda_-}$ for $\lambda\in(\mu_-,\lambda_-)$ and $\hat{\mathcal{M}}_0$ and $\hat{\mathcal{M}}_{\lambda_-}^l$ are the unique $\tau_{\lambda_-}$-minimal sets, Propositions~\ref{prop:precollision} and \ref{prop:collision} ensure that the minimal set \eqref{eq:precollision} provided by $\hat\kappa_{\lambda_-}$ is $\hat{\mathcal{M}}_0$. That is, $\hat\kappa_{\lambda_-}$ vanishes at the residual set of its continuity points. Combining this information with that provided by Theorem~\ref{th:pitchfork}(ii), we observe that: two hyperbolic minimal sets for $\lambda<\lambda_-$ close to $\lambda_-$ collide (at a residual set of base points) at $\lambda_-$ and two hyperbolic minimal sets for $\lambda>\lambda_+$ close to $\lambda_+$ collide at $\lambda_+$, giving rise to a unique nonhyperbolic minimal set on the interval $[\lambda_-,\lambda_+]$. This is the generalized local transcritical bifurcation of minimal sets around $\hat{\mathcal{M}}_0$ mentioned in the statement.

(iv) Let us define $\hat\kappa_{\mu_-}(\omega)=\lim_{\lambda\downarrow\mu_-} \hat\kappa_\lambda(\omega)$, and let us recall that $\hat\alpha_{\mu_-}(\omega)=\lim_{\lambda\downarrow\mu_-} \hat\alpha_\lambda(\omega)$ (see Proposition~\ref{prop:alphas}(ii)). Then, $\hat\kappa_{\mu_-}(\omega)\geq\hat\alpha_{\mu_-}(\omega)$ for all $\omega\in\Omega$. It follows from the definition of $\mu_-$ that the minimal set $\hat{\mathcal{M}}_{\mu_-}^l$ determined from $\hat\alpha_{\mu_-}$ by \eqref{eq:precollision} is nonhyperbolic, and hence (by Theorem~\ref{th:new}) that it coincides with that determined from $\hat\kappa_{\mu_-}$. That is, $\hat\kappa_{\mu_-}$ and $\hat\alpha_{\mu_-}$ coincide at the residual set of their common continuity points, giving rise to a nonhyperbolic $\hat\tau_{\mu_-}$-minimal set $\hat{\mathcal{M}}_{\mu_-}^l$. In addition, $\hat{\mathcal{M}}_{\mu_-}^l<\hat{\mathcal{N}}_\lambda<\hat{\mathcal{M}}_0$ for $\lambda\in(\mu_-,\lambda_-)$; and $\hat{\mathcal{M}}_{\mu_-}^l$ and $\hat{\mathcal{M}}_0$ are the unique $\hat\tau_{\mu_-}$-minimal sets, since $\hat{\mathcal{M}}_{\mu_-}^l$ is nonhyperbolic (see again Theorem~\ref{th:new}(ii)).

(v) We must fix $\lambda<\mu_-$ and check that $\hat{\mathcal{A}}_\lambda=\hat{\mathcal{M}}_0$. Since $\hat{\mathcal{M}}_0$ is an attractive hyperbolic copy of the base for $\hat\tau_\lambda$, it suffices to check that it is the unique $\hat\tau_\lambda$-minimal set and apply Theorem~3.4 of \cite{lineardissipativescalar}. Note first that $\hat\beta_\lambda\equiv0$ for $\lambda\leq\mu_-$, since $\hat\beta_\lambda\equiv0$ for $\lambda\in(\mu_-,\lambda_-)$ and $\lambda\mapsto\hat\beta_\lambda$ is nondecreasing, and hence $\hat{\mathcal{M}}_0$ is the upper $\hat\tau_\lambda$-minimal set. Let us take another one $\hat{\mathcal{N}}_\lambda\leq \hat{\mathcal{M}}_0$, and prove that $\hat{\mathcal{N}}_\lambda=\hat{\mathcal{M}}_0$. We fix $(\omega_0,x_0)\in\hat{\mathcal{N}}_\lambda$. Since $\hat\alpha_{\mu_-}$ and $\hat\kappa_{\mu_-}$ are strong $\hat\tau_\lambda$-subequilibria which coincide at their continuity points, Proposition~\ref{prop:atravesar} provides $s>0$ and $e>0$ such that $\hat\kappa_{\mu_-}(\omega_0)+e<\hat u_\lambda(s,\omega_0{\cdot}(-s),\hat\alpha_{\mu_-}(\omega_0{\cdot}(-s)))\leq\hat u_\lambda(s,\omega_0{\cdot}(-s),\hat\alpha_{\lambda}(\omega_0{\cdot}(-s)))=\hat\alpha_\lambda(\omega_0)\leq x_0$. Since $\hat\kappa_{\mu_-}(\omega)=\lim_{\lambda\downarrow\mu_-} \hat\kappa_\lambda(\omega)$, there exists $\xi_1>\mu_-$ such that $\hat\kappa_{\xi_1}(\omega_0)<x_0$. We take any $\xi\in(\xi_1,\lambda_-)$ and apply Proposition~\ref{prop:atravesar2} to the family of strong $\hat\tau_\lambda$-subequilibria $\hat\kappa_\mu$ with $\mu\in[\xi_1,\xi]$ to conclude that there exist $s_\xi>0$ and $e_\xi>0$ such that $\hat\kappa_{\xi}(\omega_0{\cdot}s)+e_\xi\leq\hat u_\lambda(s,\omega_0,\hat\kappa_{\xi_1}(\omega_0))<\hat u_\lambda(s,\omega_0,x_0)$ for all $s\geq s_\xi$. The lower semicontinuity of $\hat\kappa_\xi$ ensures that $\hat{\mathcal{N}}_\lambda$, which is the $\boldsymbol\upomega$-limit set of $(\omega_0,x_0)$, is strongly above the graph of $\hat\kappa_\xi$ for all $\xi\in[\xi_1,\lambda_-)$. So, it is above the graph of $\hat\kappa_{\lambda_-}$. Hence, $(\hat{\mathcal{N}}_\lambda)_\omega=\{0\}$ for all the points $\omega$ of the residual set at which $\kappa_{\lambda_-}$ coincides with $0$, which yields $\hat{\mathcal{N}}_\lambda=\hat{\mathcal{M}}_0$. This completes the proof of the theorem. Note that a local saddle-node bifurcation of minimal sets occurs at $\mu_-$ around $\hat{\mathcal{M}}_{\mu_-}^l$, due to the collision of $\hat{\mathcal{M}}_\lambda^l$ and $\hat{\mathcal{N}}_\lambda$ as $\lambda\downarrow\mu_-$. Also, since $\lambda\mapsto \alpha_\lambda(\omega)$ is discontinuous at $\mu_-$ for all $\omega\in\Omega$, there is a discontinuous bifurcation of attractors at $\mu_-$.
\end{proof}
Note that, if $\lambda_-=\lambda_+$, then $[\lambda_-,\lambda_+]$ reduces to a point, so Theorems~\ref{th:pitchfork}(iii) and \ref{th:transcritical} describe the only possible bifurcation diagrams of minimal sets of the problem $\eqref{eq:65}_\lambda$. These diagrams were found in Theorem~5.7 of \cite{nunezobaya} for the uniquely ergodic case. Moreover, in this case, the number of $\hat\tau_{\lambda_-}$-minimal sets characterizes what bifurcation diagram takes place: there is either one (and we are in the pitchfork case) or two (and we are in the transcritical and saddle-node case). The complete description of the ``generalized pitchfork bifurcation of minimal sets through the interval $[\lambda_-,\lambda_+]$'' with $\lambda_-<\lambda_+$ is a problem under study. Note also that, in all the cases that we have described, we find a compact $\tau_{\lambda_+}$-invariant pinched set containing $\hat{\mathcal{M}}_0$ which is not, in general, a copy of the base.


We complete this subsection by describing a simple situation in which the bifurcation diagram for $\eqref{eq:65}_\lambda$ is that of Theorem 6.3: a global pitchfork bifurcation occurs at the unique bifurcation point.
\begin{proposition}
Let $f\in C^{0,2}(\Omega\times\mathbb{R},\mathbb{R})$ be $\mathrm{(Co)}$ and $\mathrm{(DC)}$, and let $[-\lambda_+,-\lambda_-]$ be the Sacker and Sell spectrum of $f_x$ on $\hat{\mathcal{M}}_0 =\Omega\times\{0\}$, with $\lambda_-\leq\lambda_+$. Then,
\begin{itemize}
\item[\rm(i)] if $f_{xx}(\omega,0)\geq 0$ (resp.~$f_{xx}(\omega,0)\leq 0$) for all $\omega\in\Omega$, then $\alpha_\lambda$ (resp.~$\beta_\lambda$) takes the value $0$ at its continuity points for all $\lambda<\lambda_+$.
\item[\rm(ii)] If $f$ is $\mathrm{(SDC)}$ and $f_{xx}(\omega,0)=0$ for all $\omega\in\Omega$, then the bifurcation diagram is that described in Theorem~$\mathrm{\ref{th:pitchfork}(iii)}$.
\end{itemize}
\end{proposition}
\begin{proof}
(i) If $f_{xx}(\omega,0)\ge 0$, then $f_{xx}(\omega,x)\ge 0$ for all $x\le 0$. We assume for contradiction the existence of $\lambda<\lambda_+$ such that $\alpha_\lambda(\omega)<0$ at its continuity points, which is equivalent to say that $\alpha_\lambda$ is a strictly negative equilibrium (see Proposition~\ref{prop:postcollision}). Taylor's Theorem ensures that $\alpha_\lambda'(\omega{\cdot}t)/\alpha_\lambda(\omega{\cdot}t)\le f_x(\omega{\cdot}t,0)+\lambda$, and hence Birkhoff's Ergodic Theorem ensures that $0=\int_{\Omega}(\alpha_\lambda'(\omega)/\alpha_\lambda(\omega))\,dm \le\int_\Omega (f_x(\omega,0)+\lambda)\,dm$ for all $m\in \mathfrak{M}_{\rm erg}(\Omega,\sigma)$. But this contradicts the existence of $m_+$ with $\int_\Omega f_x(\omega,0)\,dm_+=-\lambda_+<-\lambda$. The proof for $\beta_\lambda$ is analogous.

(ii) Property (i) shows that $\mu_-=\mu_+=\lambda_+$ for the values $\mu_-,\,\mu_+$ defined before Theorem 6.3, whose point (iii) shows the assertion.
\end{proof}
\subsection{Transcritical bifurcation in a simple example} In what follows, we will study the bifurcation diagram for the one-parametric family
\begin{equation}\label{eq:simplexample}
x'=-x^3+(a_2(\omega{\cdot}t)+\xi)x^2\,,\quad\omega\in\Omega\,,
\end{equation}
where $a_2\in C(\Omega,\mathbb{R})$, $x\in\mathbb{R}$, and $\xi\in\mathbb{R}$. We will denote by $\check\tau_\xi$ the local skewproduct flow defined by $\eqref{eq:simplexample}_\xi$, given by the solutions $\check u_\xi(t,\omega,x_0)$. The study of such example will allow us to characterize the bifurcation diagrams of minimal sets of $x'=-x^3+a_2(\omega{\cdot}t)x^2+\lambda x$ in terms of the Sacker and Sell spectrum of $a_2$. Note that the right-hand side of $\eqref{eq:simplexample}_\xi$ is $\mathrm{(Co)}$ and $\mathrm{(SDC)}_*$. The coercivity property and Theorem~\ref{th:attractorexistence} ensure the existence of global attractor $\check{\mathcal{A}}_\xi$ for $\check\tau_\xi$, delimited by semicontinuous equilibria $\check\alpha_\xi$ and $\check\beta_\xi$, and it is not hard to check that their properties are similar to those described in Theorem~\ref{th:attractorproperties}, i.e., the maps $\xi\mapsto\check\beta_\xi(\omega)$ and $\xi\mapsto\check\alpha_\xi(\omega)$ are nondecreasing for all $\omega\in\Omega$ and, respectively, right- and left-continuous. The Sacker and Sell spectrum of $a\colon\Omega\mapsto\mathbb{R}$ is defined in Subsection~\ref{subsec:furstenbergremark}.
\begin{proposition}[Weak generalized transcritical bifurcation]\label{prop:newdiagram} Let $[-\xi_+,-\xi_-]$ be the Sacker and Sell spectrum of $a_2$, with $\xi_-\leq\xi_+$. Then, $\check{\mathcal{M}}_0=\Omega\times\{0\}$ is a nonhyperbolic $\check\tau_\xi$-minimal set for all $\xi\in\mathbb{R}$. In addition,
\begin{enumerate}[label=\rm{(\roman*)}]
\item $\check\tau_\xi$ admits exactly two minimal sets $\check{\mathcal{M}}_\xi^l<\check{\mathcal{M}}_0$ given by the graphs of $\check\alpha_\xi<0$ for $\xi<\xi_-$, and with $\check{\mathcal{M}}_\xi^l$ hyperbolic attractive; and $\check\alpha_\xi$ collides with $0$ on a $\sigma$-invariant residual set as $\xi\uparrow\xi_-$;
\item $\check{\mathcal{M}}_0$ is the unique $\check\tau_\xi$-minimal set for $\xi\in[\xi_-,\xi_+]$;
\item $\check\tau_\xi$ admits exactly two minimal sets $\check{\mathcal{M}}_0<\check{\mathcal{M}}_\xi^u$ given by the graphs of $0<\check\beta_\xi$ for $\xi>\xi_+$, and with $\check{\mathcal{M}}_\xi^u$ hyperbolic attractive; and $\check\beta_\xi$ collides with $0$ in a $\sigma$-invariant residual set as $\xi\downarrow\xi_+$.
\end{enumerate}
\end{proposition}
\begin{proof} Let us fix $\xi<\xi_+$. We will prove the existence of $\omega_0\in\Omega$ such that $\check u_\xi(t,\omega_0,x)$ is unbounded for any $x>0$, which ensures the absence of $\check\tau_\xi$-minimal sets above $\check{\mathcal{M}}_0$. Let $\check m_+\in\mathfrak{M}_\mathrm{erg}(\Omega,\sigma)$ satisfy $\int_\Omega (a_2(\omega)+\xi_+)\,d\check m_+=0$. Birkhoff's Ergodic Theorem provides $\Omega_0\subseteq\Omega$ with $\check m_+(\Omega_0)=1$ such that $\sup_{t\leq 0}\int_0^t (a_2(\omega{\cdot}s)+\xi)\, ds=\infty$. We take $\omega_0\in\Omega_0$. Let $\check v_\xi(t,\omega_0,x_0)$ solve $x'=(a_2(\omega_0{\cdot}t)+\xi)x^2$ with $\check v_\xi(0,\omega_0,x_0)=x_0>0$. That is, $\check v_\xi(t,\omega_0,x_0)=(1/x_0-\int_0^t(a_2(\omega_0{\cdot}s)+\xi)\,ds)^{-1}$, which means that $\check v_\xi(t,\omega_0,x_0)$ tends to $\infty$ as $t$ decreases, in finite time. It is easy to check that $\check u_\xi(t,\omega_0,x_0)>\check v_\xi(t,\omega_0,x_0)$ if $t<0$, from where the initial assertion follows. Analogously, it can be proved that there are no $\check\tau_\xi$-minimal sets below $\check{\mathcal{M}}_0$ for $\xi>\xi_-$. In particular, $\check{\mathcal{M}}_0$ is the unique $\check\tau_\xi$-minimal set for $\xi\in(\xi_-,\xi_+)$.

Let us now fix $\xi>\xi_+$ and prove the existence of a $\check\tau_\xi$-minimal set $\check{\mathcal{M}}_\xi^u$ strictly above $\check{\mathcal{M}}_0$. We will use this property: given $\alpha\in(0,\xi-\xi_+)$, there exists $t_\alpha>0$ such that $\int_0^{t_\alpha}(a_2(\omega{\cdot}s)+\xi)\, ds\geq\alpha t_\alpha$ for all $\omega\in\Omega$. To prove it, assume for contradiction the existence of $\{t_n\}_{n\in\mathbb{N}}\uparrow\infty$ and $\{\omega_n\}_{n\in\mathbb{N}}$ such that $(1/t_n)\int_0^{t_n}(a_2(\omega_n{\cdot}s)+\xi)\, ds<\alpha$, and use an argument analogous to that of Kryloff and Bogoliuboff's Theorem (see Theorem~9.05 of \cite{nemytskii1}) to construct $m\in\mathfrak{M}_\mathrm{inv}(\Omega,\sigma)$ such that $\int_\Omega (a_2(\omega)+\xi)\,dm\leq \alpha<\xi-\xi_+$, which is impossible (see Subsection~\ref{subsec:furstenbergremark}).

Let us take $\alpha\in(0,\xi-\xi_+)$ with $\alpha<1$. The coercivity of $-x^3+(a_2(\omega)+\xi)x^2$ ensures that $\check u_\xi(t,\omega_0,1/\alpha)$ is globally forward defined and bounded. We will check below that $\check u_\xi(t,\omega_0,1/\alpha)$ is also bounded away from zero for $t\geq0$, which shows the existence of a minimal set $\check{\mathcal{M}}_\xi^u>\check{\mathcal{M}}_0$ contained in the $\boldsymbol\upomega$-limit of $(\omega_0,1/\alpha)$. The function $w(t)=(\check u_\xi(t,\omega_0,1/\alpha))^{-1}$ solves $y'=-(a_2(\omega{\cdot}t)+\xi)+1/y$ with $w(0)=\alpha<2/\alpha$. Let us call $t_1=\sup\{t>0\colon\, w(s)\leq 2/\alpha+lt_\alpha$ for all $s\in[0,t]\}$, where $l=1/\alpha+\sup_{\omega\in\Omega}|a_2(\omega)+\xi|$. We assume for contradiction that $t_1<\infty$, and define $t_0=\inf\{t<t_1\colon\, w(s)\geq 2/\alpha$ for $s\in[t,t_1]\}$. Then, $t_0\leq t_1-t_\alpha$: otherwise, $w(t_1)= w(t_0)+\int_{t_0}^{t_1} (-(a_2(\omega{\cdot}s)+\xi)+1/w(s))\, ds<2/\alpha+lt_\alpha$, which is not the case. In particular, $w(t)\geq 2/\alpha$ for $t\in[t_1-t_\alpha,t_1]$, and hence
\begin{equation*}
\begin{split}
w(t_1)&=w(t_1-t_\alpha)-\int_0^{t_\alpha}(a_2(\omega_0{\cdot}(t_1-t_\alpha){\cdot}s)+\xi)\, ds+\int_{t_1-t_\alpha}^{t_1}\frac{1}{w(s)}\,ds\\
&\leq w(t_1-t_\alpha)+(-\alpha+\alpha/2)t_\alpha<w(t_1-t_\alpha)\,,
\end{split}
\end{equation*}
which contradicts the definition of $t_1$. This shows that $\check u_\xi(t,\omega_0,1/\alpha)\geq (2/\alpha+lt_\alpha)^{-1}$, which completes this step.

We have proved the existence of $\check{\mathcal{M}}_\xi^u>\check{\mathcal{M}}_0$ for $\xi>\xi_+$. Analogous arguments show the existence of $\check{\mathcal{M}}_\xi^l<\check{\mathcal{M}}_0$ for $\xi<\xi_-$. Theorem~\ref{th:new} and the nonhyperbolicity of $\check{\mathcal{M}}_0$ ensure that $\check\tau_\xi$ admits at most two minimal sets. Hence, since the Sacker and Sell spectrum on $\check{\mathcal{M}}_0$ reduces to $\{0\}$, Theorem~\ref{th:twominimalsets} ensures that $\check{\mathcal{M}}_\xi^u$ (for $\xi>\xi_+$) and $\check{\mathcal{M}}_\xi^l$ (for $\xi<\xi_-$) are hyperbolic attractive. Therefore, they are copies of the base: the graphs of $\check\beta_\xi$ and $\check\alpha_\xi$, respectively. Theorem~\ref{th:twominimalsets} also precludes the existence of a second minimal set for $\check\tau_{\xi_-}$ (resp. $\check\tau_{\xi_+}$), since it would be hyperbolic and hence persisting for $\xi$ close to $\xi_+$ (resp. $\xi_-$). In turn, this means that $\check\beta_{\xi_+}(\omega)=\lim_{\xi\downarrow\xi_+}\check\beta_\xi(\omega)$ (resp. $\check\alpha_{\xi_-}(\omega)=\lim_{\xi\uparrow\xi_-}\check\alpha_\xi(\omega)$) coincides with 0 at their continuity points: otherwise, Propositions~\ref{prop:postcollision} and \ref{prop:precollision} would provide a second $\check\tau_{\xi_+}$ (resp. $\check\tau_{\xi_-}$) minimal set. This completes the proof of all the assertions.
\end{proof}
Our last result in this subsection is a remarkable consequence of the previous one: the bifurcation diagram of
\begin{equation}\label{eq:laultima}
x'=-x^3+a_2(\omega{\cdot}t)x^2+\lambda x\,,\quad \omega\in\Omega\,,
\end{equation}
just depends on the Sacker and Sell spectrum of $a_2$.
\begin{corollary} Let $[-\xi_+,-\xi_-]$ be the Sacker and Sell spectrum of $a_2$. Then,
\begin{enumerate}[label=\rm{(\roman*)}]
\item if $0\in[\xi_-,\xi_+]$, then \eqref{eq:laultima} exhibits the global pitchfork bifurcation diagram of minimal sets described by Theorem~$\mathrm{\ref{th:pitchfork}(iii)}$.
\item If $0\not\in[\xi_-,\xi_+]$, then \eqref{eq:laultima} exhibits the transcritical and saddle-node bifurcation diagram of minimal sets described by Theorem~$\mathrm{\ref{th:transcritical}}$.
\end{enumerate}
\end{corollary}
\begin{proof} The Sacker and Sell spectrum of $-3x^2+2a_2(\omega)x$ on the minimal set $\Omega\times{0}$ reduces to 0. That is, with the notation of Theorems~\ref{th:pitchfork} and \ref{th:transcritical}, $\lambda_-=\lambda_+=0$. As explained after Theorem~\ref{th:transcritical}, in this case, the bifurcation diagrams of Theorems~\ref{th:pitchfork}(iii) and \ref{th:transcritical} are the only possible ones and we can distinguish between them by the number of $\hat\tau_0$-minimal sets: one in the pitchfork case and two in the transcritical and saddle-node case. Note that $\eqref{eq:laultima}_0$ corresponds to $\eqref{eq:simplexample}_0$. Proposition~\ref{prop:newdiagram} ensures that $\eqref{eq:simplexample}_0$, and therefore $\eqref{eq:laultima}_0$, admits $\hat{\mathcal{M}}_0$ as unique minimal set if and only if $0\in[\xi_-,\xi_+]$, from where the statements follow.
\end{proof}
\subsection{Bifurcation of a recurrent solution}\label{subsec:porponer}
We complete the paper by translating part of the information that we have obtained to the problem we posed at the beginning of Subsection \ref{subsec:biflambdax}. Let us assume that the function $f$ of \eqref{eq:62primera} is $C^{0,2}(\Omega\times\mathbb{R},\mathbb{R})$, (Co) and (SDC). Recall that the minimal set $\mathcal M$ for the flow $\tau$ defined by \eqref{eq:62flujo} is given by the hull of the graph of the recurrent solution $\widetilde x(t)=u(t,\omega_0,\widetilde x(0))$ of \eqref{eq:62primera}, and note that $\widetilde x(t)$ solves the equation \eqref{eq:62}$_\lambda$ corresponding to $(\omega_0,\widetilde x(0))\in\mathcal M$ for all $\lambda\in\mathbb{R}$. It is easy to check that the function $g$ of $\eqref{eq:new63}_\lambda$ is $C^{0,2}(\mathcal M\times\mathbb{R},\mathbb{R})$, (Co) and (SDC), and the Sacker and Sell spectrum of $g_x$ on the corresponding minimal set $\mathcal M\times\{0\}$ coincides with that of $f_x$ on $\mathcal M$, which we represent by $[-\lambda_+,-\lambda_-]$.

To begin with, we assume that the bifurcation diagram for $\eqref{eq:new63}_\lambda$ around $\mathcal M\times \{0\}$ is that of a generalized transcritical bifurcation described by Theorem~\ref{th:transcritical}. Then, $\widetilde x(t)$ is hyperbolic attractive for \eqref{eq:62}$_\lambda$ if $\lambda<\lambda_-$, nonhyperbolic if $\lambda\in[\lambda_-,\lambda_+]$ and hyperbolic repelling for $\lambda>\lambda_+$; and, as $\lambda$ crosses $\lambda_+$ from the left (resp.~$\lambda_-$ from the right), the solution $\widetilde x(t)$ ``bifurcates" in two solutions of \eqref{eq:62}$_{\lambda}$ corresponding to $(\omega_0,\widetilde x(0))$: $\widetilde x(t)$ and a new hyperbolic attractive (resp.~repelling) recurrent solution $\widetilde x_\lambda(t)$, which persists for all $\lambda>\lambda_+$ (resp.~in a bounded interval $(\lambda_0,\lambda_-)$). In general, we cannot assure that $\lim_{\lambda\to\lambda_\pm}\widetilde x_\lambda(t)=\widetilde x(t)$ for all $t\in\mathbb R$, but there is a certain type of approaching which we describe in what follows for $\lambda_+$. Let $\widetilde{\mathcal M}_\lambda=\{(\omega,z,\tilde\kappa_\lambda(\omega,z))\,| \;(\omega,z)\in\mathcal M\}$ be the hull of the graph of $\widetilde x_\lambda$ for $\lambda>\lambda_+$. Then, there exists the pointwise limit $\tilde\kappa_{\lambda_+}(\omega,z)=\lim_{\lambda\downarrow\lambda_+}\tilde\kappa_\lambda(\omega,z)$, and $\tilde\kappa_\lambda$ is one of the two delimiter equilibria of the global attractor for all $\lambda\ge\lambda_+$; and hence, due to the collision properties explained in Theorem \ref{th:transcritical}, $\tilde\kappa_{\lambda_+}(\omega,z)=z$ for all $(\omega,z)$ in a $\tau$-invariant residual subset $\mathcal R$ of $\mathcal M$. This means that
$\lim_{\lambda\downarrow\lambda_+}\widetilde x_\lambda(t)=\widetilde x(t)$ for all $t\in\mathbb R$ in the case that $(\omega_0,\widetilde x(0))\in\mathcal R$, which we cannot a priori know. What we can assert in all the cases is that $\widetilde x_{\lambda_+}(t)=\lim_{\lambda\downarrow\lambda_+}\widetilde x_\lambda(t)$ is a bounded globally defined solution of the equation \eqref{eq:62}$_{\lambda_+}$ corresponding to $(\omega_0,\widetilde x(0))$, and that $\liminf_{t\to\pm\infty}|\widetilde x_{\lambda_+}(t)-\widetilde x(t)|=0$. Similar limiting properties (and restrictions) can be described for $\lambda\uparrow\lambda_-$.

If the bifurcation diagram is that of Theorem \ref{th:pitchfork}(iii), with a global pitchfork bifurcation at $\lambda_+$, then the asymptotic exponential stability of $\widetilde x(t)$ for $\lambda<\lambda_-$ is also lost as $\lambda$ crosses $\lambda_-$, and two hyperbolic attractive recurrent solutions $\widetilde x_\lambda^u(t)$ and $\widetilde x^l_\lambda(t)$ appear for $\lambda>\lambda_+$, with limiting approaching to $\widetilde x(t)$ as $\lambda\downarrow\lambda_+$ similar to that described above.

The possible dynamical complexity of these cases is due to the possible existence of a pinched compact invariant set for \eqref{eq:new63}$_{\lambda_+}$ which is not a copy of the base. Examples of these ``far away from trivial"~situations can be found in \cite{nunezobaya}, even in the simplest situation of uniquely ergodic flow on $\mathcal M$.


\begin{thebibliography}{99}
\bibitem{alob3} A.I.~Alonso, R.~Obaya. The structure of the bounded trajectories set of a scalar convex differential equation. \emph{Proc. Roy. Soc. Edinburgh A} \textbf{133} (2), 237-263 (2003).
\bibitem{anagjager1} V.~Anagnostopoulou, T.~J\"{a}ger. Nonautonomous saddle-node bifurcations: random and deterministic forcing. \emph{J. Differential Equations} \textbf{253} (2), 379–399 (2012).
\bibitem{anjk} V.~Anagnostopoulou, T.~J\"{a}ger, G.~Keller. A model for the nonautonomous Hopf bifurcation. \emph{Nonlinearity} \textbf{28} (7), 2587–2616 (2015).
\bibitem{arnol} L.~Arnold. \emph{Random dynamical systems.} Springer-Verlag, Berlin, 1998.
\bibitem{brbh} B.~Braaksma, H.~Broer, G.~Huitema. Toward a quasi-periodic bifurcation theory. \emph{Memoirs Amer. Math. Soc.} \textbf{83} (421), 83-175 (1990).
\bibitem{lineardissipativescalar} J.~Campos, C.~N\'{u}\~{n}ez, R.~Obaya. Uniform stability and chaotic dynamics in nonhomogeneous linear dissipative scalar ordinary differential equations. Preprint, 2021.
\bibitem{pullbackforwards} T.~Caraballo, J.~Langa, R.~Obaya. Pullback, forward and chaotic dynamics in 1D non-autonomous linear-dissipative equations. \emph{Nonlinearity}, \textbf{30} (1), 274-299 (2017).
\bibitem{carvalho1} A.~Carvalho, J.~Langa, J.~Robinson. \emph{Attractors for infinite-dimensional non-autonomous dynamical systems.} Applied Mathematical Sciences \textbf{182}, Springer, New York, 2013.
\bibitem{chueshov3} I.D.~Chueshov. \emph{Monotone random systems theory and applications.} Lecture Notes in Mathematics \textbf{1779}. Springer-Verlag, Berlin, 2002.
\bibitem{coppel1} W.A.~Coppel. \emph{Dichotomies in stability theory}. Lecture Notes in Mathematics \textbf{629}, Springer-Verlag, Berlin-New York, 1978.
\bibitem{sinai1} I.P.~Cornfeld, S.V.~Fomin, Ya.G.~Sinai. \emph{Ergodic theory}. Springer-Verlag, New York, 1982.
\bibitem{ellis1} R.~Ellis. \emph{Lectures on topological dynamics}. W.A. Benjamin, New York, 1969.
\bibitem{fabri1} R.~Fabbri, R.~Johnson. A nonautonomous saddle-node bifurcation pattern. \emph{Stoch. Dyn.} \textbf{4} (3), 335-350 (2004).
\bibitem{fuhrmann1} G.~Fuhrmann. Non-smooth saddle-node bifurcations III: Strange attractors in continuous time. \emph{J. Differential Equations} \textbf{261} (3), 2109-2140 (2016).
\bibitem{furstenberg1} H.~Furstenberg. Strict ergodicity and transformation of the torus. \emph{Amer. J. Math.} \textbf{83}, 573-601 (1961).
\bibitem{jager1} T.H.~J\"{a}ger. Quasiperiodically forced interval maps with negative Schwarzian derivative. \emph{Nonlinearity} \textbf{16} (4), 1239-1255 (2003).
\bibitem{jager2} T.H.~J\"{a}ger. The creation of strange non-chaotic attractors in non-smooth saddle-node bifurcations. \emph{Mem. Amer. Math. Soc.} \textbf{201} (945), 2009.
\bibitem{joma} R.~Johnson, F.~Mantellini. A nonautonomous transcritical bifurcation problem with an application to quasi-periodic bubbles. \emph{Discrete Contin. Dyn. Syst.} \textbf{9} (1), 209-224 (2003).
\bibitem{johnson1} R.~Johnson, R.~Obaya, S.~Novo, C.~N\'{u}\~{n}ez, R.~Fabbri. \emph{Nonautonomous linear Hamiltonian systems: oscillation, spectral theory and control}. Developments in Mathematics \textbf{36}, Springer, 2016.
\bibitem{johnson2} R.~Johnson, K.~Palmer, G.R.~Sell. Ergodic properties of linear dynamical systems. \emph{SIAM J. Math. Anal.} \textbf{18} (1), 1-33 (1987).
\bibitem{kloeden1} P.~Kloeden. Pitchfork and transcritical bifurcations in systems with homogeneous nonlinearities and an almost periodic time coefficient. \emph{Commun. Pure Appl. Anal.} \textbf{3} (2), 161-173 (2004).
\bibitem{kloeden2} P.~Kloeden, M.~Rassmussen. Nonautonomous dynamical systems. \emph{Mathematical Surveys and Monographs} \textbf{176}, American Mathematical Society, Providence, 2011.
\bibitem{langa1} J.A.~Langa, J.C.~Robinson, A.~Su\'{a}rez. Bifurcations in non-autonomous scalar equations. \emph{J. Differential Equations} \textbf{221} (1), 1-35 (2006).
\bibitem{mane1} R.~Ma\~{n}\'{e}. \emph{Ergodic theory and differentiable dynamics}. Springer-Verlag, Berlin, 1987.
\bibitem{demelo} W.~de~Melo, S.~van~Strien. \emph{One-dimensional dynamics.} Springer-Verlag, Berlin, 1993.
\bibitem{nemytskii1} V.V.~Nemytskii, V.V.~Stepanov. \emph{Qualitative theory of differential equations.} Princeton Mathematical Series \textbf{22}, Princeton University Press, 1960.
\bibitem{almostautomorphic} S.~Novo, C.~N\'{u}\~{n}ez, R.~Obaya. Almost automorphic and almost periodic dynamics for quasimonotone non-autonomous functional differential equations. \emph{J. Dynam. Differential Equations} \textbf{17} (3), 589-619 (2005).
\bibitem{nunezobaya} C.~N\'{u}\~{n}ez, R.~Obaya. A non-autonomous bifurcation theory for deterministic scalar differential equations. \emph{Discrete Contin. Dyn. Syst. Ser. B} \textbf{9} (3-4), 701-730 (2008).
\bibitem{nuob10} C.~N\'{u}\~{n}ez, R.~Obaya. Li-Yorke chaos in nonautonomous Hopf bifurcation patterns-I. \emph{Nonlinearity} \textbf{32} (10), 3940-3980 (2019).
\bibitem{nuos4} C.~N\'{u}\~{n}ez, R.~Obaya, A.M.~Sanz. Minimal sets in monotone and concave skew-product semiflows I: a general theory. \emph{J. Differential Equations} \textbf{252} (10), 5492-5517 (2012).
\bibitem{potz} C.~P\"{o}tzsche. Nonautonomous continuation of bounded solutions. \emph{Commun. Pure Appl. Anal.} \textbf{10} (3), 937-961 (2011).
\bibitem{potz2} C.~P\"{o}tzsche. Nonautonomous bifurcation of bounded solutions I: a Lyapunov-Schmidt approach. \emph{Discrete Contin. Dyn. Syst. Ser. B} \textbf{14} (2), 739-776 (2010).
\bibitem{potz4} C.~P\"{o}tzsche. Nonautonomous bifurcation of bounded solutions II: A shovel-bifurcation pattern. \emph{Discrete Contin. Dyn. Syst.} \textbf{31} (3), 941-973 (2011).
\bibitem{rasmussen2} M.~Rasmussen. \emph{Attractivity and bifurcation for nonautonomous dynamical systems}. Lecture Notes in Mathematics \textbf{1907}. Springer, Berlin, 2007.
\bibitem{rasmussen1} M.~Rasmussen. Nonautonomous bifurcation patterns for one-dimensional differential equations. \emph{J. Differential Equations} \textbf{234} (1), 267-288 (2007).
\bibitem{sackersell3} R.J.~Sacker, G.R.~Sell. A spectral theory for linear differential systems. \emph{J. Differential Equations} \textbf{27}, 320-358 (1978).
\bibitem{tineo1} A.~Tineo. A result of Ambrosetti-Prodi type for first-order ODEs with cubic non-linearities. I, II. \emph{Ann. Mat. Pura Appl. (4)} \textbf{182} (2), 113-128, 129-141 (2003).
\bibitem{walters1} P.~Walters. \textit{An introduction to ergodic theory}. Graduate Texts in Mathematics \textbf{79}, Springer-Verlag, New York-Berlin, 1982.
\end{thebibliography}
\end{document}